\documentclass{article}

\usepackage{arxiv}

\usepackage[utf8]{inputenc} 
\usepackage[T1]{fontenc}    
\usepackage{hyperref}       
\usepackage{url}            
\usepackage{booktabs}       
\usepackage{amsfonts}       
\usepackage{nicefrac}       
\usepackage{microtype}      
\usepackage{lipsum}

\usepackage{enumerate}
\usepackage{amsmath}
\usepackage{amssymb}
\usepackage{graphicx}

\newtheorem{theorem}{Theorem}

\newtheorem{condition}[theorem]{Condition}

\newtheorem{definition}[theorem]{Definition}

\newtheorem{lemma}[theorem]{Lemma}

\newtheorem{proposition}[theorem]{Proposition}
\newtheorem{remark}[theorem]{Remark}

\newenvironment{proof}[1][Proof]{\noindent\textbf{#1.} }{\ \rule{0.5em}{0.5em}}

\allowdisplaybreaks[4]
\newcommand{\h}{\dot{H}^{\gamma}(\mathbb{R}^d)}
\newcommand{\hh}{\mathbb{H}}
\newcommand{\vv}{\mathbb{V}}
\newcommand{\e}{\mathbb{E}}
\newcommand{\be}{\bar{\mathbb{E}}}
\newcommand{\um}{u_n^m}
\newcommand{\x}{\chi_{n_1,n_2}^m}
\newcommand{\p}{\mathbf{B}^{m_1,m_2}}
\newcommand{\ur}{u^m(r,\tau-t,u_0)}
\newcommand{\us}{u^m(s,\tau-t,u_0)}
\newcommand{\ue}{u^{\varepsilon}}
\newcommand{\tu}{\tilde{u}^\varepsilon}

\newcommand{\ty}{\tilde{Y}^\varepsilon}
\newcommand{\tr}{\tilde{R}^\varepsilon}
\numberwithin{equation}{section}
\numberwithin{theorem}{section}

\def\({\left(}
\def\){\right)}

\ifpdf
  \DeclareGraphicsExtensions{.eps,.pdf,.png,.jpg}
\else
  \DeclareGraphicsExtensions{.eps}
\fi

\title{Dynamics and large deviations for  fractional stochastic partial differential  equations with L\'evy noise}

\author{Jiaohui Xu\thanks{Center
for Nonlinear Studies, School of Mathematics, Northwest University, Xi'an 710127, P. R. China  (E.mail: jxu@us.es).}
\and Tom\'as Caraballo\thanks{Dpto. Ecuaciones Diferenciales y An\'alisis Num\'erico, Universidad de Sevilla, Facultad de Matem\'aticas, c/ Tarfia s/n, 41012-Sevilla, Spain (E.mail: caraball@us.es).}
  \and Jos\'e Valero\thanks{Centro de Investigaci\'on Operativa,
Universidad Miguel Hern\'andez de Elche, Avenida de la Universidad s/n, 03202-Elche, Spain (E.mail: jvalero@umh.es).}
}

\usepackage{amsopn}

\usepackage{amsopn}

\ifpdf
\hypersetup{
  pdftitle={Dynamics and large deviations for  fractional stochastic partial differential  equations with L\'evy noise},
  pdfauthor={J.H. Xu, T. Caraballo and J. Valero}
}
\fi

\begin{document}

\maketitle

\begin{abstract} This paper is mainly concerned with a kind of fractional stochastic evolution equations driven by L\'evy noise in a bounded domain.
We first state the well-posedness of the problem via iterative approximations and energy estimates. Then, the existence and uniqueness  of weak pullback mean random attractors for the equations {are} established by defining a mean random dynamical system.  Next, we prove the existence of invariant measures  when the problem is autonomous by means of the fact that $H^\gamma(\mathcal{O})$ is compactly embedded in $L^2(\mathcal{O})$ with $\gamma\in (0,1)$. Moreover, the uniqueness of this invariant measure is presented which ensures the ergodicity of the problem. Finally,  a large deviation principle result for solutions of SPDEs perturbed by small  L\'evy noise and Brownian motion is obtained by a variational formula for positive functionals of a Poisson random measure and Brownian motion. 
Additionally, the results are  illustrated by the fractional stochastic Chafee-Infante equations.
\end{abstract}



\section{Introduction}

In this paper, we consider the following   fractional stochastic PDEs driven by L\'evy noise and Brownian motion,
\begin{alignat}{3}\label{eq1-1}
\begin{cases}
\displaystyle d u(t)+(-\Delta)^{\gamma}u(t)dt+f(u(t))dt=g(t,u(t))dW(t)+\int_{E}h(u(t-),\xi)\tilde{N}(dt,d\xi),\\
u(t,x)=0,\\
u(\tau,x)=u_0(x),\\
\end{cases}
\begin{aligned}
&\hskip-.3cm\mbox{in}~\mathcal{O}\times (\tau,\infty),\\
&\hskip-.3cm\mbox{on}~\partial\mathcal{O}\times(\tau,\infty),\\
&\hskip-.3cm\mbox{in}~\mathcal{O},\\
\end{aligned}
\end{alignat}
where $\mathcal{O}\subset \mathbb{R}^d$ $(d>1)$ is a bounded domain with smooth boundary, $\tau\in\mathbb{R}$,  the operator $(-\Delta)^{\gamma}$ with $\gamma\in(0,1)$ is the so-called  fractional Laplacian, $f:\mathbb{R}\rightarrow \mathbb{R}$ is a polynomial of odd degree with positive leading coefficient, the functions $g(t,u)$ and $h(u,\xi)$ satisfy some conditions which will be specified later. We consider problem \eqref{eq1-1} with respect to a given stochastic basis $(\Omega,\mathcal{F},\{\mathcal{F}_t\}_{t\in\mathbb{R}},\mathbb{P},W,N)$ and a Hilbert space $U$, where
 $W$ is  a  two-sided $U$-valued cylindrical Wiener process and $N$ is a Poisson measure  induced by a stationary $\mathcal{F}_t$-Poisson point process on $(\tau,T]\times E$ with a $\sigma$-finite intensity measure $L_{T-\tau}\otimes\lambda$, $L_{T-\tau}$ is the Lebesgue measure on $(\tau,T]$ and $\lambda$ is a $\sigma$-finite measure  on a measurable space $E$, $\tilde{N}(dt,d\xi):=N(dt,d\xi)-\lambda(d\xi)dt$ is the compensated Poisson random measure. Assume $W$ and $\tilde{N}$ are independent.
 
 Stochastic partial differential equations arise in many different fields since stochastic perturbations originated from many natural sources cannot be ignored in a realistic modeling. In recent decades,  stochastic PDEs driven by Brownian motion have been extensively studied theoretically \cite{L3, X3}, concerning well-posedness, existence of stationary solutions, stochastic attractors and invariant measures.  However, the fact that forcing terms may be treated stochastically does not mean that details of the stochastic treatment are arbitrary \cite{P1}. In fact, it  turns out that a process may not only be Gaussian but also can exhibit skew, fat tails and other properties usually associated with more exotic types of 
 stochastic phenomena, such as non-Gaussian L\'evy noise. For example, they have been used to develop models for neuronal activity that account for synaptic impulses occurring randomly, both in time and at different locations of a spatially extended neuron. Other applications arise in chemical reaction-diffusion systems and stochastic turbulence models \cite{G3, M3, Z2}.

 The fractional Laplacian operator, which is written as $(-\Delta)^\gamma$ with $\gamma\in(0,1)$, has multiple equivalent characterizations \cite{R3, S1}. In the present paper, we will mainly adopt the non-local {one} (see \eqref{eq2-9}).  Although the eigenfunctions of $(-\Delta)^\gamma$ are not smooth {in the sense that they are just H\"older continuous up to the boundary of $\mathcal{O}$ but not Lipschitz continuous,} {it is possible to  construct a continuous operator $A$ which involves $(-\Delta)^\gamma$ (see \eqref{eq2-12}).}
{By means of the fact that $W^{\gamma,2}(\mathcal{O})$ is compactly embedded in $L^2(\mathcal{O})$ and the Hilbert-Schmidt theorem, we can find the eigenfunctions $e_j\ (j\in\mathbb{N})$ of $A$} which forms an orthonormal basis of {$L^2(\mathcal{O})$} with corresponding eigenvalues $0<\lambda_1\leq \lambda_2\leq \cdots \rightarrow \infty$ $({Ae_j=\lambda_je_j})$. Moreover, the domain of $A^r$ is denoted by $D(A^r)$ which is equipped with the norm $\|u\|_{D(A^r)}=\|A^ru\|_{L^2(\mathcal{O})}$ for $u\in D(A^r)$. Notice that, $\{e_j/\lambda_j^r\}$ is a complete orthonormal system of $D(A^r)$. By the Riesz representation theorem, $D(A^{-r})$ is the dual space of $D(A^r)$. In this way, we know that $D(A^r)$ is continuously embedded into $L^p(\mathcal{O})$ as long as $r$ is large enough \cite{W2}.
 
The non-local or memory effects are ubiquitous in physics and engineering \cite{B6, T2, T3}. Therefore, evolutionary equations with fractional Laplacian operator can be used to model these non-local effects (see \cite{A2, CPY, C2,  D4, G1, WR1, WR2, XC} and the references therein).  Particularly, the solutions and their dynamics of fractional partial differential equations have been extensively studied by a great amount of researchers, see \cite{G3, G2,  L31, S2, W2} and the references therein.

Consequently, it is meaningful to study the dynamics of problem \eqref{eq1-1}.  To be precise, the first goal of this paper is to analyze the well-posedness of \eqref{eq1-1} in 
 $ L^2(\Omega; \mathbf{D}([\tau,T];L^2(\mathbb{R}^d)))\cap L^2(\Omega; L^2(\tau,T;W^{\gamma,2}(\mathbb{R}^d)))\cap L^{p+1}(\Omega;L^{p+1}(\tau,T;L^{p+1}(\mathbb{R}^d)))$, the existence and uniqueness of weak pullback mean random attractors for the mean random dynamical systems generated by the solution operators. The second goal is to prove  the existence of invariant measures and ergodicity to problem \eqref{eq1-1} in the autonomous case. This result holds true {since $W^{\gamma,2}(\mathbb{R}^d)$ is compactly embedded in $L^2(\mathcal{O})$ }, where $\mathcal{O}\subset \mathbb{R}^d $ is a bounded domain.
 
The third goal, which is also the main novelty of this paper, is to establish a large deviation principle to fractional stochastic PDEs \eqref{eq1-1} with L\'evy noise by a variational representation obtained in \cite{B3} and weak convergence approach. 
 The large deviation principle is an active and important topic in probability and statistics. Large deviation properties of SPDEs driven by infinite dimensional Brownian motions and Poisson random measure have been studied in \cite{ B7, B8, B3}. However, as far as the authors are aware, there are no any results about the large deviation principle to FSPDEs, and our work will fill  this gap.  To this end, we follow some ideas introduced by \cite{B3} which can be properly  adapted to our problem. This is mainly due to the fact that the eigenfunctions of the fractional Laplacian operator $(-\Delta)^\gamma$ share similar properties to the ones of  the classical  Laplacian operator $-\Delta$. By carrying out a careful analysis, we need to impose some assumptions on $\gamma$,  $d$  and $p$ (namely, $p+1\in (2,\frac{2d}{d-2\gamma}]$) such that $W^{\gamma,2}(\mathbb{R}^d)$ is continuously embedded in $L^{p+1}(\mathbb{R}^d)$, which allows us to accomplish the proposed study.

The paper is organized as follows. In Section \ref{s2}, we review the definition of fractional Laplacian operator, impose the conditions on the nonlinear terms and introduce the concept of a large deviation principle. Then, the well-posedness of problem
 \eqref{eq1-1} is established in Section \ref{s3}  by an iterative method. Section \ref{s4} is devoted to the existence and uniqueness of weak pullback mean random attractors. In Section \ref{s5}, we study the existence of invariant {measures} and ergodicity to problem \eqref{eq1-1} when it is autonomous. In Section \ref{s6}, a general large deviation result to \eqref{eq1-1} is proved by a variational formula for positive functionals of a Poisson random measure and Brownian motion. An illustrative example concerning Chafee-Infante model is exhibited in Section \ref{s7} and an appendix with the proofs of some results concludes our paper in Section \ref{s8}.
\section{Preliminaries}\label{s2}
In this section, we will introduce some basic definitions and properties of the fractional Laplacian operator,  impose proper assumptions  on nonlinear terms in \eqref{eq1-1} and recall the general criteria for a large deviation principle. 

\subsection{Fractional setting}

Let $\mathcal{S}$ be the Schwartz space of rapidly decaying $C^{\infty}$ functions on $\mathbb{R}^d$. Then the integral fractional Laplacian operator $(-\Delta)^{\gamma}$ with $0<\gamma<1$ is defined, for $u\in\mathcal{S}$, by
\begin{equation}\label{eq2-9}
(-\Delta)^{\gamma}u(x)=-\frac{1}{2}C(d,\gamma)\int_{\mathbb{R}^d}\frac{u(x+y)+u(x-y)-2u(x)}{|y|^{d+2\gamma}}dy,\qquad x\in\mathbb{R}^d,
\end{equation}
where $C(d,\gamma)$ is a positive constant given by
\begin{equation}\label{eq2-10}
C(d,\gamma)=\frac{\gamma4^{\gamma}\Gamma(\frac{d+2\gamma}{2})}{\pi^{\frac{d}{2}}\Gamma(1-\gamma)}.
\end{equation}
The reader is referred to \cite{N1} for more details on the integral fractional operators.  
Moreover, for any real $0<\gamma<1$, the fractional Sobolev space $W^{\gamma,2}(\mathbb{R}^d):=H^{\gamma}(\mathbb{R}^d)$ is defined by
$$H^{\gamma}(\mathbb{R}^d)=\bigg\{u\in L^2(\mathbb{R}^d):\int_{\mathbb{R}^d}\int_{\mathbb{R}^d}\frac{|u(x)-u(y)|^2}{|x-y|^{d+2\gamma}}dxdy<\infty\bigg\},
$$
endowed with the norm
$$\|u\|_{H^{\gamma}(\mathbb{R}^d)}={\bigg(}\int_{\mathbb{R}^d}|u(x)|^2dx+\int_{\mathbb{R}^d}\int_{\mathbb{R}^d}\frac{|u(x)-u(y)|^2}{|x-y|^{d+2\gamma}}dxdy\bigg)^{\frac{1}{2}}. $$
We  denote the Gagliardo semi-norm of $H^\gamma(\mathbb{R}^d)$ as $\|\cdot\|_{\dot{H}^\gamma(\mathbb{R}^d)}$, i.e.,
\begin{equation*}
\|u\|^2_{\dot{H}^\gamma(\mathbb{R}^d)}=\int_{\mathbb{R}^d}\int_{\mathbb{R}^d}\frac{|u(x)-u(y)|^2}{|x-y|^{d+2\gamma}}dxdy,\quad u\in H^{\gamma}(\mathbb{R}^d).
\end{equation*}
Then, for all $u\in H^{\gamma}(\mathbb{R}^d)$, we have $\|u\|^2_{H^{\gamma}(\mathbb{R}^d)}=\|u\|^2+\|u\|^2_{\h}$. Note that $H^{\gamma}(\mathbb{R}^d)$ is a Hilbert space with inner product 
$$(u,v)_{H^{\gamma}(\mathbb{R}^d)}=\int_{\mathbb{R}^d}u(x)v(x)dx+\int_{\mathbb{R}^d}\int_{\mathbb{R}^d}\frac{(u(x)-u(y))(v(x)-v(y))}{|x-y|^{d+2\gamma}}dxdy,\qquad \forall u,v\in H^{\gamma}(\mathbb{R}^d).$$
By \cite{N1}, we infer that for every fixed $\gamma\in(0,1)$ and 
$u\in H^{\gamma}(\mathbb{R}^d)$, the norm $\|u\|_{H^{\gamma}(\mathbb{R}^d)}$ is equivalent to $\left(\|u\|^2_{L^2(\mathbb{R}^d)}+{\|(-\Delta)^{\frac{\gamma}{2}}u\|_{L^2(\mathbb{R}^d)}^2}\right)^{\frac{1}{2}}$. More precisely, we have 
\begin{equation*}
\|u\|^2_{H^{\gamma}(\mathbb{R}^d)}=\|u\|^2_{L^2(\mathbb{R}^d)}+\frac{2}{C(d,\gamma)}\|(-\Delta)^{\frac{\gamma}{2}}u\|^2_{L^2(\mathbb{R}^d)},\quad \forall u\in H^{\gamma}(\mathbb{R}^d).
\end{equation*}
Since the fractional Laplacian operator $(-\Delta)^{\gamma}$ defined above is a non-local one, we here interpret the homogenous Dirichlet boundary as $u=0$ on $\mathbb{R}^d\backslash\mathcal{O}$ instead of $u=0$ only on $\partial \mathcal{O}$. Such an interpretation is consistent with the non-local nature of the integral fractional Laplacian and has been used in many publications. Based on this interpretation,
we recall $\hh$ and $\vv$ are  two Hilbert spaces given by 
$\hh=\{u\in L^2(\mathbb{R}^d): u=0~\mbox{a.e. in}~\mathbb{R}^d \backslash\mathcal{O}\}$ and 
$\vv=\{u\in H^{\gamma}(\mathbb{R}^d):u=0~\mbox{a.e. in}~\mathbb{R}^d \backslash\mathcal{O}\}$,  respectively. Then we have $\vv\hookrightarrow \hh=\hh^*\hookrightarrow \vv^*$, where $\hh^*$ is identified with $\hh$ by the Riesz representation theorem, $\vv^*$  and $\hh^*$ are the dual spaces of  $\vv$ and $\hh$, respectively.

We conclude this subsection  by introducing some notation. For $p\geq 1$, we denote by $L^p(\mathbb{R}^d)$ the usual $L^p$-space over $\mathbb{R}^d$ with the standard norm $\|\cdot\|_p$. The norm and inner product of $\hh$ and $\vv$ are denoted by $|\cdot|$ and $\|\cdot\|$, $(\cdot,\cdot)$ and $((\cdot,\cdot))$, respectively. Moreover, the norm of $\vv^*$ is denoted by $\|\cdot\|_{*}$. 
For the simplicity of notation, when no confusion may arise, we will use the unified notation $<\cdot,\cdot>$ to denote the dual relations between different spaces. For a Polish space $E$, denote by {$C([0,T];E)$} and $\mathbf{D}([0,T];E)$ the spaces of continuous functions and right continuous functions with left limits from $[0,T]$ to $E$,  respectively, endowed with the uniform topology both, if not specified.

\subsection{Stochastic setting and assumptions}

 Let $(\Omega,\mathcal{F},\{\mathcal{F}_t\}_{t\in\mathbb{R}},\mathbb{P})$ be a  filtered probability space satisfying  the usual conditions, i.e., 
 $\{\mathcal{F}_t\}_{t\in\mathbb{R}}$ is an increasing right continuous 
 family of sub-$\sigma$-algebras of $\mathcal{F}$ that contains all $\mathbb{P}$-null sets.
 The collection of all strongly-measurable, square-integrable $\hh$-valued random variables, denoted by $L^2(\Omega;\hh)$, is a Banach space equipped with the norm $\|u(\cdot)\|^2_{L^2(\Omega;\hh)}=\mathbb{E}|u(\cdot,\omega)|^2$, where the expectation $\mathbb{E}$ is defined by $\mathbb{E}u=\int_{\Omega}u(\omega)d\mathbb{P}$.
  Furthermore,  let $H$ and $U$ be two separable Hilbert spaces, $\mathcal{L}_2(U;H)$ denote the space of Hilbert-Schmidt operators from a separable Hilbert space $U$ to $H$ with norm $\|\cdot\|_{\mathcal{L}_2(U;H)}$ (see \cite{D1} for more details). 


Throughout this paper, we impose the following conditions on $f$, $g$ and $h$.\\
$\bullet$ {\bf Assumptions on nonlinear term $f$.}
 Suppose the nonlinear term $f:\mathbb{R}\rightarrow \mathbb{R}$ has the following form,
\begin{equation}\label{eq2-1}
f(u)=\sum_{k=1}^{2p}f_{2p-k}u^{k-1},\qquad f_0>0,~~p\in\mathbb{N}.
\end{equation}
In fact, no significant changes in the proofs of the results presented here are required if we consider, more generally, a continuously differentiable function $f$ on $\mathbb{R}$ satisfying
\begin{enumerate}
\item[$(f_1)$] $|f(u)|\leq l_1(1+|u|^{p})$,\qquad $(f_2)$ $u\cdot f(u)\geq -l_2+l_3|u|^{p+1}$,
\qquad $(f_3)$ $f^{\prime}(u)\geq -l_4$,
\end{enumerate}
for some  $l_j >0$, $j=1,2,3,4.$

For convenience, we fix a positive number $\delta $ and set 
\begin{equation}\label{eq2-2}
F(u)=f(u)-\delta u,\qquad \forall u\in\mathbb{R}.
\end{equation}
By \eqref{eq2-2} and conditions $(f_1)$-$(f_3)$,  after simple calculations, we find that
\begin{enumerate}
\item[$(F_1)$] $|F(u)|\leq k_1(1+|u|^{p})$,\qquad $(F_2)$ $u\cdot F(u)\geq -k_2+k_3|u|^{p+1}$,
\qquad $(F_3)$ $F^{\prime}(u)\geq -k_4:=-l_4-\delta$,
\end{enumerate}
for some  $k_j> 0$, $j=1,2,3,4$.\\
$\bullet$ {\bf Assumptions on nonlinear term $g$.} Suppose {$g:[\tau,\infty)\times \hh\rightarrow \mathcal{L}_2(U;\hh)$} is locally Lipschitz continuous   and grows linearly in its second argument uniformly for $t\in[\tau,\infty)$, that is:
\begin{enumerate}
\item[$(g_1)$] For every $r>0$, there exists a positive constant $L_g(r)$ depending on $r$, such that for all $t\in[\tau,\infty)$, $u_1,u_2\in \hh$ with $|u_1|\leq r$ and $|u_2|\leq r$,
\begin{equation}\label{eq2-3}
\|g(t,u_1)-g(t,u_2)\|^2_{\mathcal{L}_2(U;\hh)}\leq L_g(r)|u_1-u_2|^2;
\end{equation}
\item [$(g_2)$]
There exists a positive constant $C_g$, such that for all $t\in[\tau,\infty)$ and $u\in\hh$,
\begin{equation}\label{eq2-4}
\|g(t,u)\|^2_{\mathcal{L}_2(U;\hh)}\leq C_g(1+|u|^2);
\end{equation}
\item[$(g_3)$] For every fixed $u\in \hh$, $g(\cdot,u):[\tau,\infty)\rightarrow \mathcal{L}_2(U;\hh)$ is progressively measurable.
\end{enumerate}
$\bullet$ {\bf Assumptions on nonlinear term $h$.} Suppose $ h: \hh\times E\rightarrow \hh$  is locally Lipschitz continuous and grows linearly in its first argument uniformly for $\xi\in E^\prime$,  where $E^\prime\subset E$ satisfying $\lambda(E^\prime)<\infty$, precisely:
\begin{enumerate}
\item[$(h_1)$] For every $r>0$, there exists a positive constant $L_h(r)$ depending on $r$, such that, for all $u_1,u_2\in \hh$ with $|u_1|\leq r$ and $|u_2|\leq r$,
\begin{equation}\label{eq2-5}
\int_{E^\prime}|h(u_1,\xi)-h(u_2,\xi)|^2\lambda(d\xi)\leq L_h(r)|u_1-u_2|^2;
\end{equation}
\item [$(h_2)$]
There exists a positive constant $C_h$ such that, for all  $u\in\hh$,
\begin{equation}\label{eq2-6}
\int_{E^\prime}|h(u,\xi)|^2\lambda(d\xi)\leq C_h(1+|u|^2);
\end{equation}
\item[$(h_3)$] $h: \hh\times E\rightarrow \hh$ is a measurable mapping.
\end{enumerate}

In light of \eqref{eq2-2}, problem \eqref{eq1-1} can be put into the form when the boundary condition is replaced by $u=0$ on $\mathbb{R}^d\backslash\mathcal{O}$:
\begin{equation}\label{eq2-13}
\begin{split}
du(t)&+(-\Delta)^{\gamma}u(t)dt+F(u(t))dt+\delta u(t)dt\\
&=g(t,u(t))dW(t)+\int_{E}h(u(t-),\xi)\tilde{N}(dt,d\xi),\quad x\in\mathcal{O},~t>\tau,\\
\end{split}
\end{equation}
with boundary and initial conditions,
\begin{equation}\label{eq2-14}
u(t,x)=0,\quad x\in\mathbb{R}^d\backslash \mathcal{O}, ~t>\tau, \qquad \mbox{and}\qquad  u(\tau,x)=u_0(x),\quad x\in\mathcal{O}.
\end{equation}

To prove the existence and uniqueness of weak solutions to problem \eqref{eq2-13}-\eqref{eq2-14},  we follow    the ideas of \cite{W2}. To this end, 
let $a:\vv\times \vv \rightarrow \mathbb{R}$ be a bilinear form given by
\begin{equation}\label{eq2-11}
a(v_1,v_2)=\delta(v_1,v_2)+\frac{1}{2}C(d,\gamma)\int_{\mathbb{R}^d}\int_{\mathbb{R}^d}\frac{(v_1(x)-v_1(y))(v_2(x)-v_2(y))}{|x-y|^{d+2\gamma}}dxdy, \quad \forall v_1, v_2\in \vv,
\end{equation}
where $C(d,\gamma)$ is the constant in \eqref{eq2-10} and $\delta$ the one in \eqref{eq2-2}. For convenience, we associate an operator {$A: \mathbb{V}\rightarrow \mathbb{V}^*$} with $a$ in the following way:
\begin{equation}\label{eq2-12}
<A(v_1),v_2>_{(\vv^*,\vv)}=a(v_1,v_2),\qquad\mbox{for all}~ v_1,v_2\in \vv,
\end{equation}
where $<\cdot,\cdot>_{(\vv^*,\vv)}$ is the duality paring of $\vv^*$ and $\vv$. 
It follows from \cite{W2} that the inverse operator $A^{-1}:\vv^*\subset \hh\rightarrow \vv\subset\hh$ is  symmetric and compact. Therefore, the Hilbert-Schmidt theorem shows that $A$ has a family of eigenfunctions  $\{e_j\}_{j=1}^{\infty}$ such that $\{e_j\}_{j=1}^{\infty}$ forms an orthonormal basis of $\hh$. Moreover, if $\lambda_j$ is the eigenvalue corresponding to $e_j$, i.e.,
$Ae_j=\lambda_je_j,\ j=1,2,\cdots,$
then $\lambda_j$ satisfies 
$ 0<\delta<\lambda_1\leq \lambda_2\leq\cdots\leq\lambda_j\rightarrow\infty~~\mbox{as}~~j\rightarrow\infty.$ 
Actually, $e_j$ $(j\in\mathbb{N})$ are eigenfunctions of the integral fractional Laplacian operator $(-\Delta)^{\gamma}$. We will consider the fractional power of the operator $A$. Given $u\in\hh$, we have 
$ u=\sum_{j=1}^\infty a_je_j $ with $a_j=\int_{\mathbb{R}^d}u(x)e_j(x)dx.$
Then, for $r>0$, define $A^ru=\sum_{j=1}^\infty a_j\lambda_j^re_j$ provided the series is convergent for $u$ in $\hh$. The domain of $A^r$ is denoted by $D(A^r)$ which is equipped with the norm {$\|u\|_{D(A^r)}=|A^ru|$} for $u\in D(A^r)$. By the Riesz representation theorem, $D(A^{-r})$ is the dual space of $D(A^r)$.

{\subsection{Large deviation principle}
For a topological space $\mathcal{E}$, denote the corresponding Borel $\sigma$-field by $\mathcal{B}(\mathcal{E})$.   For a measure $\lambda$ on $\mathcal{E}$ and a Hilbert space $\hh$, let $L^2(\mathcal{E},\lambda;\hh)$ denote the space of measurable functions $f$ from $\mathcal{E}$ to $\hh$ such that $\int_{\mathcal{E}}|f(u)|^2\lambda(du)<\infty$. For a function $x:[0,T]\rightarrow \mathcal{E}$, we use the notation $x(t)$ to denote the evaluation of $x$ at $t\in[0,T]$. A similar convention will be followed for stochastic processes. Eventually, we say a collection $\{u^{\varepsilon}\}$ of $\mathcal{E}$-valued random variables is tight if the distributions of $u^{\varepsilon}$ are tight in $\mathrm{P}(\mathcal{E})$ (the space of probability measures on $\mathcal{E}$).
A function $I:\mathcal{E}\rightarrow [0,\infty]$ is called a rate function on $\mathcal{E}$, if for each $M<\infty$, the level set $\{u\in\mathcal{E}: I(u)\leq M\}$ is a compact subset of $\mathcal{E}$. A sequence $\{u^{\varepsilon}\}$ of $\mathcal{E}$-valued random variables is said to satisfy the Laplace principle upper bound (respectively lower bound) on $\mathcal{E}$ with rate function $I$, {if} for each $h\in C_b(\mathcal{E})$ (the space of real continuous bounded functions), 
$$ \limsup_{\varepsilon\rightarrow 0}\varepsilon \log\e\left\{\exp\left[-\frac{1}{\varepsilon}h(u^{\varepsilon})\right]\right\}\leq -\inf_{u\in\mathcal{E}}\{h(u)+I(u)\},$$
and
$$ \liminf_{\varepsilon\rightarrow 0}\varepsilon \log\e\left\{\exp\left[-\frac{1}{\varepsilon}h(u^{\varepsilon})\right]\right\}\geq -\inf_{u\in\mathcal{E}}\{h(u)+I(u)\},$$
respectively. The Laplace principle is said to hold for $\{u^{\varepsilon}\}$ with rate function $I$ if both the Laplace upper and lower bounds hold. It is well known that when $\mathcal{E}$ is a Polish space,  the family $\{u^{\varepsilon}\}$ satisfies the Laplace principle upper (respectively lower) bound with a rate function $I$ on $\mathcal{E}$ if and only if $\{u^{\varepsilon}\}$ satisfies the large derivation upper (respectively lower) bound for all closed sets (respectively open sets) with the rate function $I$. For more details, see \cite{B2} and the references therein.}

\section{Well-posedness of problem \eqref{eq2-13}-\eqref{eq2-14}}\label{s3}

We will prove the existence and uniqueness of solutions to \eqref{eq2-13}-\eqref{eq2-14} in the following sense.

\begin{definition}\label{def3-1}
Let $u_0\in L^2(\Omega; \hh)$ be $\mathcal{F}_\tau$-measurable.    An $\hh$-valued $\mathcal{F}_t$-adapted c\`adl\`ag stochastic process $u$ is called a solution of \eqref{eq2-13}-\eqref{eq2-14} if 
\begin{equation}\label{eq3-1}
u\in L^2(\Omega; \mathbf{D}([\tau,T];\hh))\cap L^2(\Omega; L^2(\tau,T;\vv))\cap L^{p+1}(\Omega;L^{p+1}(\tau,T;L^{p+1}(\mathbb{R}^d))),
\end{equation}
and, for all $t\geq \tau$, 
\begin{equation}\label{eq3-2}
\begin{split}
u(t)&+\int_{\tau}^t(-\Delta)^{\gamma}u(s)ds+\delta\int_\tau^tu(s)ds+\int_\tau^tF(u(s))ds\\
&~~=u_0+\int_\tau^tg(s,u(s))dW(s)+\int_\tau^t\int_Eh(u(s-),\xi)\tilde{N}(ds,d\xi), 
\end{split}
\end{equation}
in $(\vv\cap L^{p+1}(\mathbb{R}^d))^*$, $\mathbb{P}$-almost surely.
\end{definition}

\begin{theorem}\label{thm3-2}
Assume conditions $(F_1)$-$(F_3)$, $(g_1)$-$(g_3)$ and $(h_1)$-$(h_3)$ hold. Then, for every $\mathcal{F}_\tau$-measurable initial value $u_0\in L^2(\Omega;\hh)$, problem \eqref{eq2-13}-\eqref{eq2-14} has a unique global solution in the sense of Definition \ref{def3-1}.
Moreover, the solution $u$ depends continuously on $u_0$ from $L^2(\Omega;\hh)$ to $L^2(\Omega;\mathbf{D}([\tau,T];\hh))\cap L^2(\Omega;L^2(\tau,T;\vv))$.
\end{theorem}
\begin{proof} The proof of this theorem follows a standard scheme for example, \cite [Theorem 3.2]{L1} and \cite[Theorem 2.1]{R1}, but with particular technical difficulties caused by  L\'evy noise and fractional Laplace operator. We will split the proof into 5 steps.

{\bf Step 1.} As $\lambda$ is  $\sigma$-finite  on the Polish space $E$, there exist measurable subsets $E_m\nearrow E$ satisfying $\lambda(E_m)<\infty$ for all $m\in\mathbb{N}$ and $\cup_{m=1}^{\infty}E_m=E$. Now, for each $m\in\mathbb{N}$ and $t>\tau$, we define the function {$u^m_0\equiv u_0$} and consider recursively the equations,
\begin{equation}\label{eq3-3}
\begin{split}
u_n^m(t)&=u_0-\int_\tau^t(-\Delta)^{\gamma}u^m_n(s)ds-\delta\int_\tau^tu^m_n(s)ds-\int_\tau^tF(u^m_n(s))ds\\
&~~+\int_\tau^tg(s,u^m_{n-1}(s))dW(s)+\int_\tau^t\int_{E_m}h(u^m_{n-1}(s-),\xi)\tilde{N}(ds,d\xi),\quad \mathbb{P}\mbox{-}a.s.
\end{split}
\end{equation}
Applying the It\^o formula to $|u^m_n|^2$, we have
\begin{equation}\label{eq3-4}
\begin{split}
 |\um(t)|^2&=|u_0|^2-2\int_\tau^t<(-\Delta)^{\gamma}\um(s)+\delta\um(s)+F(\um(s)),\um(s)>ds\\
 &~~+2\int_\tau^t\left(\um(s),g(s,u^m_{n-1}(s))dW(s)\right)+\int_\tau^t\|g(s,u^m_{n-1}(s))\|^2_{\mathcal{L}_2(U;\hh)}ds\\
 &~~+2\int_\tau^t\int_{E_m}\left(\um(s),h(u^m_{n-1}(s-),\xi)\right)\tilde{N}(ds,d\xi)+\int_\tau^t\int_{E_m}|h(u^m_{n-1}(s-),\xi)|^2N(ds,d\xi).
 \end{split}
\end{equation}
By definition of the fractional Laplacian operator $(-\Delta)^{\gamma}$ and \eqref{eq2-11}-\eqref{eq2-12}, setting $\eta=\min\left\{\frac{C(d,\gamma)}{2},\delta\right\}$, we obtain
\begin{eqnarray}\label{eq3-5}
& -2\int_\tau^t<(-\Delta)^{\gamma}\um(s)+\delta\um(s),\um(s)>ds=-2\int_\tau^t
\left(|(-\Delta)^{\frac{\gamma}{2}}\um(s)|^2+\delta |\um(s)|^2\right)ds\\
&={-2\int_\tau^t\left(\frac{C(d,\gamma)}{2}\|\um(s)\|^2_{\dot{H}^{\gamma}(\mathbb{R}^d)}+\delta |\um(s)|^2\right)ds
\leq -2\eta\int_\tau^t\|\um(s)\|^2ds}.\notag
\end{eqnarray}
 Making use of condition $(F_2)$, we derive
 \begin{eqnarray}\label{eq3-6}
 -2\int_\tau^t\int_{\mathbb{R}^d}F(\um(s))\um(s)dxds&=&-2\int_\tau^t\int_{\mathcal{O}}F(\um(s))\um(s)dxds\leq2\int_\tau^t\int_{\mathcal{O}}\left(k_2-k_3|\um(s)|^{p+1}\right)dxds\notag\\
 &\leq& 2k_2|\mathcal{O}|(t-\tau)-2k_3\|\um\|^{p+1}_{L^{p+1}(\tau,t;L^{p+1}(\mathbb{R}^d))}.
 \end{eqnarray}
  Substituting \eqref{eq3-5}-\eqref{eq3-6} into \eqref{eq3-4},  taking supremum with respect to $t\in[\tau,T']$ for any {$\tau\leq T'\leq T$} and expectation, we find
\begin{equation}\label{eq3-7}
\begin{split}
&\max\left\{\e\left[\sup_{t\in[\tau,T']}|\um(t)|^2\right],~2\eta\e\int_\tau^{T'}\|\um(t)\|^2dt,~2k_3\e\int_\tau^{T'}\|\um(t)\|_{p+1}^{p+1}dt\right\}\\
&\quad \leq \e|u_0|^2+2k_2|\mathcal{O}|(T-\tau)
+2\e\left(\sup_{t\in[\tau,T']}\int_\tau^t\left(\um(s),g(s,\um(s))dW(s)\right)\right)\\
&~~~~+\e\int_\tau^{T'}\|g(t,u^m_{n-1}(t))\|^2_{\mathcal{L}_2(U;\hh)}dt+
\e\int_\tau^{T'}\int_{E_m}|h(u^m_{n-1}(t{-}),\xi)|^2N(dt,d\xi)\\
&~~~~+2\e \left(\sup_{t\in[\tau,T']}\int_\tau^t\int_{E_m}(\um(s),h(u^m_{n-1}(s-),\xi))\tilde{N}(ds,d\xi)\right).
\end{split}
\end{equation}
 By applying assumption $(g_2)$,  the Burkholder-Davis-Gundy and Young inequalities,  we have
\begin{equation}\label{eq3-8}
\begin{split}
&~\qquad\quad 2\e\left(\sup_{t\in[\tau,T']}\int_\tau^t\left(\um(s),g(s,u^m_{n-1}(s))dW(s)\right)\right)\\
&\qquad\leq 2C_b\e\left(\int_\tau^{T'}|\um(s)|^2\|g(s,u^m_{n-1}(s))\|^2_{\mathcal{L}_2(U;\hh)}ds\right)^{\frac{1}{2}}\\
&\qquad\leq 2C_b\e\left(\sup_{t\in[\tau,T']}|\um(t)|^2\int_\tau^{T'}\|g(t,u^m_{n-1}(t))\|^2_{\mathcal{L}_2(U;\hh)}dt\right)^{\frac{1}{2}}\\
&\qquad\leq \frac{1}{4}\e\left[\sup_{t\in[\tau,T']}|\um(t)|^2\right]+4C_b^2\e\int_\tau^{T'}\|g(t,u^m_{n-1}(t))\|^2_{\mathcal{L}_2(U;\hh)}dt\\
&\qquad\leq \frac{1}{4}\e\left[\sup_{t\in[\tau,T']}|\um(t)|^2\right]+4C_b^2C_g(T-\tau)+4C_b^2C_g\e\int_\tau^{T'}|u^m_{n-1}(t)|^2dt,
\end{split}
\end{equation}
here and in the sequel, $C_b$ is the constant obtained from the Burkholder-Davis-Gundy inequality for Brownian motion and the Poisson process. Similarly, we conclude from assumption $(h_2)$, the Burkholder-Davis-Gundy and Young inequalities that
\begin{equation}\label{eq3-9}
\begin{split}
&~\quad2\e \left(\sup_{t\in[\tau,T']}\int_\tau^t\int_{E_m}(\um(s),h(u^m_{n-1}(s-),\xi))\tilde{N}(ds,d\xi)\right)\\
&\leq 2C_b\e\left(\int_\tau^{T'}\int_{E_m}|\um(s)|^2|h(u^m_{n-1}(s-),\xi)|^2\lambda(d\xi)ds\right)^{\frac{1}{2}}
\\
&\leq 2C_b\e\left(\sup_{t\in[\tau,T']}|\um(t)|^2\int_\tau^{T'}\int_{E_m}|h(u^m_{n-1}(s-),\xi)|^2\lambda(d\xi)ds\right)^{\frac{1}{2}}\\
&\leq\frac{1}{4}\e\left[\sup_{t\in[\tau,T']}|\um(t)|^2\right]+4C_b^2\e\int_\tau^{T'}\int_{E_m}|h(u^m_{n-1}(s-),\xi)|^2\lambda(d\xi)ds\\
&\leq \frac{1}{4}\e\left[\sup_{t\in[\tau,T']}|\um(t)|^2\right]+4C_b^2C_h(T-\tau)+4C_b^2C_h\e\int_\tau^{T'}|u^m_{n-1}(t)|^2dt.
\end{split}
\end{equation}
By means of  $(g_2)$ and $(h_2)$, we infer
\begin{equation}\label{eq3-10}
\e\int_\tau^{T'}\|g(t,u^m_{n-1}(t))\|^2_{\mathcal{L}_2(U;\hh)}dt\leq C_g(T-\tau)+C_g\e\int_\tau^{T'}|u^m_{n-1}(t)|^2dt, \ \forall T'\in[\tau,T],
\end{equation}
and 
\begin{equation}\label{eq3-11}
\e\int_\tau^{T'}\int_{E_m}|h(u^m_{n-1}(t-),\xi)|^2\lambda(d\xi)dt\leq C_h(T-\tau)+C_h\e\int_\tau^{T'}|u^m_{n-1}(t)|^2dt,  \ \forall T'\in[\tau,T],
\end{equation}
respectively. Thus, it follows from \eqref{eq3-7}-\eqref{eq3-11} that, for all $T'\in(\tau,T]$,
\begin{equation}\label{eq3-12}
\begin{split}
\e\left[\sup_{t\in[\tau,T']}|\um(t)|^2\right]\leq 2\e|u_0|^2+(4k_2|\mathcal{O}|+2C_1)(T-\tau)+2C_1\int_\tau^{T'}\e\left[\sup_{s\in[\tau,t]}|u^m_{n-1}(s)|^2\right]dt,
\end{split}
\end{equation}
where we have denoted by $C_1:=4C_b^2C_g+4C_b^2C_h+C_g+C_h$. Let us define  
$\mathcal{U}^m_N(T'):=\sup_{n\leq N}\e \left[\sup_{t\in[\tau,T']}|\um(t)|^2\right],$ for all $N\in\mathbb{N}, T'\in[\tau,T].$
 Subsequently, inequality \eqref{eq3-12} implies for each $T'\in[\tau,T]$, 
 \begin{equation*}
 \mathcal{U}_N^m(T')\leq 2\e|u_0|^2+(4k_2|\mathcal{O}|+2C_1)(T-\tau)+2C_1\int_\tau^{T'}\mathcal{U}^m_N(t)dt.
 \end{equation*}
 The Gronwall lemma implies  for each $m\in\mathbb{N}$ and for any $\tau\leq T'\leq T$ that, 
 \begin{equation*}
 \e \left[\sup_{t\in[\tau,T']}|\um(t)|^2\right]\leq \left(2\e|u_0|^2+(4k_2|\mathcal{O}|+2C_1)(T-\tau)\right)e^{2C_1(T-\tau)}:=C_2(T),\ \forall n\in\mathbb{N},
 \end{equation*}
where $C_2(T)$ is a positive constant depending on $T$.  This, together with \eqref{eq3-7}, shows
 \begin{equation*}
\e\int_\tau^T\|\um(t)\|^2dt+\e\int_\tau^T\|\um(t)\|_{p+1}^{p+1}dt\leq C_3(T),
 \end{equation*} 
 for a positive constant $C_3(T)$. In conclusion, we proved in Step 1 that  for all $m, n\in \mathbb{N}$, there exists a constant $C_4>0$ depending on $T$, such that
 \begin{equation}\label{eq3-13}
\e \left[\sup_{t\in[\tau,T]}|\um(t)|^2\right]+\e\int_\tau^T\|\um(t)\|^2dt+\e\int_\tau^T\|\um(t)\|_{p+1}^{p+1}dt\leq C_4(T).
 \end{equation}

{\bf Step 2.} For each $m\in\mathbb{N}$ and $n_1,n_2\in\mathbb{N}$, define $\chi_{n_1,n_2}^m=u^m_{n_1}(t)-u^m_{n_2}(t)$, using similar arguments as in Step 1, by  $(F_3)$, $(g_1)$ and $(h_1)$,  the Burkholder-Davis-Gundy and Young inequalities, we obtain
\begin{eqnarray}\label{eq3-14}
&&\hskip1cm\frac{1}{2}\e\left[\sup_{t\in[\tau,T']}|\x(t)|^2\right]+2\eta\e\int_\tau^{T'}\|\x(t)\|^2dt\\
&&\leq 4k_4\int_\tau^{T'}\e\left[\sup_{s\in[\tau,t]}|\x(s)|^2\right]dt
+2C_5\int_\tau^{T'}\e\left[\sup_{s\in[\tau,t]}|\chi^m_{n_1-1,n_2-1}(s)|^2\right]dt,\ \mbox{for all}\ \tau\leq T'\leq T,\notag
\end{eqnarray}
where we have denoted by $C_5=4C_b^2L_g+L_g+4C_b^2L_h+L_h$.

 Let $\mathcal{V}^m(T')=\varlimsup\limits_{n_1,n_2\rightarrow\infty}\e\left[\sup_{t\in[\tau,T']}|\x(t)|^2\right].$
Then, {by the Fatou-Lebesgue Theorem}  and \eqref{eq3-14}, we have $\mathcal{V}^m(T')\leq 4(2k_4+C_5)\int_\tau^{T'}\mathcal{V}^m(t)dt$ for all $T'\in[\tau,T]$.
 {The Gronwall lemma implies}  that
 \begin{equation}\label{eq3-15}
 \mathcal{V}^m(T)=\varlimsup\limits_{n_1,n_2\rightarrow\infty}\e\left[\sup_{t\in[\tau,T]}|\x(t)|^2\right]=0\quad \mbox{and}\quad \varlimsup\limits_{n_1,n_2\rightarrow\infty}\e\int_\tau^T\|\x(t)\|^2dt=0.
 \end{equation}
Hence, for each $m\in\mathbb{N}$, there exists an adapted process $u^m\in L^2(\Omega;\mathbf{D}([\tau,T];\hh))\cap L^2(\Omega;L^2(\tau,T;\vv))$ such that
$\lim_{n\rightarrow \infty}\e\left[\sup_{t\in[\tau,T]}|\um(t)-u^m(t)|^2\right]=\lim_{n\rightarrow \infty}\e\int_\tau^T\|\um(t)-u^m(t)\|^2dt=0.$
Additionally, thanks to \eqref{eq3-13}, we immediately derive
\begin{equation}\label{eq3-16}
\begin{split}
\e \left[\sup_{t\in[\tau,T]}|u^m(t)|^2\right]+\e\int_\tau^T\|u^m(t)\|^2dt+\e\int_\tau^T\|u^m(t)\|_{p+1}^{p+1}dt\leq C_4(T).
\end{split}
\end{equation} 
 Now, for each $m\in\mathbb{N}$, taking limit in \eqref{eq3-3} as $n\rightarrow \infty$, by means of the continuity of function $F$, Lipschitz condition imposed on $g$ and $h$, it is easy to show that $u^m$ is the unique solution of the equation
 \begin{equation}\label{eq3-17}
 \begin{split}
u^m(t)&=u_0-\int_\tau^t(-\Delta)^{\gamma}u^m(s)ds-\delta\int_\tau^tu^m(s)ds-\int_\tau^tF(u^m(s))ds\\
&~~+\int_\tau^tg(s,u^m(s))dW(s)+\int_\tau^t\int_{E_m}h(u^m(s-),\xi)\tilde{N}(ds,d\xi),\quad \tau<t\leq T,~ \mathbb{P}\mbox{-}a.s.
\end{split}
\end{equation}
Notice that, equality \eqref{eq3-17} holds  in $(\vv\cap L^{p+1}(\mathbb{R}^d))^*$.

{\bf Step 3.}  For $m_1$, $m_2\in\mathbb{N}$ with $m_2<m_1$,  we have $\lambda(E_{m_2})<\lambda(E_{m_1})$ since $E_m$ is increasing. Define $\mathbf{B}^{m_1,m_2}(t)=u^{m_1}(t)-u^{m_2}(t)$, applying the It\^o formula to $|\p(t)|^2$, proceeding likewise as in Step 1, we obtain for all $T'\in[\tau,T]$ that
\begin{eqnarray}\label{eq3-18}
&&\max\left\{\e\left[\sup_{t\in[\tau,{T'}]}|\p(t)|^2\right],~2\eta\e\int_\tau^{T'}\|\p(t)\|^2dt\right\}\notag\\
&&\leq 2k_4\e\int_\tau^{T'}|\p(t)|^2dt\notag\\
&&~~+2\e\left(\sup_{t\in[\tau,T']}\int_\tau^t\left(\p(s),(g(s,u^{m_1}(s))-g(s,u^{m_2}(s)))dW(s)\right)\right)\notag\\
&&~~+\e\int_\tau^{T'}\|g(s,u^{m_1}(s))-g(s,u^{m_2}(s))\|^2_{\mathcal{L}_2(U;\hh)}ds\\
&&~~+2\e\left(\sup_{t\in[\tau,T']}\int_\tau^t\int_{E_{m_2}}(\p(s),h(u^{m_1}(s-),\xi)-h(u^{m_2}(s-),\xi))\tilde{N}(ds,d\xi)\right)\notag\\
&&~~+\e\int_\tau^{T'}\int_{E_{m_2}}|h(u^{m_1}(s-),\xi)-h(u^{m_2}(s-),\xi)|^2N(ds,d\xi)\notag\\
&&~~+2\e\left(\sup_{t\in[\tau,T']}\int_\tau^t\int_{E_{m_1}\backslash E_{m_2}}(\p(s),h(u^{m_1}(s-),\xi))\tilde{N}(ds,d\xi)\right)\notag\\
&&~~+\e\int_\tau^{T'}\int_{E_{m_1}\backslash E_{m_2}}|h(u^{m_1}(s-),\xi)|^2N(ds,d\xi)\notag\\
&&\leq 2k_4\int_\tau^{T'}\e\left[\sup_{s\in[\tau,t]}|\p(s)|^2\right]dt+I_1+I_2+I_3+I_4+I_5+I_6.\notag
\end{eqnarray}
For $I_1$, similar to \eqref{eq3-8}, by $ (g_1)$,  the Burkholder-Davis-Gundy and Young inequalities,  we have
\begin{eqnarray}\label{eq3-19}
&&2\e\left(\sup_{t\in[\tau,T']}\int_\tau^t\left(\p(s),(g(s,u^{m_1}(s))-g(s,u^{m_2}(s)))dW(s)\right)\right)\notag\\
&&\leq \frac{1}{4}\e\left[\sup_{t\in[\tau,T']}|\p(t)|^2\right]+4C_b^2L_g\int_\tau^{T'}\e\left[\sup_{s\in[\tau,t]}|\p(s)|^2\right]dt.
\end{eqnarray}
 For $I_3$, similar to \eqref{eq3-9}, making use of $(h_1)$, the Burkholder-Davis-Gundy and Young inequalities,  we obtain
 \begin{equation}\label{eq3-20}
 \begin{split}
 &2\e\left(\sup_{t\in[\tau,T']}\int_\tau^t\int_{E_{m_2}}(\p(s),h(u^{m_1}(s-),\xi)-h(u^{m_2}(s-),\xi))\tilde{N}(ds,d\xi)\right)\\[1.0ex]
 &\leq \frac{1}{4}\e\left[\sup_{t\in[\tau,T']}|\p(t)|^2\right]+4C_b^2L_h\int_\tau^{T'}\e\left[\sup_{s\in[\tau,t]}|\p(s)|^2\right]dt.
 \end{split}
 \end{equation}
 For $I_2$ and $I_4$, similar to \eqref{eq3-10}-\eqref{eq3-11}, by means of $(g_1)$ and $(h_1)$, we derive
 \begin{equation}\label{eq3-21}
 I_2=\e\int_\tau^{T'}\|g(s,u^{m_1}(s))-g(s,u^{m_2}(s))\|^2_{\mathcal{L}_2(U;\hh)}ds\leq L_g\int_\tau^{T'}\e\left[\sup_{s\in[\tau,t]}|\p(s)|^2\right]dt,
 \end{equation}
 and
 \begin{equation}\label{eq3-22}
 I_4=\e\int_\tau^{T'}\int_{E_{m_2}}|h(u^{m_1}(s-),\xi)-h(u^{m_2}(s-),\xi)|^2N(ds,d\xi)\leq L_h\int_\tau^{T'}\e\left[\sup_{s\in[\tau,t]}|\p(s)|^2\right]dt,
 \end{equation} 
respectively. For $I_5$,  by assumption $(h_2)$, the Burkholder-Davis-Gundy and Young inequalities,  we have
\begin{equation}\label{eq3-23}
\begin{split}
I_5&=2\e\left(\sup_{t\in[\tau,T']}\int_\tau^t\int_{E_{m_1}\backslash E_{m_2}}(\p(s),h(u^{m_1}(s-),\xi))\tilde{N}(ds,d\xi)\right)\\
&\leq \frac{1}{4}\e\left[\sup_{t\in[\tau,T']}|\p(t)|^2\right]+4C_b^2
\e\int_\tau^{T'}\int_{E_{m_1}\backslash E_{m_2}}|h(u^{m_1}(t-),\xi)|^2\lambda(d\xi)dt.
\end{split}
\end{equation} 
Consequently, substituting \eqref{eq3-19}-\eqref{eq3-23} to \eqref{eq3-18}, we find
\begin{equation}\label{eq3-24}
\begin{split}
\frac{1}{4}\e\left[\sup_{t\in[\tau,T']}|\p(t)|^2\right]&\leq (2k_4+4C_b^2L_g+L_g+4C_b^2L_h+L_h)\int_\tau^{T'}\e\left[\sup_{s\in[\tau,t]}|\p(s)|^2\right]dt\\[1.0ex]
&~~+(4C_b^2+1)\e\int_\tau^{T'}\int_{E_{m_1}\backslash E_{m_2}}|h(u^{m_1}(t-),\xi)|^2\lambda(d\xi)dt.
\end{split}
\end{equation}
It follows from the fact $E_{m_2}\subset E_{m_1}$ with $\lambda(E_{m_1}\backslash E_{m_2})\rightarrow 0$ as $m_1,m_2\rightarrow \infty$,  assumption $(h_2)$ and  the property of absolutely Lebesgue integrable function that, for any $\varepsilon>0$, there exist $M(\varepsilon)>0$ and $\delta>0$, such that for all $m_1,m_2\geq M(\varepsilon)$, 
\begin{equation}\label{eq3-25}
\lambda(E_{m_1}\backslash E_{m_2})<\delta\quad\mbox{and}\quad 
\int_{E_{m_1}\backslash E_{m_2}}{|h(u^{m_1}(t^-),\xi)|^2}\lambda(d\xi)<\varepsilon.
\end{equation}
Immediately, the Gronwall lemma, together with \eqref{eq3-24}-\eqref{eq3-25} and \eqref{eq3-18}, implies that
\begin{equation*}
\e\left[\sup_{t\in[\tau,T]}|\p(t)|^2\right]\to 0,\quad\mbox{and}\quad  \e\int_\tau^T\|\p(t)\|^2dt\to0,\qquad \mbox{as}~~m_1,~m_2\rightarrow\infty.
\end{equation*}
 Therefore, there exists an $\mathcal{F}_t$-adapted process $u\in L^2(\Omega;\mathbf{D}([\tau,T];\hh))\cap L^2(\Omega;L^2(\tau,T;\vv))$ such that $\lim_{m\rightarrow \infty}\e\left[\sup_{t\in[\tau,T]}|u^m(t)-u(t)|^2\right]=\lim_{m\rightarrow \infty}\e\int_\tau^T\|u^m(t)-u(t)\|^2dt=0,$
 which, combining  with \eqref{eq3-16}, yields 
  \begin{equation*}
\e \left[\sup_{t\in[\tau,T]}|u(t)|^2\right]+\e\int_\tau^T\|u(t)\|^2dt+\e\int_\tau^T\|u(t)\|_{p+1}^{p+1}dt\leq C_4(T).
\end{equation*} 
 Eventually, taking limit in \eqref{eq3-17} shows that $u$ is the unique solution of \eqref{eq2-13} on the interval $[\tau,T]$. 
 
 {\bf Step 4.} By repeating the above arguments, we obtain the existence of the unique solution of \eqref{eq2-13} on the interval $[T,2T-\tau]$, which finally leads to the completion of the global existence and uniqueness of solutions to \eqref{eq2-13}-\eqref{eq2-14} by further iterations.

 {\bf Step 5.} Continuity of solutions with respect to initial data. Let $u_{0,1}$, $u_{0,2}\in L^2(\Omega;\hh)$ be two  $\mathcal{F}_\tau$-measurable initial data, $u_1$ and $u_2$ are the corresponding solutions of \eqref{eq2-13}-\eqref{eq2-14} on $[\tau,T]$ for any $T>\tau$, respectively. Denote by $\bar{u}=u_1-u_2$ and $\bar{u}_0=u_{0,1}-u_{0,2}$. Then by the It\^o formula,  definition of fractional Laplacian operator $(-\Delta)^{\gamma}$, \eqref{eq2-11}-\eqref{eq2-12} and condition $(F_3)$, similar to \eqref{eq3-5}, we have  for every $m\in\mathbb{N}$,
 \begin{eqnarray}\label{eq3-26}
 |\bar{u}^m(t)|^2+2\eta\int_{\tau}^t\|\bar{u}^m(s)\|^2ds&\leq& |\bar{u}_0|^2+2\int_{\tau}^t\left(\bar{u}^m(s),(g(s,u^m_1(s))-g(s,u^m_2(s)))dW(s)\right)\notag\\
 &&+2k_4\int_{\tau}^t|\bar{u}^m(s)|^2ds+\int_\tau^t\|g(s,u^m_1(s))-g(s,u^m_2(s))\|^2_{\mathcal{L}_2(U;\hh)}ds\notag\\
 &&+2\int_\tau^t\int_{E_m}(\bar{u}^m(s),h(u^m_1(s-),\xi)-h(u^m_2(s-),\xi))\tilde{N}(ds,d\xi)\notag\\
 &&+\int_\tau^t\int_{E_m}|h(u_1^m(s-),\xi)-h(u_2^m(s-),\xi)|^2N(ds,d\xi).
 \end{eqnarray}
On the one hand, by \eqref{eq3-26}, we find that for all $\tau\leq T'\leq T$,
\begin{eqnarray}\label{eq3-27}
\e\left[\sup_{\tau\leq t\leq T'}|\bar{u}^m(t)|^2\right]&\leq& \e|\bar{u}_0|^2+2\e\left(\sup_{\tau\leq t\leq T'}\int_\tau^t(\bar{u}^m(s),(g(s,u^m_1(s))-g(s,u^m_2(s)))dW(s))\right)\notag\\
&&+2k_4\int_\tau^{T'}\e\left[\sup_{\tau\leq s\leq t}|\bar{u}^m(s)|^2\right]dt+\e\int_\tau^{T'}\|g(s,u^m_1(s))-g(s,u^m_2(s))\|^2_{\mathcal{L}_2(U;\hh)}ds\notag\\
&&+2\e\left(\sup_{\tau\leq t\leq T'}\int_\tau^t\int_{E_m}(\bar{u}^m(s),h(u^m_1(s-),\xi)-h(u^m_2(s-),\xi))\tilde{N}(ds,d\xi)\right)\\
&&+\e\int_\tau^{T'}\int_{E_m}|h(u_1^m(s-),\xi)-h(u_2^m(s-),\xi)|^2N(ds,d\xi).\notag
\end{eqnarray}
 Similar to \eqref{eq3-8}-\eqref{eq3-9}, by the Burkholder-Davis-Gundy and Young inequalities,  $(g_1)$ and $(h_1)$, we have
 \begin{equation}\label{eq3-28}
 \begin{split}
&\quad 2\e\left(\sup_{\tau\leq t\leq T'}\int_\tau^t(\bar{u}^m(s),(g(s,u^m_1(s))-g(s,u^m_2(s)))dW(s))\right)\\
 &\leq \frac{1}{4}\e\left[\sup_{\tau\leq t\leq T'}|\bar{u}^m(t)|^2\right]+4C_b^2L_g\int_\tau^{T'}\e\left[\sup_{\tau\leq s\leq t}|\bar{u}^m(s)|^2\right]dt,
 \end{split}
 \end{equation}
 and
\begin{equation}\label{eq3-29}
 \begin{split}
&\quad 2\e\left(\sup_{\tau\leq t\leq T'}\int_\tau^t\int_{E_m}(\bar{u}^m(s),h(u^m_1(s-),\xi)-h(u^m_2(s-),\xi))\tilde{N}(ds,d\xi)\right)\\
 &\leq \frac{1}{4}\e\left[\sup_{\tau\leq t\leq T'}|\bar{u}^m(t)|^2\right]+4C_b^2L_h\int_\tau^{T'}\e\left[\sup_{\tau\leq s\leq t}|\bar{u}^m(s)|^2\right]dt,
 \end{split}
 \end{equation} 
respectively. It follows from \eqref{eq3-27}-\eqref{eq3-29}, conditions $(g_1)$ and $(h_1)$ that
 \begin{equation*}
\e\left[\sup_{\tau\leq t\leq T'}|\bar{u}^m(t)|^2\right]\leq 2\e|\bar{u}_0|^2+2(2k_4+4C_b^2L_g+4C_b^2L_h+L_g+L_h)\int_\tau^{T'}\e\left[\sup_{\tau\leq s\leq t}|\bar{u}^m(s)|^2\right]dt.
\end{equation*}
 Applying the Gronwall lemma  to the above inequality, we obtain in particular for $T'=T$,
 \begin{equation}\label{eq3-30}
 \e\left[\sup_{\tau\leq t\leq T}|\bar{u}^m(t)|^2\right]\leq 2e^{C_6(T-\tau)}\e|\bar{u}_0|^2,
 \end{equation}
 where $C_6:=2\left(2k_4+4C_b^2L_g+4C_b^2L_h+L_g+L_h\right)$. By \eqref{eq3-26} and \eqref{eq3-30}, for some $C_7:=C_7(T)$, we derive
 \begin{equation}\label{eq3-31}
 \e\int_\tau^T\|\bar{u}^m(t)\|^2dt\leq C_7\e{|\bar{u}_0|^2}.
 \end{equation}
In Step 3, we have proved for each initial value $u_0\in L^2(\Omega;\hh)$ and every $m\in\mathbb{N}$, the corresponding solution sequence $\{u^m\}$ is Cauchy in $L^2(\Omega;\mathbf{D}([\tau,T];\hh))\cap L^2(\Omega;L^2(\tau,T;\vv))$. Therefore, $u_1$, $u_2\in L^2(\Omega;\mathbf{D}([\tau,T];\hh))\cap L^2(\Omega;L^2(\tau,T;\vv))$ and satisfy \eqref{eq3-30}-\eqref{eq3-31}. Moreover, there exists a positive constant $C(T)$ such that,
\begin{equation}\label{eq3-32}
\|u_1-u_2\|_{L^2(\Omega;\mathbf{D}([\tau,T];\hh))}^2+\|u_1-u_2\|^2_{L^2(\Omega;L^2(\tau,T;\vv))}\leq C\|u_{0,1}-u_{0,2}\|^2_{L^2(\Omega;\hh)}.
\end{equation} 
Namely, the solution depends continuously on initial data. The proof of this theorem is complete. 
\end{proof}

\begin{remark}\label{remark}
Notice that, under assumptions of Theorem \ref{thm3-2}, if $u$ is the unique solution to problem \eqref{eq2-13}-\eqref{eq2-14} corresponding to the initial value $u_0\in L^2(\Omega;\hh)$, then there exists a sequence $u^m \in L^2(\Omega;\mathbf{D}([\tau,T];\hh))\cap L^2(\Omega;L^2(\tau,T;\vv))\cap L^{p+1}(\Omega;L^{p+1}(\tau,T;L^{p+1}(\mathbb{R}^{d}))) (m\geq1)$, which converges to $u$ in $L^2(\Omega;\mathbf{D}([\tau,T];\hh))\cap L^2(\Omega;L^2(\tau,T;\vv))$ and satisfies \eqref{eq3-17}. In other words, each $u^m$ is solution to problem \eqref{eq2-13}-\eqref{eq2-14} but replacing $E$ by $E_m$. 
This fact has been used in Step 5 in the previous proof and will be used repeatedly in the following sections.
\end{remark}
 
\section{Existence of weak mean random attractors}\label{s4}

This section is devoted to the existence and uniqueness of weak mean random attractors for the non-autonomous fractional stochastic differential  equations \eqref{eq2-13}-\eqref{eq2-14}. To this end,  we first define a mean random dynamical system for \eqref{eq2-13}-\eqref{eq2-14}, then prove the existence and uniqueness of weak pullback mean random attractors. 

Observe that, it follows from Theorem \ref{thm3-2} that for every $\tau\in\mathbb{R}$ and every  $\mathcal{F}_\tau$-measurable initial datum $u_0\in L^2(\Omega;\hh)$, problem \eqref{eq2-13}-\eqref{eq2-14} has a unique c\`adl\`ag $\hh$-valued $\mathcal{F}_t$-adapted solution $u(t,\tau,u_0)$ with initial condition $u_0$ at $\tau$ in the sense of Definition \ref{def3-1}. Theorem \ref{thm3-2} presented that $u(\cdot,\tau,u_0)\in L^2(\Omega;\mathbf{D}([\tau,\infty);\hh))$, which implies that $u\in \mathbf{D}([\tau,\infty);L^2(\Omega;\hh))$. In this way, we are able to define a cocycle generated by the problem under consideration. Given $t\in\mathbb{R}^+$ and $\tau\in\mathbb{R}$, let $\Phi(t,\tau)$ be a mapping from $L^2(\Omega,\mathcal{F}_\tau;\hh)$ to $L^2(\Omega,\mathcal{F}_{t+\tau};\hh)$ defined by
$ \Phi(t,\tau)(u_0)=u(t+\tau,\tau,u_0),$
where $u_0\in L^2(\Omega,\mathcal{F}_{\tau};\hh)$. The uniqueness of solution to \eqref{eq2-13}-\eqref{eq2-14} implies  that for every $t,s>0$ and $\tau\in\mathbb{R}$,
$\Phi(t+s,\tau)=\Phi(t,s+\tau)\circ\Phi(s,\tau).$
This cocycle $\Phi$ is called the mean random dynamical system generated by \eqref{eq2-13}-\eqref{eq2-14} on $L^2(\Omega,\mathcal{F};\hh)$. We will study the existence and uniqueness of weak pullback random attractors for $\Phi$.

Let $B=\{B(\tau)\subset L^2(\Omega,\mathcal{F}_\tau;\hh):\tau\in\mathbb{R}\}$ be a family of nonempty bounded sets such that
\begin{equation}\label{eq4-0}
\lim_{\tau\rightarrow -\infty}e^{(2\delta-L_1)\tau}\|B(\tau)\|^2_{L^2(\Omega,\mathcal{F}_\tau;\hh)}=0,
\end{equation}
where $2\delta>L_1$ with $L_1=C_g+C_h$ and $\|\mathcal{Q}\|_{L^2(\Omega,\mathcal{F}_\tau;\hh)}=\sup_{u\in \mathcal{Q}}\|u\|_{L^2(\Omega,\mathcal{F}_\tau;\hh)}$ for a subset $\mathcal{Q}$ in $L^2(\Omega,\mathcal{F}_\tau;\hh)$. We will use $\mathcal{D}$ to denote the collection of all families of nonempty bounded sets satisfying \eqref{eq4-0}.

We will first derive uniform estimates on the solutions of \eqref{eq2-13}-\eqref{eq2-14}, then construct a $\mathcal{D}$-pullback absorbing set for the system $\Phi$.

\begin{lemma}\label{lem4-1}
Suppose $(F_1)$-$(F_3)$, $(g_1)$-$(g_3)$ and $(h_1)$-$(h_3)$ hold. In addition, assume $2\delta>L_1:=C_g+C_h$. Then for every $\tau\in\mathbb{R}$ and $B\in\mathcal{D}$, there exists $T=T(\tau,B)>0$ such that for all $t\geq T$, the solution $u$ of \eqref{eq2-13}-\eqref{eq2-14} satisfies 
\begin{equation*}
\e|u(\tau,\tau-t,u_0)|^2\leq \frac{2k_2|\mathcal{O}|+L_1}{2\delta-L_1}+1,\qquad \forall u_0\in B(\tau-t).
\end{equation*}
\end{lemma}
\begin{proof} We will split the proof into two steps. 

{\bf Step 1.}  As $\lambda$ is $\sigma$-finite on the Polish space $E$, there exist measurable subsets $E_m\nearrow E$ satisfying $\lambda(E_m)<\infty$ for all $m\in\mathbb{N}$ and $\cup_{m=1}^{\infty}E_m=E$. Taking into account Remark \ref{remark}, for each $m\in\mathbb{N}$,
applying   It\^o's  formula to $e^{(2\delta-L_1)t}|\ur|^2$ with $r\geq \tau-t$ (see, for example, \cite{C4}), we obtain 
\begin{eqnarray}\label{eq4-1}
&&e^{(2\delta-L_1)r}|\ur|^2+2\int_{\tau-t}^re^{(2\delta-L_1)s}|(-\Delta)^{\frac{\gamma}{2}}\us|^2ds\notag\\
&&+2\delta\int_{\tau-t}^re^{(2\delta-L_1)s}|\us|^2ds\notag\\
&&+2\int_{\tau-t}^re^{(2\delta-L_1)s}<F(\us),\us>ds\\
&&=e^{(2\delta-L_1)(\tau-t)}|u_0|^2+(2\delta-L_1)\int_{\tau-t}^re^{(2\delta-L_1)s}|\us|^2ds\notag\\
&&+2\int_{\tau-t}^re^{(2\delta-L_1)s}(\us,g(s,\us)dW(s))\notag\\
&&+\int_{\tau-t}^re^{(2\delta-L_1)s}\|g(s,\us)\|^2_{\mathcal{L}_2(U;\hh)}ds\notag\\
&&+2\int_{\tau-t}^r\int_{E_m}e^{(2\delta-L_1)s}(\us,h(u(s-,\tau-t,u_0),\xi))\tilde{N}(ds,d\xi)\notag\\
&&+\int_{\tau-t}^r\int_{E_m}e^{(2\delta-L_1)s}|h(u^m(s-,\tau-t,u_0),\xi)|^2N(ds,d\xi),\notag
\end{eqnarray}
which implies that for all $r\geq \tau-t$,
\begin{eqnarray}\label{eq4-2}
&& e^{(2\delta-L_1)r}\e|\ur|^2+C(d,\gamma)\int_{\tau-t}^re^{(2\delta-L_1)s}\e\|\us\|_{\dot{\mathbb{H}}^{\gamma}(\mathbb{R}^d)}^2ds\notag\\
&&+2\delta\int_{\tau-t}^re^{(2\delta-L_1)s}\e|\us|^2ds\notag\\
&&+2\int_{\tau-t}^re^{(2\delta-L_1)s}\e<F(\us),\us>ds\notag\\
&&=e^{(2\delta-L_1)(\tau-t)}\e|u_0|^2+(2\delta-L_1)\int_{\tau-t}^re^{(2\delta-L_1)s}\e|\us|^2ds\\
&&+\int_{\tau-t}^re^{(2\delta-L_1)s}\e\|g(s,\us)\|^2_{\mathcal{L}_2(U;\hh)}ds\notag\\
&&+\int_{\tau-t}^r\int_{E_m}e^{(2\delta-L_1)s}\e|h(u^m(s-,\tau-t,u_0),\xi)|^2N(ds,d\xi).\notag
\end{eqnarray}
We now do estimates one by one for \eqref{eq4-2}. On the one hand, by $(F_2)$, we have
\begin{equation}\label{eq4-4}
\begin{split}
&2\e<F(\ur),\ur>=2\e\int_{\mathbb{R}^d}F(\ur)\ur dx\\
&=2\e\int_{\mathcal{O}}F(\ur)\ur dx\geq 2\e\int_{\mathcal{O}}(-k_2+k_3|\ur|^{p+1})dx\\
&=-2k_2|\mathcal{O}|+2k_3\e\|\ur\|^{p+1}_{p+1}.
\end{split}
\end{equation}
On the other hand, by means of assumptions $(g_2)$ and $(h_2)$, we obtain
\begin{equation}\label{eq4-5}
\e\|g(r,\ur)\|^2_{\mathcal{L}_2(U;\hh)}\leq C_g+C_g\e|\ur|^2,
\end{equation}
and
\begin{equation}\label{eq4-6}
\int_{E_m}\e|h(u^m(r-,\tau-t,u_0),\xi)|^2\lambda(d\xi)\leq C_h+C_h\e|\ur|^2,
\end{equation}
separately. Substituting  \eqref{eq4-4}-\eqref{eq4-6} into \eqref{eq4-2}, ignoring the second term of the left-hand side of \eqref{eq4-2} and the second term of right-hand side of estimate to \eqref{eq4-4}, for every $m\in\mathbb{N}$, we find
\begin{equation*}
e^{(2\delta-L_1)r}\e|\ur|^2\leq e^{(2\delta-L_1)(\tau-t)}\e|u_0|^2+e^{(2\delta-L_1)r}\frac{2k_2|\mathcal{O}|+L_1}{2\delta-L_1}.
\end{equation*}
Therefore,  we infer  on the interval $(\tau-t,\tau)$ that
\begin{equation}\label{eq4-7}
\e|u^m(\tau,\tau-t,u_0)|^2\leq e^{-(2\delta-L_1)t}\e|u_0|^2+\frac{2k_2|\mathcal{O}|+L_1}{2\delta-L_1},\qquad \forall m\in\mathbb{N}.
\end{equation}

{\bf Step 2.} Let us proceed likewise Step 3 in Theorem \ref{thm3-2}. For $m_1$, $m_2\in\mathbb{N}$ with $m_2<m_1$, then we have $\lambda(E_{m_2})<\lambda(E_{m_1})$. Define $\mathbf{R}_{m_1,m_2}(r,\tau-t,u_0)=u^{m_1}(r,\tau-t,u_0)-u^{m_2}(r,\tau-t,u_0)$, similar arguments  as \eqref{eq3-18}-\eqref{eq3-25} imply (replace $\mathbf{B}^{m_1,m_2}$ and $t\in[\tau,T]$ by $\mathbf{R}_{m_1,m_2}$ and $s\in[\tau-t,r]$, respectively)  for every $\tau\in\mathbb{R}$, $t>\tau$ and $r>\tau-t$ that, 
\begin{equation*}
\e\left[\sup_{s\in[\tau-t,r]}|\mathbf{R}_{m_1,m_2}(s)|^2\right]\to0,\ \ \e\int_{\tau-t}^r\|\mathbf{R}_{m_1,m_2}(s)\|^2ds\to 0\qquad \mbox{as}~~m_1,~m_2\rightarrow\infty.
\end{equation*}
 Therefore, for every $r>\tau-t$, there exists an $\mathcal{F}_t$-adapted process $u\in L^2(\Omega;\mathbf{D}([\tau-t,r];\hh))\cap L^2(\Omega;L^2(\tau-t,r;\vv))$ {(thanks to the uniqueness of solution, this limit is denoted by the same $u$}) such that
 \begin{equation*}
\lim_{m\rightarrow \infty}\e\left[\sup_{s\in[\tau-t,r]}|u^m(s)-u(s)|^2\right]=\lim_{m\rightarrow \infty}\e\int_{\tau-t}^r\|u^m(s)-u(s)\|^2ds=0,
\end{equation*}
 which, together with \eqref{eq4-7}, yields
  \begin{equation}\label{eq4-8}
\e|u(\tau,\tau-t,u_0)|^2\leq e^{-(2\delta-L_1)t}\e|u_0|^2+\frac{2k_2|\mathcal{O}|+L_1}{2\delta-L_1}.
\end{equation} 
Since $u_0\in B(\tau-t)$ and $B\in\mathcal{D}$, one has
\begin{equation*}
e^{-(2\delta-L_1)t}\e|u_0|^2\leq e^{-(2\delta-L_1)t}\|B(\tau-t)\|^2\rightarrow 0\qquad\mbox{as}\quad t\rightarrow \infty,
\end{equation*}
which along with \eqref{eq4-8} concludes the proof. 
\end{proof} 

We will present now the existence of weakly compact $\mathcal{D}$-pullback absorbing sets to problem \eqref{eq2-13}-\eqref{eq2-14}.

\begin{lemma}\label{lem4-2}
Under assumptions of Lemma \ref{lem4-1}, the mean random dynamical system $\Phi$ related to \eqref{eq2-13}-\eqref{eq2-14} has a weakly compact $\mathcal{D}$-pullback absorbing set $K=\{K(\tau):\tau\in\mathbb{R}\}\in\mathcal{D}$, which is given, for each $\tau\in\mathbb{R}$, by
\begin{equation}\label{eq4-9}
K(\tau)=\{u\in L^2(\Omega,\mathcal{F}_\tau;\hh):\e|u(\tau)|^2\leq R\},\ \mbox{where}\ R:=\frac{2k_2|\mathcal{O}|+L_1}{2\delta-L_1}+1.
\end{equation}
\end{lemma}
\begin{proof} For  each $\tau\in\mathbb{R}$, it is easy to see that $K(\tau)$  given by \eqref{eq4-9} is a bounded closed convex subset of $L^2(\Omega,\mathcal{F}_\tau;\hh)$. Therefore, it is a weakly compact subset of $L^2(\Omega,\mathcal{F}_\tau;\hh)$. Moreover, it follows from Lemma~\ref{lem4-1} that, for every $\tau\in\mathbb{R}$ and {$B=\{B(\tau-t)\}\in\mathcal{D}$}, there exists $T=T(\tau,B)>0$ such that $\Phi(t,\tau-t,B(\tau-t))\subset K(\tau),$ for all $t\geq T$. 
On the other hand, 
\begin{equation*}
\lim_{\tau\rightarrow-\infty}e^{(2\delta-L_1)\tau}\|K(\tau)\|^2_{L^2(\Omega,\mathcal{F}_\tau;\hh)}\leq\lim_{\tau\rightarrow-\infty}e^{(2\delta-L_1)\tau}\left(\frac{2k_2|\mathcal{O}|+L_1}{2\delta-L_1}+1\right)=0.
\end{equation*}
Hence, we have verified $K(\tau)\in\mathcal{D}$. Namely, $K$ is a weakly compact $\mathcal{D}$-pullback absorbing set for $\Phi$. 
\end{proof}

Now, we are in a  position to address the existence and uniqueness  of weak $\mathcal{D}$-pullback mean random attractors to problem \eqref{eq2-13}-\eqref{eq2-14} (see \cite{W3} for the definition of this kind of attractors).

\begin{theorem}\label{thm4-3}
Assume the conditions of Lemma \ref{lem4-1} hold. Then the mean random dynamical system $\Phi$ to problem \eqref{eq2-13}-\eqref{eq2-14} has a unique weak $\mathcal{D}$-pullback mean random attractor $\mathcal{A}=\{\mathcal{A}(\tau):\tau\in\mathbb{R}\}\in\mathcal{D}$ in $L^2(\Omega,\mathcal{F};\hh)$ over $(\Omega,\mathcal{F},\{\mathcal{F}_t\}_{t\in\mathbb{R}},\mathbb{P})$.
\end{theorem}
\begin{proof} The existence and uniqueness of the weak $\mathcal{D}$-pullback mean random attractors $\mathcal{A}\in\mathcal{D}$ of $\Phi$ are  immediate consequences of  Lemma \ref{lem4-2} and \cite[Theorem 2.13]{W3}. 
\end{proof}

\section{Invariant measures and ergodicity}\label{s5}
In this section, we establish the existence  of invariant measures and ergodicity to the following autonomous problem, for $x\in\mathcal{O}$ and $t>0$,
\begin{equation}\label{eq5-1}
du(t)+(-\Delta)^{\gamma}u(t)dt+F(u(t))dt+\delta u(t)dt
=g(u(t))dW(t)+\int_{E}h(u(t-),\xi)\tilde{N}(dt,d\xi),
\end{equation}
with boundary and initial conditions, 
\begin{equation}\label{eq5-2}
u(t,x)=0,\quad x\in\mathbb{R}^d\backslash \mathcal{O},~~t\geq 0, \qquad \mbox{and}\qquad  u(0,x)=u_0,\quad x\in\mathcal{O},
\end{equation}
respectively, where $\delta$ is a positive constant as in \eqref{eq2-2}, $W$ and $\tilde{N}$ are independent real-valued standard Wiener process and compensated Poisson random measure, respectively. 

In the remaining of this section, we will assume the nonlinear functions $g:\mathbb{R}\rightarrow \mathbb{R}$, $h:\mathbb{R}\times E\rightarrow \mathbb{R}$ are globally Lipschitz continuous with linearly growing rate. More precisely, there exist positive constants $\alpha_1$, $\alpha_2$, $\beta_1$ and $\beta_2$  such that, for every $s,r\in\mathbb{R}$,
\begin{equation}\label{eq5-3}
|g(s)-g(r)|^2\leq \alpha_1|s-r|^2,\qquad |g(s)|^2\leq \beta_1+\frac{\delta}{2}|s|^2,
\end{equation}
and for $\lambda(E^\prime)<\infty$ with $E^\prime \subset E$,  suppose 
\begin{equation}\label{eq5-4}
\int_{E^\prime}|h(s,\xi)-h(r,\xi)|^2\lambda(d\xi)\leq \alpha_2|s-r|^2,\qquad 
\int_{E^\prime}|h(s,\xi)|^2\lambda(d\xi)\leq \beta_2+\frac{\delta}{2}|s|^2.
\end{equation}
Under assumptions {$(F_1)$-$(F_3)$, \eqref{eq5-3}-\eqref{eq5-4},}  Theorem \ref{thm3-2} shows that for every $\mathcal{F}_0$-measurable initial value $u_0$ in $L^2(\Omega;\hh)$, problem 
\eqref{eq5-1}-\eqref{eq5-2} possesses a unique c\`adl\`ag $\hh$-valued $\mathcal{F}_t$-adapted solution $u(t,u_0)$ in the sense of Definition \ref{def3-1}.

\subsection{Introduction to invariant measures and ergodicity}
We first provide the definitions of invariant measures and ergodicity (for more details, see \cite{D2, M1,M2} and the references therein).  Let $\mathbb{X}$ be a Polish space,
$P_t$ and $P(t,x,\Gamma)$, $t\geq 0$, $x\in \mathbb{X}$, $\Gamma\in\mathcal{B}(\mathbb{X})$, be the corresponding transition semigroup and transition function.
\begin{definition}\label{def5-1}
A probability measure $\mu$ on $(\mathbb{X},\mathcal{B}(\mathbb{X}))$ is said  to be an invariant measure or a stationary measure for a given transition probability function $P(t,x,dy)$ if it satisfies 
$$\mu(A)=\int_{\mathbb{X}}P(t,x,A)\mu(dx),\ \forall A\in\mathcal{B}(\mathbb{X}),\forall t>0.$$
Equivalently, if for all $\varphi\in C_b(\mathbb{X})$ {(the space of bounded and continuous  Borel  functions on $\mathbb{X}$)} and  $t\geq 0$,
$$ \int_{\mathbb{X}}\varphi(x)\mu (dx)=\int_{\mathbb{X}}(P_t\varphi)(x)\mu(dx),$$
where the Markov semigroup $(P_t)_{t\geq 0}$ is defined by
$$ P_t\varphi(x)=\int_{\mathbb{X}}\varphi(y)P(t,x,dy).$$
\end{definition}

\begin{definition}
Let $\mu$ be an invariant measure for $(P_t)_{t\geq 0}$. We say that the measure $\mu$ is an ergodic measure, if for all $\varphi\in L^2(\mathbb{X};\mu)$, 
$$ \lim_{T\rightarrow \infty}\frac{1}{T}\int_0^T(P_t\varphi)(x)dt=\int_{\mathbb{X}}\varphi(x)\mu(dx),~~\mbox{in}~ L^2(\mathbb{X};\mu).$$
\end{definition}

The next lemma is crucial to prove the existence of ergodicity.
\begin{lemma}\cite[Theorem 3.2.6]{D2}\label{lem5-3}
If $\mu\in\mathcal{M}(\mathbb{X})$ (the set of all probability measures defined on $(\mathbb{X},\mathcal{B}(\mathbb{X})$) is the unique invariant measure  for the semigroup $P_t (t>0)$, then it is ergodic.
\end{lemma}

Now, we will introduce the transition operators to problem \eqref{eq5-1}-\eqref{eq5-2}. By the construction of the solutions, we find that the transition operators are homogeneous. Therefore,
let $u(t,u_0)$ be the unique solution to \eqref{eq5-1}-\eqref{eq5-2},
then for  any $\varphi\in C_b(\hh)$,  $ t\geq 0$ and $ u_0\in \hh$, the corresponding Markov transition operator $P_{t}$ can be defined as (\cite[Theorem 9.8, p.244]{D1}),
\begin{equation}
(P_{t}\varphi)(u_0)=\e[\varphi(u(t,u_0))].
\end{equation}
 As usual, for $\Gamma\in\mathcal{B}(\hh)$, $t\geq 0$ and $u_0\in\hh$, we set
$ P(t,u_0,\Gamma)=(P_t \chi_{\Gamma})(u_0)=\mathbb{P}\{\omega\in\Omega:u(t,u_0)\in\Gamma\},$
where $\chi_{\Gamma}$ denotes the characteristic (or indicator) function of $\Gamma$. Then $P(t,u_0,\cdot)$ is the probability distribution of $u(t,u_0)$.
 In addition,  Theorem \ref{thm3-2} (cf. \eqref{eq3-32}) proved the solution to \eqref{eq5-1}-\eqref{eq5-2} depends continuously on initial value, which implies the Feller property of $P_t$ for $t\geq 0$. Similarly arguments as in \cite{W1} show that the solution $u(t,u_0)$ of problem \eqref{eq5-1}-\eqref{eq5-2} is an $\hh$-valued Markov process, which implies  that if $\varphi:\hh\rightarrow \mathbb{R}$ is a bounded Borel function, then $ (P_t\varphi)(u_0)=(P_r(P_{t-r}\varphi))(u_0),$ for all $u_0\in\hh,$ for any $0\leq r\leq t$, $\mathbb{P}$-a.s.

\subsection{Existence of invariant measures}

We will prove the existence of invariant measures to problem \eqref{eq5-1}-\eqref{eq5-2}   by using the Krylov-Bogolyubov method. To this end, we first derive uniform estimates on the solutions via the following lemma.
 
\begin{lemma}\label{lem5-4}
Assume $(F_1)$-$(F_3)$ and \eqref{eq5-3}-\eqref{eq5-4} hold.  Then, the solution $u(t,u_0)$ of \eqref{eq5-1}-\eqref{eq5-2} satisfies
\begin{equation*}
\e|u(t)|^2+C(d,\gamma)\int_0^te^{\delta(s-t)}\e\|u(s)\|^2_{\dot{\mathbb{H}}^{\gamma}(\mathbb{R}^d)}ds\leq e^{-\delta t}\e|u_0|^2+\frac{2k_2|\mathcal{O}|+\beta_1+\beta_2}{\delta}, \ \mbox{for all}\ t\geq 0,
\end{equation*}
and 
\begin{equation*}
\frac{\eta_1}{t}\e\int_0^t\|u(s)\|^2ds\leq  \frac{1}{T_0}\e|u_0|^2+{2k_2|\mathcal{O}|+\beta_1+\beta_2},\ \mbox{for all}\  t>T_0>0, \ \mbox{where}\ \eta_1=\min\{C(d,\gamma),\delta\}.
\end{equation*}
\end{lemma}
\begin{proof}   As $\lambda$ is $\sigma$-finite  on the Polish space $E$, there exist measurable subsets $E_m\nearrow E$ satisfying $\lambda(E_m)<\infty$ for all $m\in\mathbb{N}$ and $\cup_{m=1}^\infty E_m=E$. Now, on account of Remark \ref{remark}, for each $m\in\mathbb{N}$,
similar to \eqref{eq4-2},  applying It\^o's formula to $e^{\delta t}\e|u^m(t)|^2$, we obtain
\begin{equation}\label{eq5-5}
\begin{split}
&e^{\delta t}\e|u^m(t)|^2+C(d,\gamma)\int_0^te^{\delta s}\e\|u^m(s)\|^2_{\dot{\mathbb{H}}^{\gamma}(\mathbb{R}^d)}ds+2\delta\int_0^te^{\delta s}\e|u^m(s)|^2ds\\
&~~+2\int_0^te^{\delta s}\e<F(u^m(s)),u^m(s)>ds\\
&=\e|u_0|^2+\delta \int_0^te^{\delta s}\e|u^m(s)|^2ds+\int_0^te^{\delta s}\e|g(u^m(s))|^2ds+\int_0^t\int_{E_m}e^{\delta s}\e|h(u^m(s^-),\xi)|^2\lambda(d\xi)ds.
\end{split}
\end{equation}
Making use of the same estimates as in \eqref{eq4-4}, we have
\begin{equation}\label{eq5-6}
-2\e<F(u^m(s)),u^m(s)>\leq 2k_2|\mathcal{O}|-2k_3\e\|u^m(s)\|^{p+1}_{p+1}.
\end{equation}
By means of \eqref{eq5-3}-\eqref{eq5-4}, we find
\begin{equation}\label{eq5-7}
\e|g(u^m(s))|^2\leq \beta_1+\frac{\delta}{2}\e|u^m(s)|^2,\quad\mbox{and}\quad \int_{E_m}\e|h(u^m(s-),\xi)|^2\lambda(d\xi)\leq \beta_2+\frac{\delta}{2}\e|u^m(s)|^2.
\end{equation}
It follows from \eqref{eq5-5}-\eqref{eq5-7} that
\begin{equation}\label{eq5-8}
\begin{split}
e^{\delta t}\e|u^m(t)|^2&+C(d,\gamma)\int_0^te^{\delta s}\e\|u^m(s)\|^2_{\dot{\mathbb{H}}^{\gamma}(\mathbb{R}^d)}ds
\leq\e|u_0|^2+\frac{(2k_2|\mathcal{O}|+\beta_1+\beta_2)}{\delta}e^{\delta t}.
\end{split}
\end{equation}
Thus, 
\begin{equation}\label{eq5-9}
\e|u^m(t)|^2+C(d,\gamma)\int_0^te^{\delta(s-t)}\e\|u^m(s)\|^2_{\dot{\mathbb{H}}^{\gamma}(\mathbb{R}^d)}ds\leq e^{-\delta t}\e|u_0|^2+\frac{2k_2|\mathcal{O}|+\beta_1+\beta_2}{\delta}.
\end{equation}
On the other hand, using It\^o's formula to $|u^m(t)|^2$, we derive
\begin{equation*}
\begin{split}
\e|u^m(t)|^2&+C(d,\gamma)\int_0^t\e\|u^m(s)\|^2_{\dot{\mathbb{H}}^{\gamma}(\mathbb{R}^d)}ds+2\delta\int_0^t\e|u^m(s)|^2ds+2\int_0^te^{\delta s}\e<F(u^m(s)),u^m(s)>ds\\
&\leq \e|u_0|^2+\int_0^t\e|g(u^m(s))|^2ds
+\int_0^t\int_{E_m}\e|h(u^m(s^-),\xi)|^2\lambda(d\xi)ds.
\end{split}
\end{equation*}
Let $\eta_1  =\min\{C(d,\gamma),\delta\}$, by means of \eqref{eq5-6}-\eqref{eq5-7}, we obtain
\begin{equation*}
\begin{split}
\e|u^m(t)|^2&+\eta_1\int_0^t\e\|u^m(s)\|^2ds\leq \e|u_0|^2+(2k_2|\mathcal{O}|+\beta_1+\beta_2)t,\quad\forall t>0,
\end{split}
\end{equation*}
which implies, 
\begin{equation}\label{eq5-10}
\frac{\eta_1}{t}\e\int_0^t\|u^m(s)\|^2ds\leq \frac{\e|u_0|^2}{T_0}+2k_2|\mathcal{O}|+\beta_1+\beta_2,\qquad \forall t>T_0>0.
\end{equation}
Since $\{u^m\}_{m=1}^{\infty}$ converges to $u$ in $L^2(\Omega;\mathbf{D}([0,T];\hh))\cap L^2(\Omega;L^2(0,T;\vv))$, then $u$ satisfies energy estimates \eqref{eq5-5} (see also Step 3 of Theorem \ref{thm3-2}), which, together with \eqref{eq5-9}-\eqref{eq5-10}, concludes the proof of this lemma.
\end{proof}

\begin{theorem}\label{thm5-5}
Under  conditions of Lemma \ref{lem5-4}, problem \eqref{eq5-1}-\eqref{eq5-2} has an invariant measure on $\hh$.
\end{theorem}
\begin{proof}  Using the Chebyshev inequality and Lemma \ref{lem5-4}, we infer for  $T_0>0$ and $R>0$,
\begin{equation}\label{eq5-12}
\begin{split}
\sup_{t\geq T_0}\frac{1}{t}\int_0^t\mathbb{P}(\{\| u(s,u_0)\|>R\})ds&\leq \sup_{t\geq T_0}\frac{1}{tR^2}\int_0^t\e\|u(s,u_0)\|^2ds\leq \frac{\e|u_0|^2}{T_0R^2\eta_1}
+\frac{2k_2|\mathcal{O}|+\beta_1+\beta_2}{ R^2 \eta_1}.
\end{split}
\end{equation}
The above inequality implies, for all $t>T_0$ and  every $\varepsilon >0$, there exists $R_0:=\sqrt{\frac{\e|u_0|^2}{T_0\eta_1\varepsilon}+\frac{(2k_2|\mathcal{O}|+\beta_1+\beta_2)}{\eta_1\varepsilon}}>0$ such that, for any $R\geq R_0$, 
\begin{equation}\label{eq5-13}
\begin{split}
\mu_{t,u_0}(\Gamma)&:=\frac{1}{t}\int_0^tP(s,u_0,\Gamma)ds=\frac{1}{t}\int_0^t
\mathbb{P}(\{\omega\in\Omega: u(s,u_0)\in\Gamma\})ds\\
&\geq \frac{1}{t}\int_0^t\mathbb{P}(\{\omega\in\Omega: \|u(s,u_0)\|\leq R_0\})ds\\
& =1-\frac{1}{t}\int_0^t\mathbb{P}(\{\omega\in\Omega: \|u(s,u_0)\|>R_0\})ds= 1-\varepsilon,
\qquad \Gamma: =B(0,R),
\end{split}
\end{equation} 
where $B(0,R)$ is the ball centered at $0$ with radius $R$ in $\vv$. Since $\vv$ is compactly embedded  in $\hh$ $(\vv\hookrightarrow\hookrightarrow\hh)$, \eqref{eq5-13} shows for every $\varepsilon>0$,
there exists a compact set $\mathcal{K}\in\hh$ such that $\mu_{t,u_0}(\mathcal{K})>1-\varepsilon$ for all $t>T_0$. Hence, the sequence of probability measure $\mu_{t,u_0}$ is tight on $\hh$. 

As a result,
an application of the Krylov-Bogoliubov theorem (see \cite{C3}) shows that there exists a sequence $t_n\rightarrow \infty$ as $n\rightarrow \infty$ such that $\mu_{t_n,u_0}\rightarrow \mu $ weakly as $n\rightarrow \infty$. Moreover, $\mu$ is an invariant measure for this transition operator $P_t$, defined by
$ (P_t\varphi)(u_0)=\e[\varphi(u(t,u_0))],$ for all $\varphi\in C_b(\hh).$ 
Thus, the proof of this theorem is complete. 
\end{proof} 

\subsection{Ergodicity}

 We are now interested in the ergodicity of problem \eqref{eq5-1}-\eqref{eq5-2}. Lemma \ref{lem5-3} states that if $\mu$ is the unique invariant measure for $P_t$, then it is ergodic. Thus, in what follows, we will prove the invariant measure $\mu$ in Theorem \ref{thm5-5} is unique. To this end, the following lemma is needed.
\begin{lemma}\label{lem5-6}
Under  assumptions of Lemma \ref{lem5-4},  additionally, suppose $2\delta>2k_4+\alpha_1+\alpha_2$. Let $u$ and $v$ be two solutions of problem \eqref{eq5-1}-\eqref{eq5-2} with initial data $u_0$, $v_0\in L^2(\Omega;\hh)$, respectively. Then, we have
\begin{equation*}
\e|u(t)-v(t)|^2\leq \e|u_0-v_0|^2 e^{-(2\delta-(2k_4+\alpha_1+\alpha_2))t},\qquad \forall t\geq 0.
\end{equation*}
\end{lemma}
\begin{proof}   Denote by $z=u-v$,   $z_0=u_0-v_0$, $z^m=u^m-v^m$. According to Remark \ref{remark}, by  similar computations to \eqref{eq3-26}, applying   It\^o's formula to $e^{(2\delta-(2k_4+\alpha_1+\alpha_2))t}|z^m(t)|^2$ and making use of condition $(F_3)$, we obtain
\begin{equation}\label{eq5-14}
\begin{split}
&e^{(2\delta-(2k_4+\alpha_1+\alpha_2))t}|z^m(t)|^2+2\delta\int_0^te^{(2\delta-(2k_4+\alpha_1+\alpha_2))s}|z^m(s)|^2ds\\
&\leq |z_0|^2+(2\delta -\alpha_1-\alpha_2)\int_{0}^te^{(2\delta-(2k_4+\alpha_1+\alpha_2))s}|z^m(s)|^2ds\\
 &~~+2\int_{0}^te^{(2\delta-(2k_4+\alpha_1+\alpha_2))s}(z^m(s),(g(u^m(s))-g(v^m(s)))dW(s)\\
 &~~+\int_0^te^{(2\delta-(2k_4+\alpha_1+\alpha_2))s}|g(u^m(s))-g(v^m(s))|^2ds\\
 &~~+2\int_0^t\int_{E_m}e^{(2\delta-(2k_4+\alpha_1+\alpha_2))s}(z^m(s),h(u^m(s-),\xi)-h(v^m(s-),\xi))\tilde{N}(ds,d\xi)\\
 &~~+\int_0^t\int_{E_m}e^{(2\delta-(2k_4+\alpha_1+\alpha_2))s}|h(u^m(s-),\xi)-h(v^m(s-),\xi)|^2N(ds,d\xi).
\end{split}
\end{equation}
Taking expectation on both sides of \eqref{eq5-14} and thanks to  \eqref{eq5-3}-\eqref{eq5-4}, we have
\begin{equation*}
\e|z^m(t)|^2\leq \e|z_0|^2 e^{-(2\delta-(2k_4+\alpha_1+\alpha_2))t},\qquad \forall t\geq 0.
\end{equation*}
Since the sequences $\{u^m\}_{m=1}^{\infty}$ and $\{v^m\}_{m=1}^{\infty}$ are converging  in $L^2(\Omega;\mathbf{D}([0,T];\hh))\cap L^2(\Omega;L^2(0,T;\vv))$, so is  $\{z^m\}_{m=1}^{\infty}$. By taking limit on both sides of \eqref{eq5-14}, we can conclude the proof of this lemma. 
\end{proof} 

Let us now establish the uniqueness of invariant measure to \eqref{eq5-1}-\eqref{eq5-2} which ensures its ergodicity.

\begin{theorem}
Suppose the conditions of Lemma \ref{lem5-6} hold true. Then, for any initial value $u_0\in L^2(\Omega;\hh)$, there exists a unique invariant measure $\mu$ to problem \eqref{eq5-1}-\eqref{eq5-2}. In addition, this measure $\mu$ is ergodic. 
\end{theorem}
\begin{proof}   Assume there is another invariant measure $\tilde{\mu}$ for transition operator $(P_t)_{t\geq 0}$. Then, for every $\varphi\in \mbox{Lip}(\hh)$ ($\varphi$ is a Lipschitz function with Lipschitz constant $L_{\varphi}$) and initial data $u_0$, $v_0\in L^2(\Omega;\hh)$. By means of Definition \ref{def5-1} and
 Lemma \ref{lem5-6}, we have  
\begin{equation*}
\begin{split}
&\left| \int_{\hh}\varphi(u_0) \mu(d u_0)-\int_{\hh}\varphi(v_0)\tilde{\mu}(dv_0)\right|
=\left|\int_{\hh}(P_t\varphi)(u_0)\mu(du_0)-\int_{\hh}(P_t\varphi)(v_0)\tilde{\mu}(dv_0)\right|\\
&=\left|\int_{\hh}\int_{\hh}[(P_t\varphi)(u_0)-(P_t\varphi)(v_0)]\mu(du_0)\tilde{\mu}(dv_0)\right|=\left|\int_{\hh}\int_{\hh}\e[\varphi(u(t,u_0))-\e\varphi(v(t,v_0))]
\mu(d u_0)\tilde{\mu}(dv_0)\right|\\
&\leq L_{\varphi} e^{-\frac{(2\delta-(2k_4+\alpha_1+\alpha_2))t}{2}}
\int_{\hh}\int_{\hh}\e|u_0-v_0|\mu(du_0)\tilde{\mu}(dv_0)\rightarrow 0~~~\mbox{as}~~t\rightarrow \infty.
\end{split}
\end{equation*}
Since $\mu$ is the unique invariant measure for transition operator $(P_t)_{t\geq 0}$, by the density of Lip$(\mathbb{H})$ in $C_b(\hh)$, we know  $\mu$ is ergodic.  
\end{proof}

\section{Large deviation Principle}\label{s6}

In this section,  we will consider the following stochastic perturbations of the fractional partial differential equation, where $\varepsilon$ is a small parameter,
\begin{equation}\label{eq6-1}
\begin{split}
d\ue(t)&+(-\Delta)^{\gamma}\ue(t)dt+F(\ue(t))dt+\delta \ue(t)dt\\
&=\sqrt{\varepsilon}g(t,\ue(t))dW(t)+\varepsilon\int_{E}{h(t,\ue(t-),\xi)}\tilde{N}^{\varepsilon^{-1}}(dt,d\xi),\quad x\in\mathcal{O},~t>0,\\
\end{split}
\end{equation}
with boundary and initial conditions
\begin{equation}\label{eq6-2}
\ue(t,x)=0,\quad x\in\mathbb{R}^d\backslash \mathcal{O},~~t>0, \qquad \mbox{and}\qquad  \ \ue(0,x)=u_0(x),\quad x\in\mathcal{O},
\end{equation}
respectively. {Assume $\mathcal{O}$ is a bounded domain in $\mathbb{R}^d$ with smooth boundary satisfying  {$2\gamma<d$,    and $F$ satisfies conditions $(F_1)$-$(F_3)$ with $p+1\in (2,\frac{2d}{d-2\gamma}]$} such that $\vv$ is continuously embedded in $L^{p+1}(\mathcal{O})$ (see \cite[Theorem 6.7]{N1})}. Here $E$ is a locally compact Polish space, $W$ is a cylindrical Brownian motion in $U$, $N^{\varepsilon^{-1}}$ is a Poisson random measure on $[0,T]\times E$ with a $\sigma$-finite intensity measure $\varepsilon^{-1}L_{T}\otimes\lambda$, $L_{T}$ is the Lebesgue measure on $[0,T]$ and $\lambda$  a $\sigma$-finite measure on $E$.
 $\tilde{N}^{\varepsilon^{-1}}([0,t]\times O)=N^{\varepsilon^{-1}}([0,t]\times O)-\varepsilon^{-1}t\lambda(O)$, $\forall O\in\mathcal{B}(E)$ with $\lambda(O)<\infty$, is the compensated Poisson random measure. 
We emphasize that in this section,  $\mathbf{D}([0,T];E)$ denotes the space of right continuous functions with left limits from $[0,T]$ to $E$ endowed with the Skorokhod topology.

Let $\{\ue,\varepsilon>0\}$ be a family of random variables defined on a probability space $(\Omega,\mathcal{F},\mathbb{P})$ taking values in a  Polish space $\mathcal{E}$.  The large deviation principle of problem \eqref{eq6-1}-\eqref{eq6-2} is concerned with the exponential decay of $\mathbb{P}(\ue\in O)$ as $\varepsilon\rightarrow 0$.

\subsection{Controlled Poisson random measure}

Let $E$ be a locally compact Polish space. Denote by $\mathcal{M}_{FC}(E)$  the space of all measures $\lambda$ on $(E,\mathcal{B}(E))$ such that $\lambda(K)<\infty$ for every compact $K$ in $E$. Endow $\mathcal{M}_{FC}(E)$ with the weakest topology such that for every $f\in C_c(E)$ (the space of continuous functions with compact support), the function $\lambda\rightarrow 
\langle f,\lambda\rangle=\int_{E}
f(u)d\lambda( u) $ $(\lambda \in\mathcal{M}_{FC}(E)) $ is continuous.
This topology can be metrized such that $\mathcal{M}_{FC}(E)$ is a Polish space (see, for example, \cite{B3}). Fix $T\in(0,\infty)$ and let $E_T=[0,T]\times E$. Fix a measure $\lambda\in\mathcal{M}_{FC}(E)$, we denote by $\lambda_T=L_T\otimes \lambda$. 

We recall that a Poisson random measure $\mathbf{n}$ on $E_T$ with intensity measure $\lambda_T$ is an $\mathcal{M}_{FC}(E_T)$ valued random variable, such that for each $B\in\mathcal{B}(E_T)$ with $\lambda_T(B)<\infty$, $\mathbf{n}(B)$ is Poisson distributed with mean $\lambda_T(B)$ and for disjoint $B_1,\cdots,B_k\in\mathcal{B}(E_T)$, $\mathbf{n}(B_1),\cdots,\mathbf{n}(B_k)$ are mutually independent random variables (see, for example, \cite{I1}). Denote by $\mathbb{P}$ the measure induced by $\mathbf{n}$ on $(\mathcal{M}_{FC}(E_T),\mathcal{B}(\mathcal{M}_{FC}(E_T)))$.
Then letting $\mathbb{M}=\mathcal{M}_{FC}(E_T)$, $\mathbb{P}$ is the unique probability measure on $(\mathbb{M},\mathcal{B}(\mathbb{M}))$ under which the canonical map, $N:\mathbb{M}\rightarrow \mathbb{M}$, $N(m)\doteq m$, is a Poisson random measure with intensity measure $\lambda_T$. With applications to large deviations in mind, we also consider, for $\theta>0$, the probability measure $\mathbb{P}_\theta$ on $(\mathbb{M},\mathcal{B}(\mathbb{M}))$ under which $N$ is a Poisson random measure with intensity $\theta\lambda_T$. The corresponding expectation operators will be denoted by $\mathbb{E}$ and $\mathbb{E}_\theta$, respectively.

Set $F=E\times [0,\infty)$ and $F_T=[0,T]\times F$. Similarly, let $\bar{\mathbb{M}}=\mathcal{M}_{FC}(F_T)$ and let $\bar{\mathbb{P}}$ be the unique probability measure on $(\bar{\mathbb{M}},\mathcal{B}(\bar{\mathbb{M}})))$ under which the canonical map, $\bar{N}:\bar{\mathbb{M}}\rightarrow \bar{\mathbb{M}}$, $\bar{N}(m):= m$, is a Poisson random measure with intensity measure $\bar{\lambda}_T=L_T\otimes \lambda\otimes L_{\infty}$, with $L_\infty$ being Lebesgue measure on $[0,\infty)$. The corresponding expectation operator will be denoted by $\bar{\mathbb{E}}$. Let $\mathcal{F}_t:= \sigma\{\bar{N}((0,s]\times O):0\leq s\leq t, O\in\mathcal{B}(F)\}$ and  denote by $\bar{\mathcal{F}}_t$ the completion under $\bar{\mathbb{P}}$.  Set $\bar{\mathcal{P}}$ the predictable $\sigma$-field on $[0,T]\times \bar{\mathbb{M}}$ with the filtration $\{\bar{\mathcal{F}}_t:0\leq t\leq T\}$ on $(\bar{\mathbb{M}},\mathcal{B}(\bar{\mathbb{M}}))$.  Let $\bar{\mathcal{A}}$ be the class of all $(\bar{\mathcal{P}}\times\mathcal{B}({E}))/\mathcal{B}[0,\infty)$-measurable maps $\varphi: E_T\times \bar{\mathbb{M}}\rightarrow [0,\infty)$. For $\varphi\in\bar{\mathcal{A}}$, define a continuous process $N^{\varphi}$ on $E_T$ by
\begin{equation}\label{eq6-3}
N^{\varphi}((0,t]\times {K})(\cdot)=\int_{(0,t]\times {K}}\int_{(0,\infty)} 1_{[0,\varphi(s,x,\cdot)]}(r)\bar{N}(dsdxdr),\qquad t\in[0,T], ~{K}\in\mathcal{B}(E).
\end{equation}
$N^{\varphi}$ is the controlled random measure, with $\varphi$ selecting the intensity for the points at location $x$ and time $s$, in a possibly random but nonanticipating way. When $\varphi(s,x,\bar{m})\equiv\theta\in(0,\infty)$, we write $N^{\varphi}=N^{\theta}$. Note that $N^{\theta}$ has the same distribution with respect to $\bar{\mathbb{P}}$ as $N$ has with respect to $\mathbb{P}_\theta$.

\subsection{Poisson random measure and Brownian motion}
Set {$\mathbb{W}=C([0,T];U)$}, $\mathbb{U}=\mathbb{W}\times \mathbb{M}$ and $\bar{\mathbb{U}}=\mathbb{W}\times\bar{\mathbb{M}}$. Then, for $(w,m)\in\mathbb{U}$, let the mapping $N^{\mathbb{U}}:\mathbb{U}\rightarrow \mathbb{M}$ be defined by $N^{\mathbb{U}}(w,m)=m$, and let ${W^{\mathbb{U}}}:\mathbb{U}\rightarrow C([0,T];U)$ be defined by
$W^{\mathbb{U}}(w,m)(t)=\sum_{i=1}^{\infty}(w(t),a_i)a_i$, {recalling that the sequence $\{a_i\}_{i=1}^{\infty}$ is an orthonormal basis of the separable Hilbert space $U$.} The mappings $\bar{N}^{\bar{\mathbb{U}}}:\bar{\mathbb{U}}\rightarrow \bar{\mathbb{M}}$ and $\bar{W}^{\bar{\mathbb{U}}}:\bar{\mathbb{U}}\rightarrow C([0,T];U)$ are defined analogously. 
For every $t\in[0,T]$, define the $\sigma$-filtration 
$$\mathcal{G}^{\mathbb{U}}_t:=\sigma\left(\left\{(W(s), N^{\mathbb{U}}((0,s]\times O)): 0\leq s\leq t, O\in\mathcal{B}(E)\right\}\right).$$
For every $\theta>0$, {for a given $\lambda \in \mathcal{M}_{FC}(E)$}, it follows from \cite[Sec. I.8]{I1} that there exists a unique probability measure  $\mathbb{P}_\theta^{\mathbb{U}}$ on $(\mathbb{U},\mathcal{B}(\mathbb{U}))$ such that:\\
$~~~(i)$ $W$ is a cylindrical Brownian motion in $U$;\\
$~~(ii)$ $N^{\mathbb{U}}$ is a Poisson random measure  with intensity measure $\theta \lambda_T$;\\
$(iii)$ $W$ and $N$ are independent.\\
Analogously, we define $(\bar{\mathbb{P}}^{\bar{\mathbb{U}}}_\theta,\bar{\mathcal{G}}_t^{\bar{\mathbb{U}}})$ and denote $\bar{\mathbb{P}}_{\theta=1}^{\bar{\mathbb{U}}}$ by $\bar{\mathbb{P}}^{\bar{\mathbb{U}}}$. We denote by $\{\bar{\mathcal{F}}_t^{\bar{\mathbb{U}}}\}$ the $\bar{\mathbb{P}}^{\bar{\mathbb{U}}}$-completion of $\{\bar{\mathcal{G}}_t^{\bar{\mathbb{U}}}\}$ and $\bar{\mathcal{P}}^{\bar{\mathbb{U}}}$ the predictable $\sigma$-field on $[0,T]\times\bar{\mathbb{U}}$ with the filtration $\{\bar{\mathcal{F}}_t^{\bar{\mathbb{U}}}\}$ on $(\bar{\mathbb{U}},\mathcal{B}(\bar{\mathbb{U}}))$.  Let $\mathbb{A}_1$ be the class of all $(\bar{\mathcal{P}}^{\bar{\mathbb{U}}}\otimes \mathcal{B}(E))/\mathcal{B}[0,\infty)$-measurable maps $\varphi: E_T\times \bar{\mathbb{U}}\rightarrow [0,\infty)$. 

On the one hand, define $l:[0,\infty)\rightarrow [0,\infty)$ by
$ l(r)=r\log r-r+1,\ r\in[0,\infty).$
Then, for any $\varphi\in\mathbb{A}_1$, the quantity
\begin{equation}\label{eq6-4}
\tilde{L}_T(\varphi)(\omega)=\int_{E_T}l(\varphi(t,x,\omega))\lambda_T(dtdx),\qquad \omega \in\bar{\mathbb{U}},
\end{equation}
is well defined as a $[0,\infty]$-valued random variable.
On the other hand, define
\begin{equation}\label{eq6-5}
\mathbb{A}_2:=\left\{\psi: \psi~\mbox{is}~ \bar{\mathcal{P}}^{\bar{\mathbb{U}}}/\mathcal{B}(U)~\mbox{measurable and}~\int_0^T\|\psi(t)\|_U^2dt<\infty, ~\mbox{a.s.}\mbox{-}\bar{\mathbb{P}}^{\bar{\mathbb{U}}}\right\}.
\end{equation}
Set $\mathcal{U}=\mathbb{A}_2\times \mathbb{A}_1$. Define $\hat{L}_T(\psi):=\frac{1}{2}\int_0^T\|\psi(t)\|_U^2dt$ for $\psi\in\mathbb{A}_2$, and $L_T(u):=\hat{L}_T(\psi)+\tilde{L}_T(\varphi)$ for $u=(\psi,\varphi)\in\mathcal{U}$.

\subsection{A general criterium}
In this subsection, we recall a general criterium for a large deviation principle established in \cite{B3}. Let $\{\mathcal{G}^{\varepsilon}\}_{\varepsilon>0}$ be a family of measurable maps from $\bar{\mathbb{U}}$ to $\Xi$, where $\bar{\mathbb{U}}$ is introduced in Section 6.2 and $\Xi$ is some Polish space. We present below a sufficient condition for large deviation principle (LDP in abbreviation) to hold for the family $Z^{\varepsilon}=\mathcal{G}^{\varepsilon}(\sqrt{\varepsilon}W,\varepsilon N^{\varepsilon^{-1}})$ as $\varepsilon \rightarrow 0$.

Define for each $\Upsilon\in\mathbb{N}$,
$ \tilde{\mathrm{S}}^{\Upsilon}=\{\rho:E_T\rightarrow [0,\infty): \tilde{L}_T(\rho)\leq \Upsilon\},$
and 
$\hat{\mathrm{S}}^{\Upsilon}=\{\sigma\in L^2(0,T;U): \hat{L}_T(\sigma)\leq \Upsilon\}.$ 
A function $\rho\in \tilde{\mathrm{S}}^{\Upsilon}$ can be identified with a measure $\lambda_T^\rho\in\mathbb{M}$, defined by
$$ \lambda_T^\rho(O)=\int_{O}\rho(s,\xi)\lambda_T(dsd\xi),\qquad O\in\mathcal{B}(E_T).$$
{This identification induces a topology on $\tilde{\mathrm{S}}^{\Upsilon}$ under which $\tilde{\mathrm{S}}^{\Upsilon}$ is a compact space,} see  Appendix of \cite{B2}. Throughout this paper we use this topology on $\tilde{\mathrm{S}}^{\Upsilon}$. Set $\mathrm{S}^{\Upsilon}=\hat{\mathrm{S}}^{\Upsilon}\times \tilde{\mathrm{S}}^{\Upsilon}$.  Define $\mathbb{S}=\cup_{\Upsilon\geq 1}\mathrm{S}^{\Upsilon}$, and let 
$\mathcal{U}^{\Upsilon}=\{u=(\psi,\varphi)\in\mathcal{U}: u(\omega)\in\mathrm{S}^{\Upsilon},~\bar{\mathbb{P}}^{\bar{\mathbb{U}}}\mbox{-a.e.}~ \omega\}.$
The following condition will be sufficient to establish an LDP for a family $\{Z^{\varepsilon}\}_{\varepsilon>0}$ defined by $Z^{\varepsilon}=\mathcal{G}^{\varepsilon}(\sqrt{\varepsilon}W,\varepsilon N^{\varepsilon^{-1}})$.

\begin{condition}\label{con6-1}
There exists a measurable map $\mathcal{G}^0:{\mathbb{U}}\rightarrow \Xi $ such that the following hold.
\begin{enumerate}
\item[$(a)$] For any $\Upsilon\in\mathbb{N}$, let $(\sigma_n,\rho_n)$, $(\sigma, \rho)\in \mathrm{S}^{\Upsilon}$ be such that $(\sigma_n,\rho_n)\rightarrow (\sigma,\rho)$ as $n\rightarrow \infty$. Then,
$$\mathcal{G}^0\left(\int_0^{\cdot}\sigma_n(s)ds, \lambda_T^{\rho_n}\right)\rightarrow \mathcal{G}^0\left(\int_0^{\cdot}\sigma(s)ds,\lambda_T^\rho\right),\qquad \mbox{in}~~\Xi.$$
\item[$(b)$] For any $\Upsilon\in\mathbb{N}$, let $u_\varepsilon=(\psi_\varepsilon,\varphi_\varepsilon)$, $u=(\psi,\varphi)\in\mathcal{U}^{\Upsilon}$ such that  $u_{\varepsilon}$ converges in distribution to $u$ as $\varepsilon\rightarrow 0$. Then, 
$$\mathcal{G}^{\varepsilon}\left(\sqrt{\varepsilon}W+\int_0^{\cdot}\psi_\varepsilon(s)ds,\varepsilon N^{\varepsilon^{-1}\varphi_\varepsilon}\right)\Rightarrow \mathcal{G}^0\left(\int_0^{\cdot}\psi(s)ds, \lambda_T^{\varphi}\right), \qquad \mbox{in}~~\Xi.$$
\end{enumerate}
\end{condition}
We use the symbol ``$\Rightarrow$" to denote the convergence in distribution.

For $\phi\in\Xi$, define $\mathbb{S}_\phi=\{(\sigma,\rho)\in\mathbb{S}: \phi=\mathcal{G}^0(\int_0^{\cdot}\sigma(s)ds, \lambda_T^\rho)\}$. Let $I:\Xi \rightarrow [0,\infty]$ be defined by
\begin{equation}\label{eq6-6}
I(\phi)=\inf_{\pi=(\sigma,\rho)\in \mathbb{S}_\phi}\{L_T(\pi)\},\qquad \phi\in \Xi.
\end{equation}
By convention, $I(\phi)=\infty$ if $\mathbb{S}_\phi=\emptyset$.

The following criterium was established in \cite[Theorem 4.2]{B3}.
\begin{theorem}\label{thm6-2}
For $\varepsilon>0$, let $Z^{\varepsilon}$ be defined by $Z^{\varepsilon}=\mathcal{G}^{\varepsilon}(\sqrt{\varepsilon}W, \varepsilon N^{\varepsilon^{-1}})$ and suppose that condition \ref{con6-1} holds.
Then, $I$ defined as in \eqref{eq6-6} is a rate function on $\Xi$ and the family $\{Z^\varepsilon\}_{\varepsilon>0}$ satisfies an LDP with rate function $I$.
\end{theorem}

For applications, the following strengthened form of Theorem \ref{thm6-2} is useful. Let $\{E_m\subset E,m=1,2,\cdots\}$ be an increasing sequence of compact sets such that $\cup_{m=1}^{\infty}E_m=E$. For each $m$, let 
\begin{equation*}
\begin{split}
\mathbb{A}_1^{b,m}&:=\left\{\varphi\in\mathbb{A}_1:~\mbox{for all}~(t,\omega)\in[0,T]\times\bar{\mathbb{M}},\right.\\
&\quad\left. 1/m\leq \varphi(t,x,\omega)\leq m~\mbox{if}~x\in E_m~\mbox{and}~\varphi(t,x,\omega)=1~\mbox{if}~x\in E_m^c\right\},
\end{split}
\end{equation*}
and let $\mathbb{A}_1^b=\cup_{m=1}^{\infty}\mathbb{A}_1^{b,m}$. Define
$\tilde{\mathcal{U}}^{\Upsilon}=\mathcal{U}^{\Upsilon}\cap\{(\psi,\varphi):\varphi\in\mathbb{A}_1^b\}$.

\begin{theorem}\label{thm6-3}
Suppose Condition \ref{con6-1} holds with $\mathcal{U}^{\Upsilon}$ replaced by $\tilde{\mathcal{U}}^{\Upsilon}$. Then,   Theorem \ref{thm6-2} holds true.
\end{theorem}

\subsection{Hypotheses}

In addition to assumptions $(F_1)$-$(F_3)$ and $(g_1)$-$(g_3)$ stated above, we impose the following conditions on the jump noise term.
Let {$ h:[0,T]\times\hh\times E\rightarrow \hh$} be a measurable map.
\begin{condition}\label{con6-4}
For the locally compact Polish space $E$, there exist  $L^{\prime}_h>0$ and $C_h^{\prime}>0$, such that
\begin{enumerate}
\item[$(i)$] $\int_{E}|h(t,u_1,\xi)-h(t,u_2,\xi)|^2\lambda(d\xi)\leq L_h^\prime|u_1-u_2|^2,$ for all $t\in[0,T]$ and $u_1,u_2\in\hh$;
\item [$(ii)$] $\int_{E}|h(t,u,\xi)|^2\lambda(d\xi)\leq C^\prime_h(1+|u|^2),$ for all $t\in[0,T]$ and $u\in\hh$.
\end{enumerate}
\end{condition}

Define
\begin{equation*}
\|h(t,\xi)\|_{0,\hh}=\sup_{u\in\hh}\frac{|h(t,u,\xi)|}{1+|u|},\qquad (t,\xi)\in[0,T]\times E,
\end{equation*}
and
\begin{equation*}
\|h(t,\xi)\|_{1,\hh}=\sup_{u_1,u_2\in\hh,u_1\neq u_2}\frac{|h(t,u_1,\xi)-h(t,u_2,\xi)|}{|u_1-u_2|},\qquad  (t,\xi)\in[0,T]\times E.
\end{equation*}
\begin{condition}({\bf Exponential integrability})\label{con6-5}
For $i=0,1$, there exists $\delta_1^i>0$ such that for all $\mathbb{X}\in\mathcal{B}([0,T]\times E)$ satisfying $\lambda_T(\mathbb{X})<\infty$, it follows that
$ \int_{\mathbb{X}}e^{\delta_1^i\|h(s,\xi)\|^2_{i,\hh}}\lambda(d\xi)ds<\infty.$
\end{condition}
\begin{remark}\label{rem6-6}\cite[Remark 3.1]{Z1}
 Condition \ref{con6-5} implies that, for every $\delta_2^i>0$ ($i=0,1$) and for all $\mathbb{X}\in\mathcal{B}([0,T]\times E)$ satisfying $\lambda_T(\mathbb{X})<\infty$, we have 
 $ \int_{\mathbb{X}}e^{\delta_2^i\|h(s,\xi)\|_{i,\hh}}\lambda(d\xi)ds<\infty. $
\end{remark}

\begin{lemma}\label{lem6-7}\cite[Lemma 3.1]{Z1}
Under Conditions \ref{con6-4} and \ref{con6-5}:
\begin{enumerate}
\item[$(i)$] For $i=0,1$ and every $\Upsilon\in\mathbb{N}$,
\begin{equation}\label{eq6-8}
C_{i,1}^{\Upsilon}:=\sup_{\rho\in\tilde{\mathrm{S}}^{\Upsilon}}\int_{E_T}\|h(s,\xi)\|_{i,\hh}|\rho(s,\xi)-1|\lambda(d\xi)ds<\infty,
\end{equation}
\begin{equation}\label{eq6-7}
C_{i,2}^{\Upsilon}:=\sup_{\rho\in\tilde{\mathrm{S}}^{\Upsilon}}\int_{E_T}\|h(s,\xi)\|^2_{i,\hh}(\rho(s,\xi)+1)\lambda(d\xi)ds<\infty;
\end{equation}
\item[$(ii)$] For every $\eta>0$, there exists $\delta>0$ such that for $A\subset[0,T]$ satisfying $\lambda_T(A)<\delta$,
\begin{equation}\label{eq6-9}
\sup_{\rho\in\tilde{\mathrm{S}}^{\Upsilon}}\int_A\int_E\|h(s,\xi)\|_{i,\hh}|\rho(s,\xi)-1|\lambda(d\xi)ds\leq \eta.
\end{equation}
\end{enumerate}
\end{lemma}
\begin{lemma}\label{lem6-8}\cite[Lemma 3.2]{Z1}
(i) If $\sup_{t\in[0,T]}|Y(t)|<\infty$, for any $\pi=(\sigma,\rho)\in\mathbb{S}$, then
$$ g(\cdot, Y(\cdot))\sigma(\cdot)\in L^1(0,T;\hh),\quad \int_Eh(\cdot,Y(\cdot),\xi)(\rho(\cdot,\xi)-1)\lambda(d\xi)\in L^1(0,T;\hh);$$
(ii) If the family of mappings $\{Y_n:[0,T]\rightarrow \hh,n\geq 1\}$ satisfies 
$C=\sup_{n\in\mathbb{N}}\sup_{t\in[0,T]}|Y_n(t)|<\infty$, then
$$C_{\Upsilon}:=\sup_{\pi=(\sigma,\rho)\in \mathrm{S}^{\Upsilon}}\sup_{n\in\mathbb{N}}\left[
\int_0^T\left | \int_E h(s,Y_n(s),\xi)(\rho(s,\xi)-1)\lambda(d\xi)\right|ds+\int_0^T|g(s,Y_n(s))\sigma(s)|ds\right]<\infty.
$$
\end{lemma}

We also need the following lemma, the proof of which can be found in \cite[Lemma 3.11]{B2}.
\begin{lemma}\label{lem6-9}
Let $d:[0,T]\times E\rightarrow \mathbb{R}$ be a measurable function such that
$ \int_{E_T}|d(s,\xi)|^2\lambda(d\xi)ds<\infty,$
and for all $\delta_3\in(0,\infty)$, $\mathbb{X}\in\mathcal{B}([0,T]\times E)$ satisfying $\lambda_T(\mathbb{X})<\infty$, it follows
$ \int_{\mathbb{X}}\exp(\delta_3|d(s,\xi)|)\lambda(d\xi)ds<\infty.$
\begin{enumerate}
\item[$(i)$] Fix $\Upsilon\in\mathbb{N}$ and let $\rho_n,\rho\in\tilde{\mathrm{S}}^{\Upsilon}$ be such that $\rho_n\rightarrow \rho$ as $n\rightarrow \infty$. Then
$$ \lim_{n\rightarrow \infty}\int_{E_T}d(s,\xi)(\rho_n(s,\xi)-1)\lambda(d\xi)ds=\int_{E_T}d(s,\xi)(\rho(s,\xi)-1)\lambda(d\xi)ds;$$
\item[$(ii)$] Fix $\Upsilon\in\mathbb{N}$. Given $\varepsilon>0$, there exists a compact set $K_{\varepsilon}\subset E$, such that
$$ \sup_{\rho\in\tilde{\mathrm{S}}^{\Upsilon}}\int_0^T\int_{K_\varepsilon^c}|d(s,\xi)||\rho(s,\xi)-1|\lambda(d\xi)ds\leq \varepsilon;$$
\item[$(iii)$] For every compact $K\subset E$,
$ \lim_{M\rightarrow\infty}\sup_{\rho\in\tilde{\mathrm{S}}^{\Upsilon}}\int_0^T\int_K|d(s,\xi)|1_{\{d\geq M\}}\rho(s,\xi)\lambda(d\xi)ds=0.$
\end{enumerate}
\end{lemma}

The previous lemmas will be used together with the following compactness result, which represents a variant of the criterium for compactness stated in \cite[Section 5, Chapter I]{L2} and \cite[Section 13.3]{T1}, to prove main results  of this section. 
Given $p>1$, $\alpha\in (0,1)$, let $W^{\alpha,p}(0,T;\hh)$ be the Sobolev space of all $u\in L^p(0,T;\hh)$ such that,
$\int_0^T\int_0^T\frac{|u(t)-u(s)|^p}{|t-s|^{1+\alpha p}}dtds<\infty,$
endowed with the norm
$\|u\|^p_{W^{\alpha,p}(0,T;\hh)}=\int_0^T|u(t)|^pdt+\int_0^T\int_0^T\frac{|u(t)-u(s)|^p}{|t-s|^{1+\alpha p}}dtds.$

\begin{lemma}\label{lemsoblev}{\cite[Theorem 2.1]{F1}}
Let $B_0\subset B \subset B_1$ be Banach spaces, $B_0$ and $B_1$ reflexive with compact embedding of $B_0$ in $B$. Let $p\in(1,\infty)$ and $\alpha\in (0,1)$ be given. Let $X$ be the space 
$ X=L^p(0,T;B_0)\cap W^{\alpha,p}(0,T;B_1),$
endowed with the natural norm. Then, the embedding of $X$ into $L^p(0,T;B)$ is compact.
\end{lemma}

\subsection{Main results}

In this section, 
assume $u_0$ is deterministic. Let $\ue$ be the $\hh$-valued solution to \eqref{eq6-1}-\eqref{eq6-2} with initial value $u_0$. In the sequel, we will establish an LDP for $\{\ue\}$ as $\varepsilon \rightarrow 0$. We start with the following definition.

\begin{definition}\label{def6-10}
Let $(\bar{\mathbb{U}},\mathcal{B}(\bar{\mathbb{U}}),\bar{\mathbb{P}}^{\bar{\mathbb{U}}},\{\bar{\mathcal{F}}_t^{\bar{\mathbb{U}}}\}_{t\geq0})$ be a filtered probability space. Suppose  $u_0$ is an $\bar{\mathcal{F}}_0$-measurable $\hh$-valued random variable such that ${\bar{\mathbb{E}}}|u_0|^2<\infty$. A stochastic process $\{\ue(t)\}_{t\in[0,T]}$ defined on $\bar{\mathbb{U}}$ is said to be a $\hh$-valued solution to \eqref{eq6-1}-\eqref{eq6-2} with initial value $u_0$, if
\begin{enumerate}
\item[$\bullet$] $\ue(t)$ is a $\hh$-valued $\bar{\mathcal{F}}_t^{\bar{\mathbb{U}}}$-measurable random variable for all $t\in[0,T]$;
\item[$\bullet$] $\ue(t)\in L^2(\Omega;\mathbf{D}([0,T];\hh))\cap L^2(\Omega;L^2(0,T;\vv))\cap L^{p+1}(\Omega;L^{p+1}(0,T;L^{p+1}(\mathbb{R}^d)))$, $\bar{\mathbb{P}}^{\bar{\mathbb{U}}}$\mbox{-}a.s.;
\item[$\bullet$] For all $t\in[0,T]$,
\begin{equation}\label{eq6-10}
\begin{split}
\ue(t)&=u_0-\int_0^t(-\Delta)^{\gamma}\ue(s)ds-\delta\int_0^t\ue(s)ds
-\int_0^tF(\ue(s))ds\\
&~~+\sqrt{\varepsilon}\int_0^tg(s,\ue(s))dW(s)+\varepsilon\int_0^t\int_Eh(s,\ue(s-),\xi)\tilde{N}^{\varepsilon^{-1}}(ds,d\xi),
\end{split}
\end{equation}
in $L^2(0,T;\mathbb{V}^*)+ L^q(0,T;L^q(\mathbb{R}^d))$, where $q$ is the conjugate number of $p+1$.
\end{enumerate}
\end{definition}

\begin{definition}\label{def6-11}
The stochastic fractional partial differential equation \eqref{eq6-1}-\eqref{eq6-2} is said to satisfy the pathwise uniqueness   property,  if any two $\hh$-valued solutions $\ue_1$ and $\ue_2$, defined on the same filtered probability space, with respect to the same Poisson random measure and Brownian motion, starting from the same initial condition $u_0$, coincide almost surely.
\end{definition}

We begin by introducing the mapping $\mathcal{G}^0$ that will be used to define the rate function and  verify Condition \ref{con6-1}. Recall that $\mathbb{S}=\cup_{\Upsilon \geq 1}\mathrm{S}^{\Upsilon}$. As a first step, we show that under the conditions stated below, for every $\pi=(\sigma,\rho)\in\mathbb{S}$, the deterministic integral equation 
\begin{equation}\label{eq6-11}
\begin{split}
u^\pi(t)&=u_0-\int_0^t Au^\pi(s)ds-\int_0^t F(u^\pi(s))ds+\int_0^tg(s,u^\pi(s))\sigma(s)ds\\
&~~+\int_0^t\int_Eh(s,u^\pi(s),\xi)(\rho(s,\xi)-1)\lambda(d\xi)ds,
\end{split}
\end{equation}
has a unique continuous solution. Here $\pi=(\sigma,\rho)$ plays the role of a control. 

\begin{theorem}\label{thm6-12}
Let $u_0\in\hh$ and $\pi=(\sigma,\rho)\in\mathbb{S}$. Suppose  $(F_1)$-$(F_3)$, $(g_1)$-$(g_3)$, Conditions \ref{con6-4} and \ref{con6-5} hold. Then, there exists a unique $u^\pi\in C([0,T];\hh)\cap L^2(0,T;\vv)\cap L^{p+1}(0,T;L^{p+1}(\mathbb{R}^d))$, such that
\begin{equation}\label{eq6-12}
\begin{split}
u^\pi(t)&=u_0-\int_0^tA u^\pi(s)ds-\int_0^tF(u^\pi(s))ds+\int_0^tg(s,u^\pi(s))\sigma(s)ds\\
&~~+\int_0^t\int_Eh(s,u^\pi(s),\xi)(\rho(s,\xi)-1)\lambda(d\xi)ds,~~\mbox{in}~~{L^2(0,T;\mathbb{V}^*)+ L^q(0,T;L^q(\mathbb{R}^d))}.
\end{split}
\end{equation}
Moreover, for fixed $\Upsilon\in\mathbb{N}$, there exists a constant $C_{\Upsilon}>0$ (which depends on $\Upsilon$) such that
\begin{equation}\label{eq6-13}
\sup_{\pi\in \mathrm{S}^{\Upsilon}}\left(\sup_{t\in[0,T]}|u^\pi(t)|^2+\int_0^T\|u^\pi(t)\|^2dt\right)\leq C_\Upsilon.
\end{equation}
\end{theorem}
\begin{proof} {\it Existence of solutions.} Given $n\in\mathbb{N}$, similar to \cite{W2}, let $X_n$ be the space spanned by $\{e_j,j=1,2,\cdots,n\}$ and $P_n:\hh\rightarrow X_n$ be the projection given by 
$ P_nu^\pi=\sum_{j=1}^{n}(u^\pi,e_j)e_j,$ for all  $u^\pi\in\hh.$
We can extend $P_n$ to $\mathbb{V}^*$ and $(L^p(\mathbb{R}^d))^*$ by
$ P_n\phi=\sum_{j=1}^n(\phi,e_j)e_j,$ for all $ \phi\in\mathbb{V}^*~\mbox{or}~\phi\in(L^p(\mathbb{R}^d))^*.$
Consider the following Fadeo-Galerkin approximations: $u_n^\pi(t)\in X_n$ denotes the solution of
\begin{equation}\label{eq6-14}
\begin{split}
du^\pi_n(t)&=-Au_n^\pi(t)dt-P_nF(u_n^\pi(t))dt+P_ng(t,u_n^\pi(t))\sigma(t)dt\\[1.0ex]
&~~+\int_E P_nh(t,u_n^\pi(t),\xi)(\rho(t,\xi)-1)\lambda(d\xi), \qquad \forall t\in(0,T],
\end{split}
\end{equation}
with initial condition $u_n^\pi(0):=P_nu_0=u_{n,0}$. We will state the existence and uniqueness of solutions to problem \eqref{eq6-14} for each $n\in\mathbb{N}$. Let $v_0(t)=P_nu_0$ with $t\in[0,T]$. Suppose $v_k$ has been defined for $m-1\geq k\geq 1$. Define $v_m\in C([0,T];\hh)\cap L^2(0,T;\vv)\cap L^{p+1}(0,T;L^{p+1}(\mathbb{R}^d))$ as the unique solution of
\begin{equation}\label{eq6-15}
\begin{split}
dv_m(t)&=-Av_m(t)dt-P_nF(v_m(t))dt+P_ng(t,v_{m-1}(t))\sigma(t)dt\\[1.0ex]
&~~+\int_EP_nh(t,v_{m-1}(t),\xi)(\rho(t,\xi)-1)\lambda(d\xi),
\end{split}
\end{equation}
and $v_m(0)=P_nu_0$. By slightly modifying the proof of \cite[Theorem 3.1]{X1}, 
one can verify that the limit $u_n^\pi$ of $v_m$, as $m\rightarrow \infty$, is the unique solution of \eqref{eq6-14} satisfying $C([0,T];\hh)\cap L^2(0,T;\vv)\cap L^{p+1}(0,T;L^{p+1}(\mathbb{R}^d))$.

Next we will prove that there exists a constant $C>0$ depending on  $\Upsilon$, such that
\begin{equation}\label{eq6-16}
\sup_{n\geq 1}\left(\sup_{t\in[0,T]}|u_n^\pi(t)|^2+\int_0^T\|u_n^\pi(t)\|^2dt+\int_0^T\|u_n^\pi(t)\|_{p+1}^{p+1}dt\right)\leq C,
\end{equation}
and for $\alpha\in(0,1/2)$, there exists $C_\alpha>0$ depending on  $\Upsilon$ and $\alpha$, such that
\begin{equation}\label{eq6-17}
\sup_{n\geq 1}\|u_n^\pi\|^2_{W^{\alpha,q}(0,T;V^*)}\leq C_\alpha.
\end{equation}
Let us first show \eqref{eq6-16}. By means of energy estimates, proceeding likewise as in Theorem \ref{thm3-2}, let $\eta=\min\{\frac{C(d,\gamma)}{2},\delta\}$,  with the help of $(F_2)$,  we have 
\begin{equation}\label{eq6-18}
\begin{split}
&|u_n^\pi(t)|^2+2\eta\int_0^t\|u_n^\pi(s)\|^2ds+2k_3\int_0^t\|u_n^\pi(s)\|_{p+1}^{p+1}ds\leq |u_{n,0}|^2+2\int_0^t\left(P_ng(s,u_n^\pi(s))\sigma(s),u_n^\pi(s)\right)ds\\
&~~+2k_2|\mathcal{O}|t+2\int_0^t\left(P_n\int_Eh(s,u_n^\pi(s),\xi)(\rho(s,\xi)-1)\lambda(d\xi),u_n^\pi(s)\right)ds\\
&\leq |u_{n,0}|^2+2k_2|\mathcal{O}|t+2\int_0^t\|g(s,u_n^\pi(s))\|_{\mathcal{L}_2(U;\hh)}\|\sigma(s)\|_U|u^\pi_n(s)|ds\\
&~~+2\int_0^t\int_E|h(s,u_n^\pi(s),\xi)||\rho(s,\xi)-1||u_n^{\pi}(s)|\lambda(d\xi)ds,\qquad \forall t\in[0,T].
\end{split}
\end{equation}
On the one hand, it follows from assumption $(g_2)$ and inequality $x\sqrt{1+x^2}\leq 1+2x^2$ that
\begin{equation}\label{eq6-19}
\begin{split}
2\int_0^t\|g(s,u_n^\pi(s))\|_{\mathcal{L}_2(U;\hh)}\|\sigma(s)\|_U|u_n^\pi(s)|ds
&\leq 2\sqrt{C_g}\int_0^t\sqrt{1+|u_n^\pi(s)|^2}\|\sigma(s)\|_U|u_n^\pi(s)|ds\\
&\leq 2\sqrt{C_g}\int_0^t\|\sigma(s)\|_Uds+4\sqrt{C_g}\int_0^t|u^\pi_n(s)|^2\|\sigma(s)\|_Uds.
\end{split}
\end{equation}
On the other hand, 
\begin{equation}\label{eq6-20}
\begin{split}
&~\quad 2\int_0^t\int_E|h(s,u_n^\pi(s),\xi)||\rho(s,\xi)-1||u_n^\pi(s)|\lambda(d\xi)ds\\
&=2\int_0^t\int_E\frac{|h(s,u_n^\pi(s),\xi)|}{1+|u^\pi_n(s)|}(1+|u_n^\pi(s)|)|\rho(s,\xi)-1||u_n^\pi(s)|\lambda(d\xi)ds\\
&\leq 2\int_0^t\int_E\|h(s,\xi)\|_{0,\hh}|\rho(s,\xi)-1|(1+2|u_n^\pi(s)|^2)\lambda(d\xi)ds
\end{split}
\end{equation}
\begin{equation*}
\begin{split}
&\leq 2\int_0^t\int_E\|h(s,\xi)\|_{0,\hh}|\rho(s,\xi)-1|\lambda(d\xi)ds\\
&~~+4\int_0^t|u_n^\pi(s)|^2\left(\int_E\|h(s,\xi)\|_{0,\hh}|\rho(s,\xi)-1|\lambda(d\xi)\right)ds.
\end{split}
\end{equation*}
Hence, it follows from \eqref{eq6-18}-\eqref{eq6-20} that
\begin{equation}\label{eq6-21}
\begin{split}
& |u_n^\pi(t)|^2+2\eta\int_0^t\|u_n^\pi(s)\|^2ds+2k_3\int_0^t\|u_n^\pi(s)\|_{p+1}^{p+1}ds\\&
\leq|u_{n,0}|^2+2k_2|\mathcal{O}|t+2\sqrt{C_g}\int_0^t\|\sigma(s)\|_Uds
+2\int_0^t\int_E\|h(s,\xi)\|_{0,\hh}|\rho(s,\xi)-1|\lambda(d\xi)ds\\
&+4\int_0^t|u_n^\pi(s)|^2\left(\sqrt{C_g}\|\sigma(s)\|_U+\int_E\|h(s,\xi)\|_{0,\hh}|\rho(s,\xi)-1|\lambda(d\xi)\right)ds.
\end{split}
\end{equation}
Taking supremum with respect to $t\in[0,T]$, we have 
\begin{equation*}
\begin{split}
\sup_{t\in[0,T]}|u_n^\pi(t)|^2&\leq
|u_{n,0}|^2+2k_2|\mathcal{O}|T+2\sqrt{C_g}\int_0^T\|\sigma(s)\|_Uds
+2\int_0^T\int_E\|h(s,\xi)\|_{0,\hh}|\rho(s,\xi)-1|\lambda(d\xi)ds\\
&~~+4\int_0^T\sup_{t\in[0,s]}|u_n^\pi(t)|^2\left(\sqrt{C_g}\|\sigma(s)\|_U+\int_E\|h(s,\xi)\|_{0,\hh}|\rho(s,\xi)-1|\lambda(d\xi)\right)ds.
\end{split}
\end{equation*}
Using the fact that $\sigma \in L^2(0,T;U)$ and Lemma \ref{lem6-7},  applying the Gronwall Lemma to the above inequality, combining with \eqref{eq6-21}, we  can prove \eqref{eq6-16}.

Now we will check \eqref{eq6-17}. Notice that,
\begin{equation*}
\begin{split}
u^\pi_n(t)&=P_nu_0-\int_0^tAu^\pi_n(s)ds-\int_0^tP_nF(u^\pi_n(s))ds+\int_0^tP_ng(s,u_n^\pi(s))\sigma(s)ds\\
&~~+\int_0^t\int_EP_nh(s,u_n^\pi(s),\xi)(\rho(s,\xi)-1)\lambda(d\xi)ds\\[0.6ex]
&:=J_n^1+J_n^2(t)+J_n^3(t)+J_n^4(t)+J_n^5(t),\qquad \forall t\in[0,T].
\end{split}
\end{equation*}
Using the same arguments as \cite[Theorem 3.1]{L2},  there exists a positive constant $\mathcal{C}_1$ such that
\begin{equation}\label{eq6-22}
|J_n^1|^2=|u_{n,0}|^2\leq\mathcal{C}_1.
\end{equation}
For $J_n^2(t)$, thanks to the H\"older inequality,  we infer there exists a constant $\mathcal{C}_2>0$ such that
\begin{eqnarray}\label{eq6-23}
\|J_n^2\|_{W^{1,2}(0,T;\vv^*)}^2&=&\int_0^T\|J_n^2(t)\|^2_*dt+\int_0^T\left\|\frac{dJ_n^2(t)}{dt}\right\|^2_*dt=\int_0^T\left\|\int_0^tAu_n^\pi(s)ds\right\|^2_*dt+\int_0^T\|Au_n^\pi(t)\|_*^2dt\notag\\
&\leq& (1+T^2)\int_0^T\|Au^\pi_n(t)\|_*^2dt\leq\mathcal{C}_2.
\end{eqnarray}
For $J_n^3(t)$, on the one hand, taking into account of $(F_1)$, we find 
\begin{equation*}
\begin{split}
|F(u_n^\pi)|_q^q=\int_{\mathbb{R}^d}|F(u_n^\pi(x))|^qdx&=\int_{\mathcal{O}}|F(u^\pi_n(x))|^qdx\leq k_1^q\int_{\mathcal{O}}\left(1+|u_n^\pi(x)|^p\right)^qdx\\[0.8ex]
&\leq 2^{q-1}k_1^q|\mathcal{O}|+2^{q-1}k_1^q\|u^\pi_n\|^{p+1}_{p+1}.
\end{split}
\end{equation*}
Since $P_n:(L^{p+1}(\mathbb{R}^d))^*=L^q(\mathbb{R}^d)\rightarrow X_n\subset\hh$, we know ${P_nF(u_n^\pi)}\in \hh$ for a.a. $t\in[0,T]$.  
 The above estimate, the H\"older inequality and \eqref{eq6-16} imply that there exists a constant $\mathcal{C}_3>0$ such that
\begin{eqnarray}\label{eq6-24}
\|J_n^3\|^q_{W^{1,q}(0,T;\hh)}&=&\int_0^T|J_n^3(t)|^qdt+\int_0^T\left|\frac{dJ_n^3(t)}{dt}\right|^qdt=\int_0^T\left|\int_0^tP_nF(u^\pi_n(s))ds\right|^qdt
+\int_0^T|P_nF(u_n^\pi(t))|^qdt\notag\\
&\leq& \left(T^{\frac{1}{p}+1}+1\right)\int_0^T\|F(u_n^\pi(s))\|_q^qds\leq \mathcal{C}_3.
\end{eqnarray}
We will do estimate for $J_n^4$ now. For $0\leq s< t\leq T$, it follows from assumption $(g_2)$, the H\"older inequality and the fact $u_n^\pi\in C([0,T];\hh)$ that
\begin{equation*}
\begin{split}
|J_n^4(t)-J_n^4(s)|^2&=\left|\int_s^tP_ng(r,u_n^\pi(r))\sigma(r)dr\right|^2\leq \left(\int_s^t|P_ng(r,u^\pi_n(r))\sigma(r)|dr\right)^2\\
&\leq \left(\int_s^t\|g(r,u_n^\pi(r))\|_{\mathcal{L}_2(U;\hh)}\|\sigma(r)\|_Udr\right)^2\leq \left(\int_s^t\sqrt{C_g(1+|u_n^\pi(r)|^2)}\|\sigma(r)\|_Udr\right)^2\\
&\leq C_g(t-s)\left(1+\sup_{t\in[0,T]}|u^\pi_n(t)|^2\right)\int_s^t\|\sigma(r)\|^2_Udr.
\end{split}
\end{equation*}
By means of the above estimate and the H\"older inequality,  for $\alpha\in (0,1/2)$, there exists  $\mathcal{C}_4>0$ such that
\begin{eqnarray}\label{eq6-25}
 &&\|J_n^4\|^2_{W^{\alpha,2}(0,T;\hh)}=\int_0^T|J_n^4(t)|^2dt+\int_0^T\int_0^T\frac{|J_n^4(t)-J_n^4(s)|^2}{|t-s|^{1+2\alpha}}dsdt\notag\\
&&\leq\int_0^T\left|\int_0^tP_ng(s,u_n^\pi(s))\sigma(s)ds\right|^2dt+C_g\left(1+\sup_{t\in[0,T]}|u^\pi_n(t)|^2\right)\int_0^T\int_0^T\int_s^t
\frac{\|\sigma(r)\|_U^2}{|t-s|^{2\alpha}}drdsdt\notag\\
&&\leq C_g\left(1+\sup_{t\in[0,T]}|u^\pi_n(t)|^2\right)\left(T^2+\frac{T^{2-2\alpha}}{(1-2\alpha)(2-2\alpha)}\int_0^T\|\sigma(t)\|_U^2dt\right)\leq \mathcal{C}_4.
\end{eqnarray}
For $J_n^5$, with the help of Lemma \ref{lem6-7} and the fact that $u_n^{\pi}\in C([0,T];\hh)$,  for $0\leq s<t\leq T$, we derive
\begin{equation*}
\begin{split}
|J_n^5(t)-J_n^5(s)|^2&=\left|\int_s^tP_n\int_Eh(r,u^\pi_n(r),\xi)(\rho(r,\xi)-1)\lambda(d\xi)dr\right|^2\leq \left(\int_s^t\int_E|h(r,u^\pi_n(r),\xi)||\rho(r,\xi)-1|\lambda(d\xi)dr\right)^2\\
&\leq \left(\int_s^t\int_E\|h(r,\xi)\|_{0,\hh}(1+|u_n^\pi(r)|)|\rho(r,\xi)-1|\lambda(d\xi)dr\right)^2\\
&\leq C_{0,1}^\Upsilon\left(1+\sup_{t\in[0,T]}|u_n^\pi(t)|\right)^2\int_0^T\int_E
\|h(r,\xi)\|_{0,\hh}|\rho(r,\xi)-1|\lambda(d\xi)dr.
\end{split}
\end{equation*}
{Using the similar arguments as for $J_n^4$ (cf. \eqref{eq6-25}), we deduce that there exists  $\mathcal{C}_5>0$ such that}
\begin{equation}\label{eq6-26}
\|J_n^5\|^2_{W^{\alpha,2}(0,T;\hh)}\leq \mathcal{C}_5.
\end{equation}
Moreover, since $ \vv\subset L^{p+1}(\mathcal{O})\subset \hh:=\hh\subset L^q(\mathcal{O})\subset\vv^*$. Immediately,   we conclude \eqref{eq6-17} by \eqref{eq6-22}-\eqref{eq6-26}.

The estimates \eqref{eq6-16}-\eqref{eq6-17} ensure  the existence of an element $u^\pi\in L^2(0,T;\vv)\cap L^{\infty}(0,T;\hh)\cap L^{p+1}(0,T;L^{p+1}(\mathbb{R}^d))$, and a subsequence $u^\pi_{n^\prime}$ such that, as $n^\prime \rightarrow\infty$,
\begin{equation}\label{eq6-27}
\left\{
\begin{array}{rcl}
\begin{aligned}
 &u^\pi_{n^\prime} \rightarrow u^\pi ~~\mbox{weak-star in}~~  L^{\infty}(0,T;\hh); \\
 &u^\pi_{n^\prime}\rightarrow u^\pi~~ \mbox{weakly in}~~ L^{2}(0,T;\vv);\\
 &u^\pi_{n^\prime}\rightarrow u^\pi~~ \mbox{weakly in}~~ L^{p+1}(0,T;L^{p+1}(\mathbb{R}^{d}));\\
 &u^\pi_{n^\prime}\rightarrow u^\pi~~ \mbox{strongly in}~~ L^{q}(0,T;\hh);\\
  &F(u^\pi_{n^\prime})\rightarrow F(u^\pi)~~ \mbox{weakly in}~~ L^{q}(0,T;L^q(\mathbb{R}^d)),
 \end{aligned}
\end{array}
\right.
\end{equation}
where the  strong convergence holds thanks to Lemma \ref{lemsoblev}, and the last weak convergence follows from  the same arguments  as \cite[Theorem 2.7]{X2}, respectively. 

Next, we will show $u^\pi$ is the solution of \eqref{eq6-11}.
Let $\psi$ be a continuously differentiable function on $[0,T]$ with $\psi(T)=0$. For each fixed $n\in\mathbb{N}$, we multiply \eqref{eq6-14} by $\psi(t)e_j$ and then integrate by parts. This leads to the following equation,
\begin{equation*}
\begin{split}
&-\int_0^T(u_{n^{\prime}}^\pi(t),e_j)\psi^\prime(t)dt+\int_0^Ta(u_{n^{\prime}}^\pi(t),e_j)\psi(t)dt+\int_0^T(F(u_{n^{\prime}}^\pi(t)),e_j)\psi(t)dt
=(u_{n^{\prime},0},e_j)\psi(0)\\
&+\int_0^T(g(t,u_{n^{\prime}}^\pi(t))\sigma(t),e_j)\psi(t)dt+\int_0^T\left(\int_Eh(t,u_{n^{\prime}}^\pi(t),\xi)(\rho(t,\xi)-1)\lambda(d\xi),e_j\right)\psi(t)dt.
\end{split}
\end{equation*}
Taking limit when $n^\prime\rightarrow\infty$ and using \eqref{eq6-27}, we deduce 
\begin{eqnarray}\label{eq6-28}
&&\hskip-.5cm\lim_{n^\prime\rightarrow\infty}\left[-\int_0^T(u^\pi_{n^\prime}(t),e_j)\psi^{\prime}(t)dt+\int_0^Ta(u^\pi_{n^\prime}(t),e_j)\psi(t)dt\right.+\left.\int_0^T(F(u^\pi_{n^{\prime}}(t)),e_j)\psi(t)dt-(u_{n^{\prime},0},e_j)\psi(t)\right]\notag\\
&&\qquad=-\int_0^T(u^\pi(t),e_j)\psi^{\prime}(t)dt+\int_0^Ta(u^\pi(t),e_j)\psi(t)dt+\int_0^T(F(u^\pi(t)),e_j)\psi(t)dt-(u_0,e_j)\psi(t).
\end{eqnarray}
Hence, we only need to check
\begin{equation}\label{eq6-29}
\lim_{n^\prime\rightarrow\infty}\int_0^T\left|g(t,u^\pi_{n^{\prime}}(t))\sigma(t)-g(t,u^\pi(t))\sigma(t)\right|dt=0,
\end{equation}
and
\begin{equation}\label{eq6-30}
\lim_{n^\prime\rightarrow\infty}\int_0^T\int_E\left|h(t,u^\pi_{n^\prime}(t),\xi)(\rho(t,\xi)-1)-h(t,u^\pi(t),\xi)(\rho(t,\xi)-1)\right|\lambda(d\xi)dt=0.
\end{equation}
On the one hand, for every $\varepsilon>0$, let $A_{n^{\prime},\varepsilon}=\{t\in[0,T]:|u^\pi_{n^\prime}(t)-u^\pi(t)|>\varepsilon\}$, then by \eqref{eq6-27} and the Chebyshev inequality, we have
\begin{equation}\label{eq6-31}
\lim_{n^\prime\rightarrow\infty}L_T(A_{n^{\prime},\varepsilon})\leq\lim_{n^\prime\rightarrow\infty}\frac{\int_0^T|u^\pi_{n^\prime}(t)-u^\pi(t)|^2dt}{\varepsilon^2}=0.
\end{equation}
Consider $M=\sup_{i\in\mathbb{N}}\sup_{t\in[0,T]}|u^\pi_i(t)|\vee \sup_{t\in[0,T]}|u^\pi(t)|<\infty$, this assertion holds true thanks to \eqref{eq6-27}. Thus, due to assumption $(g_1)$ and the H\"older inequality,  we derive
\begin{equation*}
\begin{split}
&\int_0^T|g(t,u^\pi_{n^\prime}(t))\sigma(t)-g(t,u^\pi(t))\sigma(t)|dt\leq \int_0^T\|g(t,u^\pi_{n^\prime}(t))-g(t,u^\pi(t))\|_{\mathcal{L}_2(U;\hh)}\|\sigma(t)\|_Udt\\
&\leq\sqrt{L_g}\int_0^T|u^\pi_{n^\prime}(t)-u^\pi(t)|\|\sigma(t)\|_Udt\leq 2M\sqrt{L_g}\int_{A_{n^{\prime},\varepsilon}}\|\sigma(t)\|_Udt+\varepsilon\sqrt{L_g}\int_{A_{n^{\prime},\varepsilon}^c}\|\sigma(t)\|_Udt\\
&\leq 2M\sqrt{L_T(A_{n^{\prime},\varepsilon})}\sqrt{L_g}\left(\int_{A_{n^{\prime},\varepsilon}}\|\sigma(t)\|_U^2dt\right)^{1/2}+\varepsilon\sqrt{L_gT}\left(\int_{A_{n^{\prime},\varepsilon}^c}\|\sigma(t)\|^2_Udt\right)^{1/2}.
\end{split}
\end{equation*}
Thanks to \eqref{eq6-31}, the fact $\sigma\in\hat{\mathrm{S}}^\Upsilon$ and the above inequality, \eqref{eq6-29} holds.
On the other hand, since
\begin{eqnarray*}
&&\int_0^T\int_E|h(t,u^\pi_{n^\prime}(t),\xi)-h(t,u^\pi(t),\xi)||\rho(t,\xi)-1|\lambda(d\xi)dt\\
&&=\int_0^T\int_E\frac{|h(t,u^\pi_{n^\prime}(t),\xi)-h(t,u^\pi(t),\xi)|}{|u^\pi_{n^\prime}(t)-u^\pi(t)|}|u^\pi_{n^\prime}(t)-u^\pi(t)||\rho(t,\xi)-1|\lambda(d\xi)dt\\
&&\leq \int_0^T\int_E\|h(t,\xi)\|_{1,\hh}|u^\pi_{n^\prime}(t)-u^\pi(t)||\rho(t,\xi)-1|\lambda(d\xi)dt\\
&&\leq 2M\int_{A_{n^{\prime},\varepsilon}}\int_E\|h(t,\xi)\|_{1,\hh}|\rho(t,\xi)-1|\lambda(d\xi)dt+\varepsilon\int_{A_{n^{\prime},\varepsilon}^c}\int_E\|h(t,\xi)\|_{1,\hh}|\rho(t,\xi)-1|\lambda(d\xi)dt.
\end{eqnarray*}
Taking into account \eqref{eq6-31} and Lemma \ref{lem6-7}, together with the above inequality, \eqref{eq6-30} is also proved.
Therefore, it follows from \eqref{eq6-28}-\eqref{eq6-30} that, when $n^\prime\rightarrow\infty$,
\begin{equation}\label{eqlimit}
\begin{split}
&-\int_0^T(u^\pi(t),e_j)\psi^\prime(t)dt+\int_0^Ta(u^\pi(t),e_j)\psi(t)dt+\int_0^T(F(u^\pi(t)),e_j)\psi(t)dt=(u_0,e_j)\psi(0)\\
&+\int_0^T(g(t,u^\pi(t))\sigma(t),e_j)\psi(t)dt+\int_0^T\left(\int_Eh(t,u^\pi(t),\xi)(\rho(t,\xi)-1)\lambda(d\xi),e_j\right)\psi(t)dt,
\end{split}
\end{equation}
which implies $u^\pi$ is  solution of \eqref{eq6-11}. Moreover, by means of Lemma \ref{lem6-8} and using the same arguments as in the proof of \cite[Theorem 2.3]{W2}, we also obtain
$$ \frac{du^\pi(t)}{dt}\in L^2(0,T;\vv^*)+L^q(0,T;L^q(\mathbb{R}^d))+L^1(0,T;\hh).$$
Hence, it follows from {\cite[Lemma 1.2]{L2}} that $u^\pi\in C([0,T];\hh)$ and 
\begin{equation*}
\frac{1}{2}\frac{d}{dt}|u^\pi(t)|^2=\left(\frac{du^\pi(t)}{dt},u^\pi(t)\right)_{(\vv^*+L^q(\mathbb{R}^d),\vv\cap L^{p+1}(\mathbb{R}^d))}.
\end{equation*}
At last,  as $u^\pi$ is solution of \eqref{eq6-11}, \eqref{eq6-27} and \eqref{eq6-16} imply \eqref{eq6-13} holds. 

{\it Uniqueness of solution.} Eventually, we show that $u^\pi$ is the unique solution of \eqref{eq6-11}. To this end, assume that $u_1^\pi$ and $u_2^\pi$ are two solutions of \eqref{eq6-11} with the same initial value $u_0$. Let $\mathcal{W}=u_1^\pi-u_2^\pi$, we have
\begin{equation}\label{eq6-32}
\begin{split}
&\frac{d|\mathcal{W}(t)|^2}{dt}+2\eta\|\mathcal{W}(t)\|^2+2(F(u_1^\pi(t))-F(u^\pi_2(t)),\mathcal{W}(t))\\[1.0ex]
&\leq2(g(t,u^\pi_1(t))\sigma(t)-g(t,u^\pi_2(t))\sigma(t),\mathcal{W}(t))\\[1.0ex]
&~~+2\int_E(h(t,u_1^\pi(t),\xi)-h(t,u_2^\pi(t),\xi),\mathcal{W}(t))(\rho(t,\xi)-1)\lambda(d\xi).
\end{split}
\end{equation}
With the help of assumption of $(F_3)$, we arrive at
\begin{equation}\label{eq6-33}
2(F(u_1^\pi(t))-F(u_2^\pi(t)),\mathcal{W}(t))=2\int_{\mathcal{O}}(F(u_1^\pi(t,x))-F(u_2^\pi(t,x)))\mathcal{W}(t,x)dx\geq -2k_4|\mathcal{W}(t)|^2.
\end{equation}
By condition $(g_1)$, we derive
\begin{equation}\label{eq6-34}
\begin{split}
2(g(t,u^\pi_1(t))\sigma(t)-g(t,u^\pi_2(t))\sigma(t),\mathcal{W}(t))&\leq 2\|g(t,u^\pi_1(t))-g(t,u^\pi_2(t))\|_{\mathcal{L}_2(U;\hh)}\|\sigma(t)\|_U|\mathcal{W}(t)|\\[0.6ex]
&\leq 2\sqrt{L_g}|\mathcal{W}(t)|^2\|\sigma\|_U.
\end{split}
\end{equation}
For the last term, we have 
\begin{eqnarray}\label{eq6-35}
&&\hskip-1cm2\int_E(h(t,u_1^\pi(t),\xi)-h(t,u_2^\pi(t),\xi),\mathcal{W}(t))(\rho(t,\xi)-1)\lambda(d\xi)\\
&&\hskip-1cm=2\int_E\frac{|h(t,u_1^\pi(t),\xi)-h(t,u_2^\pi(t),\xi)|}{|u_1^\pi(t)-u_2^\pi(t)|}|\mathcal{W}(t)|^2|\rho(t,\xi)-1|\lambda(d\xi)\notag\\
&&\hskip-1cm\leq 2|\mathcal{W}(t)|^2\int_E\|h(t,\xi)\|_{1,\hh}|\rho(t,\xi)-1|\lambda(d\xi).\notag
\end{eqnarray}
Substituting  \eqref{eq6-33}-\eqref{eq6-35} into \eqref{eq6-32}, we obtain
\begin{equation*}
\frac{d|\mathcal{W}|^2}{dt}+2\eta\|\mathcal{W}\|^2\leq2\left(k_4+\sqrt{L_g}\|\sigma\|_U+\int_E\|h(t,\xi)\|_{1,\hh}|\rho(t,\xi)-1|\lambda(d\xi)\right)|\mathcal{W}|^2.
\end{equation*}
The Gronwall lemma and the fact $\sigma\in\hat{\mathrm{S}}^\Upsilon$, together with Lemma \ref{lem6-7}, conclude the proof of uniqueness of solution to \eqref{eq6-11}. 
\end{proof}

We now prove the main result. Recall that for $\pi=(\sigma,\rho)\in\mathbb{S}$, $\lambda_T^{\rho}(dtd\xi)=\rho(t,\xi)\lambda(d\xi)dt$. Define
\begin{equation}\label{eq6-36}
\mathcal{G}^0\left(\int_0^{\cdot}\sigma(s)ds, \lambda_T^{\rho}\right)=u^\pi,\qquad \mbox{for}~\pi=(\sigma,\rho)\in\mathbb{S}~\mbox{as given in Theorem \ref{thm6-12}}.
\end{equation}
Let $I:\mathbf{D}([0,T];\hh)\rightarrow [0,\infty]$ be defined as \eqref{eq6-6}.

\begin{theorem}({\bf Main Theorem})\label{thm6-13}
Suppose  $(F_1)$-$(F_3)$, $(g_1)$-$(g_3)$ and Conditions \ref{con6-4}-\ref{con6-5} hold. 
Then, the family of solutions $\{\ue\}_{\varepsilon>0}$ satisfies a large deviation principle on $\mathbf{D}([0,T];\hh)$ with the good rate function $I$ with respect to the topology of uniform convergence.
\end{theorem}

\begin{proposition}({\bf Verifying Condition \ref{con6-1}$(a)$})\label{pro6-14}
Fix $\Upsilon\in\mathbb{N}$. Let $\pi_n=(\sigma_n,\rho_n)$, $\pi=(\sigma,\rho)\in\mathrm{S}^\Upsilon$ be such that $\pi_n\rightarrow\pi$ as $n\rightarrow \infty$. Then, for $\mathcal{G}^0$ defined as in \eqref{eq6-36}, we have
\begin{equation*}
\mathcal{G}^0\left(\int_0^{\cdot}\sigma_n(s)ds,\lambda_T^{\rho_n}\right)\rightarrow\mathcal{G}^0\left(\int_0^{\cdot}\sigma(s)ds,\lambda_T^\rho\right),\qquad \mbox{in}~~C([0,T];\hh).
\end{equation*}
\end{proposition}
\begin{proof}  By definition \eqref{eq6-36}, we know that $\mathcal{G}^0\left(\int_0^{\cdot}\sigma_n(s)ds,\lambda_T^{\rho_n}\right)=u^{\pi_n}$. Since $\pi_n\in\mathrm{S}^\Upsilon\subset\mathbb{S}$, using the similar arguments as for \eqref{eq6-16}-\eqref{eq6-17}, we deduce that there exist two positive constants $C_\Upsilon$ and $C_{\alpha,\Upsilon}$, such that
\begin{equation}\label{eq6-37}
\sup_{t\in[0,T]}|u^{\pi_n}(t)|^2+\int_0^T\|u^{\pi_n}(t)\|^2dt+\int_0^T\|u^{\pi_n}(t)\|_{p+1}^{p+1}dt\leq C_\Upsilon,
\end{equation}
and for $\alpha\in(0,\frac{1}{2})$,
\begin{equation}\label{eq6-38}
\|u^{\pi_n}\|^2_{W^{\alpha,q}(0,T;V^*)}\leq C_{\alpha,\Upsilon}.
\end{equation}
Hence,  it follows from Lemma \ref{lemsoblev} that there exist   an element $\bar{u}\in L^2(0,T;\vv)\cap L^\infty(0,T;\hh)\cap L^{p+1}(0,T;L^{p+1}(\mathbb{R}^d))$ and a subsequence $u^{\pi_{n}}$ (relabeled the same) such that, as $n\rightarrow\infty$,
\begin{equation}\label{eq6-39}
\left\{
\begin{array}{rcl}
\begin{aligned}
 &u^{\pi_n}\rightarrow \bar{u} ~~\mbox{weak-star in}~~  L^{\infty}(0,T;\hh); \\
 &u^{\pi_n}\rightarrow \bar{u} ~~ \mbox{weakly in}~~ L^{2}(0,T;\vv);\\
 &u^{\pi_n}\rightarrow \bar{u} ~~ \mbox{weakly in}~~ L^{p+1}(0,T;L^{p+1}(\mathbb{R}^{d}));\\
 &u^{\pi_n}\rightarrow \bar{u} ~~ \mbox{strongly in}~~ L^{q}(0,T;\hh);\\
  &F(u^{\pi_n})\rightarrow F(\bar{u}) ~~ \mbox{weakly in}~~ L^{q}(0,T;L^q(\mathbb{R}^d)).
 \end{aligned}
\end{array}
\right.
\end{equation}

Next, we will prove $\bar{u}=u^\pi$. Let $\psi$ be a continuously differentiable function on $[0,T]$ with $\psi(T)=0$. We multiply $u^{\pi_n}(t)$ by $\psi(t)e_j$, then use  integration by parts to obtain
\begin{equation}\label{eq6-40}
\begin{split}
&-\int_0^T(u^{\pi_n}(t),e_j)\psi^{\prime}(t)dt+\int_0^Ta(u^{\pi_n}(t),e_j)\psi(t)dt+\int_0^T(F(u^{\pi_n}(t)),e_j)\psi(t)dt=(u_0,e_j)\psi(0)\\
&+\int_0^T(g(t,u^{\pi_n}(t))\sigma_n(t),e_j)\psi(t)dt+\int_0^T\left(\int_Eh(t,u^{\pi_n}(t),\xi)(\rho_n(t,\xi)-1)\lambda(d\xi),e_j\right)\psi(t)dt.
\end{split}
\end{equation}
Set
\begin{eqnarray*} I_n^1(T)&=&\int_0^T\int_E\left(h(t,u^{\pi_n}(t),\xi)(\rho_n(t,\xi)-1),e_j\right)\psi(t)\lambda(d\xi)dt,\\
I_n^2(T)&=&\int_0^T\int_E\left(h(t,u^{\pi_n}(t),\xi)(\rho(t,\xi)-1),e_j\right)\psi(t)\lambda(d\xi)dt,\\
 I(T)&=&\int_0^T\int_E\left(h(t,\bar{u}(t),\xi)(\rho(t,\xi)-1),e_j\right)\psi(t)\lambda(d\xi)dt.
 \end{eqnarray*}
Thus, we have
\begin{equation}\label{eq6-41}
I_n^1(T)-I(T)=I_n^1(T)-I_n^2(T)+I_n^2(T)-I(T).
\end{equation}
It follows from \eqref{eq6-30} that
\begin{equation}\label{eq6-42}
\lim_{n\rightarrow \infty}(I_n^2(T)-I(T))=0.
\end{equation}
To obtain the result, it is enough to prove that there exists a subsequence $\{m\}$ of $\{n\}$ such that
\begin{equation}\label{eq6-43}
\lim_{m\rightarrow \infty}(I_m^1(T)-I_m^2(T))=0.
\end{equation}
Thanks to \eqref{eq6-37}, Lemmas \ref{lem6-7} and \ref{lem6-9}, we  infer that for any given $\varepsilon>0$, there exists a compact subset $K_\varepsilon\subset E$ such that
\begin{eqnarray}\label{eq6-44}
&&\int_0^T\int_{K_\varepsilon^c}(h(t,u^{\pi_m}(t),\xi)(\rho_m(t,\xi)-1),e_j)\psi(t)\lambda(d\xi)dt\leq \int_0^T\int_{K_\varepsilon^c}|h(t,u^{\pi_m}(t),\xi)||\rho_m(t,\xi)-1||\psi(t)|\lambda(d\xi)dt\notag\\
&&\hskip2cm \leq\int_0^T\int_{K_\varepsilon^c}\|h(t,\xi)\|_{0,\hh}(1+|u^{\pi_m}(t)|)|\rho_m(t,\xi)-1||\psi(t)|\lambda(d\xi)dt\\
&&\leq \left(1+\sup_{t\in[0,T]}|u^{\pi_m}(t)|\right)\sup_{t\in[0,T]}|\psi(t)|
\int_0^T\int_{K_\varepsilon^c}\|h(t,\xi)\|_{0,\hh}|\rho_m(t,\xi)-1|\lambda(d\xi)dt\leq \sup_{t\in[0,T]}|\psi(t)|(1+C_\Upsilon)\varepsilon,\notag
\end{eqnarray}
and
\begin{equation}\label{eq6-45}
\int_0^T\int_{K_\varepsilon^c}(h(t,u^{\pi_m}(t),\xi)(\rho(t,\xi)-1),e_j)\psi(t)\lambda(d\xi)dt\leq \sup_{t\in[0,T]}|\psi(t)|(1+C_\Upsilon)\varepsilon.
\end{equation}
To prove \eqref{eq6-43}, applying a diagonal principle, it suffices to show that, for every compact $K\subset E$ and $\eta=2\sup_{t\in[0,T]}|\psi(t)|(1+C_\Upsilon)\varepsilon>0$, there exists a subsequence $\{m\}$ (denoted the same) such that
\begin{equation}\label{eq6-46}
\begin{split}
&\lim_{m\rightarrow\infty}\left|\int_0^T\int_K(h(t,u^{\pi_m}(t),\xi)(\rho_m(t,\xi)-1),e_j)\psi(t)\lambda(d\xi)dt\right.\\
&\qquad\qquad\left.-\int_0^T\int_K(h(t,u^{\pi_m}(t),\xi)(\rho(t,\xi)-1),e_j)\psi(t)\lambda(d\xi)dt\right|\\
&\qquad=\lim_{m\rightarrow\infty}\left|\int_0^T\int_K(h(t,u^{\pi_m}(t),\xi)\rho_m(t,\xi),e_j)\psi(t)\lambda(d\xi)dt\right.\\
&\qquad\qquad\left.-\int_0^T\int_K(h(t,u^{\pi_m}(t),\xi)\rho(t,\xi),e_j)\psi(t)\lambda(d\xi)dt\right|\leq \eta.
\end{split}
\end{equation}
Denote $A_M=\{(t,\xi)\in[0,T]\times K:\|h(t,\xi)\|_{0,\hh}\geq M\}$.  By Lemma \ref{lem6-9}$(iii)$ and \eqref{eq6-37}, for any $\varepsilon>0$, there exists $M>0$ such that
\begin{eqnarray}\label{eq6-47}
&&\int_0^T\int_K|(h(t,u^{\pi_m}(t),\xi)\rho_m(t,\xi),e_j)| |\psi(t)| 1_{A_M}\lambda(d\xi)dt\leq \int_0^T\int_K|h(t,u^{\pi_m}(t),\xi)|\rho_m(t,\xi)|\psi(t)|1_{A_M}\lambda(d\xi)dt\notag\\
&&\hskip2cm \leq \left(1+\sup_{t\in[0,T]}|u^{\pi_m}(t)|\right)\sup_{t\in[0,T]}|\psi(t)|
\int_0^T\int_K \rho_m(t,\xi)\|h(t,\xi)\|_{0,\hh}1_{A_M}\lambda(d\xi)dt\\
&&\hskip2cm\leq \sup_{t\in[0,T]}|\psi(t)|(1+C_\Upsilon)\varepsilon,\notag
\end{eqnarray}
and
\begin{equation}\label{eq6-48}
 \int_0^T\int_K|(h(t,u^{\pi_m}(t),\xi)\rho(t,\xi),e_j)| |\psi(t)| 1_{A_M}\lambda(d\xi)dt\leq \sup_{t\in[0,T]}|\psi(t)|(1+C_\Upsilon)\varepsilon.
\end{equation}
Denote 
$H_m(t,\xi)=(h(t,u^{\pi_m}(t),\xi),e_j)\psi(t)$ and $ H(t,\xi)=(h(t,\bar{u}(t),\xi),e_j)\psi(t)$.
Then, we have 
\begin{equation}\label{eq6-49}
\begin{split}
|H_m(t,\xi) 1_{A_M^c}|&\leq |h(t,u^{\pi_m}(t),\xi)||\psi(t)|1_{A_M^c}\leq \left(1+\sup_{t\in[0,T]}|u^{\pi_m}(t)|\right)\sup_{t\in[0,T]}|\psi(t)|\|h(t,\xi)\|_{0,\hh}1_{A_M^c}\\
&\leq \sup_{t\in[0,T]}|\psi(t)|(1+C_\Upsilon)M,
\end{split}
\end{equation}
and
$
|H(t,\xi) 1_{A_M^c}|\leq \sup_{t\in[0,T]}|\psi(t)|(1+C_\Upsilon)M.
$

Let $\Theta(\cdot)=\frac{\lambda_T(\cdot\cap [0,T]\times K)}{\lambda_T([0,T]\times K)}$ be a probability measure on $[0,T]\times E$. It follows from 
\cite[Proposition 4.1]{Z1} that there exists a subsequence, denoted the same, such that
$ \lim_{m\rightarrow \infty}H_m=H,\ \Theta\mbox{-a.s.}$
Therefore, using  similar arguments  to those in  \cite[Lemma 3.4]{B2}, together with \cite[Lemma 2.8]{B4} and \eqref{eq6-49}, we know there exists a subsequence $\{m^{\prime}\}$ of $\{m\}$, such that
\begin{equation}\label{eq6-50}
\lim_{m^\prime\rightarrow\infty}\int_0^T\int_KH_{m^\prime}(t,\xi)\rho_{m^\prime}(t,\xi)1_{A_M^c}\lambda(d\xi)dt
=\int_0^T\int_KH(t,\xi)\rho(t,\xi)1_{A_M^c}\lambda(d\xi)dt,
\end{equation}
and
\begin{equation*}
\lim_{m^\prime\rightarrow\infty}\int_0^T\int_KH_{m^\prime}(t,\xi)\rho(t,\xi)1_{A_M^c}\lambda(d\xi)dt
=\int_0^T\int_KH(t,\xi)\rho(t,\xi)1_{A_M^c}\lambda(d\xi)dt.
\end{equation*}
Hence, the above inequality, \eqref{eq6-47}-\eqref{eq6-48} and \eqref{eq6-50} imply \eqref{eq6-46}.
 Moreover, by \eqref{eq6-44}-\eqref{eq6-46}, we obtain
\begin{equation*}
\lim_{m^\prime\rightarrow\infty}|I^1_{m^\prime}(T)-I^2_{m^\prime}(T)|\leq 4  \sup_{t\in[0,T]}|\psi(t)|(1+C_\Upsilon)\varepsilon.
\end{equation*}
Thus, \eqref{eq6-43} follows immediately and there exists a subsequence of $\{m^\prime\}$ (still denoted the same) such that
\begin{equation}\label{eq6-51}
\lim_{m^\prime\rightarrow \infty} I^1_{m^\prime}(T)=I(T).
\end{equation}
Let us proceed likewise as before, we infer that
\begin{equation}\label{eq6-52}
\lim_{m^\prime\rightarrow\infty}\int_0^T(g(t,u^{\pi_{m^\prime}}(t))\sigma_{m^\prime}(t),e_j)\psi(t)dt=\int_0^T(g(t,\bar{u}(t))\sigma(t),e_j)\psi(t)dt.
\end{equation}
By \eqref{eq6-40} and \eqref{eq6-51}-\eqref{eq6-52}, using the same arguments as \eqref{eqlimit}, we see $\bar{u}$ satisfies 
\begin{equation}\label{eq6-53}
\begin{split}
&-\int_0^T(\bar{u}(t),e_j)\psi^{\prime}(t)dt+\int_0^Ta(\bar{u}(t),e_j)\psi(t)dt+\int_0^T(F(\bar{u}(t)),e_j)\psi(t)dt=(u_0,e_j)\psi(0)\\
&+\int_0^T(g(t,\bar{u}(t))\sigma(t),e_j)\psi(t)dt+\int_0^T\left(\int_Eh(t,\bar{u}(t),\xi)(\rho(t,\xi)-1)\lambda(d\xi),e_j\right)\psi(t)dt.
\end{split}
\end{equation}
Based on the uniqueness of solution to problem \eqref{eq6-11}, we conclude that $\bar{u}=u^\pi$.

At last, we will prove $u^{\pi_n}\rightarrow u^\pi$ in $C([0,T];\hh)$ as $n\rightarrow \infty$. Let $\mathcal{W}_n=u^{\pi_n}-u^\pi$, then
\begin{equation*}
\begin{split}
&\quad \frac{d|\mathcal{W}_n(t)|^2}{dt}+2\eta\|\mathcal{W}_n(t)\|^2+2(F(u^{\pi_n}(t))-F(u^\pi(t)),\mathcal{W}_n(t))\\[1.0ex]
&\leq 2\left(g(t,u^{\pi_n}(t))\sigma_n(t)-g(t,u^\pi(t))\sigma(t),\mathcal{W}_n(t)\right)\\[1.0ex]
&~~+2\int_E\left(h(t,u^{\pi_n}(t),\xi)(\rho_n(t,\xi)-1)-h(t,u^\pi(t),\xi)(\rho(t,\xi)-1),\mathcal{W}_n(t)\right)\lambda(d\xi)\\[1.0ex]
&\leq 2\left(g(t,u^{\pi_n}(t))\sigma_n(t)-g(t,u^{\pi_n}(t))\sigma(t),\mathcal{W}_n(t)\right)+2\left((g(t,u^{\pi_n}(t))\sigma(t)-g(t,u^\pi(t))\sigma(t),\mathcal{W}_n(t)\right)\\[1.0ex]
&~~+2\int_E\left((h(t,u^{\pi_n}(t),\xi)(\rho_n(t,\xi)-1)-h(t,u^{\pi_n}(t),\xi)(\rho(t,\xi)-1),\mathcal{W}_n(t)\right)\lambda(d\xi)\\
&~~+2\int_E\left(h(t,u^{\pi_n}(t),\xi)(\rho(t,\xi)-1)-h(t,u^\pi(t),\xi)(\rho(t,\xi)-1),\mathcal{W}_n(t)\right)\lambda(d\xi)\\[1.0ex]
&:=I_1^n(t)+I_2^n(t)+I_3^n(t)+I_4^n(t).
\end{split}
\end{equation*}
Similar to estimates \eqref{eq6-33}-\eqref{eq6-35}, we have
\begin{eqnarray*}
-2(F(u^{\pi_n}(t))-F(u^\pi(t)),\mathcal{W}_n(t))&\leq& 2k_4|\mathcal{W}_n(t)|^2,\\
I_2^n(t)&\leq& 2\sqrt{L_g}\|\sigma(t)\|_U|\mathcal{W}_n(t)|^2,\\
I_4^n(t)&\leq& 2|\mathcal{W}_n(t)|^2\int_E\|h(t,\xi)\|_{1,\hh}|\rho(t,\xi)-1|\lambda(d\xi).
\end{eqnarray*}
Subsequently, collecting all the estimates above, we  obtain
\begin{equation}\label{eq6-54}
\frac{d|\mathcal{W}_n(t)|^2}{dt}+2\eta\|\mathcal{W}_n(t)\|^2\leq \aleph(t)|\mathcal{W}_n(t)|^2+I_1^n(t)+I_3^n(t),
\end{equation}
where we have used the notation 
$$ \aleph(t)=2k_4+2\sqrt{L_g}\|\sigma(t)\|_U+2\int_E\|h(t,\xi)\|_{1,\hh}|\rho(t,\xi)-1|\lambda(d\xi)\in L^1(0,T).$$
Multiplying \eqref{eq6-54} by $e^{-\int_0^t\aleph(s)ds}$ and integrating it from $0$ to $t$, we obtain 
\begin{equation*}
e^{-\int_0^t\aleph(s)ds}|\mathcal{W}_n(t)|^2\leq \int_0^te^{-\int_0^s\aleph(r)dr}(I_1^n(s)+I_3^n(s))ds\leq \int_0^t(|I_1^n(s)|+|I_3^n(s)|)ds,
\end{equation*}
which implies
\begin{equation}\label{eq6-55}
\sup_{t\in[0,T]}|\mathcal{W}_n(t)|^2\leq \exp\left(\int_0^T\aleph(t)dt\right)\int_0^T(|I_1^n(t)|+|I_3^n(t)|)dt.
\end{equation}
Since $u \in C([0,T];\hh)$ and $u^{\pi_n}\in C([0,T];\hh)$ for each $n\in\mathbb{N}$ (see Theorem \ref{thm6-12}), together with the facts that $\sup_{n\in\mathbb{N}}\sup_{t\in[0,T]}|u^{\pi_n}(t)|\leq C_\Upsilon$,  $\sup_{t\in[0,T]}|u^{\pi}(t)|\leq C_\Upsilon$ and \eqref{eq6-39}, we know $u^{\pi_n}\rightarrow u^\pi $ in $L^2(0,T;\hh)$ as $n\rightarrow \infty$.
 By condition $(g_2)$ and the H\"older inequality, we find 
 \begin{equation*}
\begin{split}
\int_0^T|I_1^n(t)|dt&\leq 2\int_0^T\|g(t,u^{\pi_n}(t))\|_{\mathcal{L}_2(U;\hh)}\|\sigma_n(t)-\sigma(t)\|_U|\mathcal{W}_n(t)|dt\\
&\leq 2\sqrt{C_g}\int_0^T\sqrt{1+|u^{\pi_n}(t)|^2}(\|\sigma_n(t)-\sigma(t)\|_U)|\mathcal{W}_n(t)|dt\\
&\leq 2\sqrt{C_g(1+C_\Upsilon^2)}\left(\int_0^T\left(\|\sigma_n(t)\|_U+\|\sigma(t)\|_U\right)^2dt\right)^{\frac{1}{2}}\left(\int_0^T|\mathcal{W}_n(t)|^2dt\right)^{\frac{1}{2}}\longrightarrow 0,\ \mbox{as}\ n\rightarrow \infty.
\end{split}
\end{equation*}
Similarly, for $I_3^n(t)$, it follows from Condition \ref{con6-4}$(ii)$, Lemma \ref{lem6-7}$(i)$  and the Lebesgue dominated theorem that,
\begin{equation*}
\begin{split}
{\int_0^T}|I_3^n(t)|dt&\leq 
2\int_0^T\int_E(h(t,u^{\pi_n}(t),\xi)(\rho_n(t,\xi)-1),\mathcal{W}_n(t))\lambda(d\xi)dt\\
&~~+2\int_0^T\int_E(h(t,u^{\pi_n}(t),\xi)(\rho(t,\xi)-1),\mathcal{W}_n(t))\lambda(d\xi)dt\\
&\leq 2\int_0^T\int_E\|h(t,\xi)\|_{0,\hh}(1+|u^{\pi_n}(t)|)|\rho_n(t,\xi)-1||\mathcal{W}_n(t)|\lambda(d\xi)dt\\
&~~+2\int_0^T\int_E\|h(t,\xi)\|_{0,\hh}(1+|u^{\pi}(t)|)|\rho(t,\xi)-1||\mathcal{W}_n(t)|\lambda(d\xi)dt\\
&\leq 2(1+C_\Upsilon)\int_0^T\int_E\|h(t,\xi)\|_{0,\hh}|\rho_n(t,\xi)-1||\mathcal{W}_n(t)|\lambda(d\xi)dt\\
&~~+2(1+C_\Upsilon)\int_0^T\int_E\|h(t,\xi)\|_{0,\hh}|\rho(t,\xi)-1||\mathcal{W}_n(t)|\lambda(d\xi)dt\longrightarrow 0,\ \mbox{as}\ n\rightarrow \infty.
\end{split}
\end{equation*}
Hence, by \eqref{eq6-55}, we obtain
$ \lim_{n\rightarrow\infty}\sup_{t\in[0,T]}|u^{\pi_n}(t)-u^{\pi}(t)|=0,$
which implies $u^{\pi_n}\rightarrow u^\pi$ in $C([0,T];\hh)$. The proof of this proposition is complete. 
\end{proof} 

\begin{theorem}
Assume $(F_1)$-$(F_3)$, $(g_1)$-$(g_3)$ and Condition \ref{con6-4} hold. If $u_0\in\hh$, there exists a unique $\hh$-valued progressively measurable process, such that $\ue\in L^2(0,T;\vv)\cap \mathbf{D}([0,T];\hh)\cap L^{p+1}(0,T;L^{p+1}(\mathbb{R}^d))$ for any $T>0$, and 
\begin{equation}\label{eq6-56}
\begin{split}
\ue(t)&=u_0-\int_0^t(-\Delta)^\gamma\ue(s)ds-\delta\int_0^t\ue(s)ds
-\int_0^tF(\ue(s))ds\\
&~~+\sqrt{\varepsilon}\int_0^tg(s,\ue(s))dW(s)+\varepsilon\int_0^t\int_Eh(s,\ue(s-),\xi)\tilde{N}^{\varepsilon^{-1}}(dsd\xi),\quad \mbox{a.s.}
\end{split}
\end{equation}
\end{theorem}
This theorem can be proved similarly  as \cite[Theorem 1.2]{B5} since $(-\Delta)^{\gamma}$ is a linear operator, showing \eqref{eq6-56} admits a strong solution (in the probability sense).
In particular, for every $\varepsilon >0$, there exists a measurable map $\mathcal{G}^\varepsilon:\bar{\mathbb{U}}\rightarrow D([0,T];\hh)$ such that, for any Poisson random measure $n^{\varepsilon^{-1}}$ on $[0,T]\times E$ with intensity measure $\varepsilon^{-1}L_T\otimes \lambda$ given in some probability space, $\mathcal{G}^{\varepsilon}(\sqrt{\varepsilon}W,\varepsilon n^{\varepsilon^{-1}})$ is the unique solution of \eqref{eq6-56} with $\tilde{N}^{\varepsilon^{-1}}$ replaced by $\tilde{n}^{\varepsilon^{-1}}$.

We have the following lemma introduced in \cite[Lemma 2.3]{B3}.
\begin{lemma}
Let $\phi_\varepsilon=(\psi_\varepsilon,\varphi_\varepsilon)\in\tilde{\mathcal{U}}^\Upsilon$ and $\lambda_\varepsilon=\frac{1}{\varphi_\varepsilon}$.  Then,
\begin{equation*}
\begin{split}
\tilde{\mathcal{E}}_t^{\varepsilon}(\lambda_\varepsilon):=&\exp\left\{\int_{[0,t]\times E\times [0,\varepsilon^{-1}]}\log(\lambda_\varepsilon(s,x))\bar{N}(dsdxdr)+\int_{[0,t]\times E\times[0,\varepsilon^{-1}]}(-\lambda_\varepsilon(s,x)+1)\bar{\lambda}_T(dsdxdr)\right\},
\end{split}
\end{equation*}
and 
\begin{equation*}
\hat{\mathcal{E}}^\varepsilon_t(\psi_\varepsilon):=\exp\left\{\frac{1}{\sqrt{\varepsilon}}\int_0^t\psi_\varepsilon(s)dW(s)-\frac{1}{2\varepsilon}\int_0^t\|\psi_\varepsilon(s)\|_U^2ds\right\},
\end{equation*}
are $\{\bar{\mathcal{F}}_t^{\bar{\mathbb{U}}}\}$-martingales. Set
$ \mathcal{E}_t^\varepsilon(\psi_\varepsilon,\lambda_\varepsilon):=\hat{\mathcal{E}}_t^\varepsilon(\psi_\varepsilon)\tilde{\mathcal{E}}_t^\varepsilon(\lambda_\varepsilon).$
Then
$
\mathbb{Q}_t^\varepsilon(G)=\int_G\mathcal{E}^\varepsilon_t(\psi_\varepsilon,\lambda_\varepsilon)d\bar{\mathbb{P}}^{\bar{\mathbb{U}}}, G\in\mathcal{B}(\bar{\mathbb{U}}),
$
defines a probability measure on $\bar{\mathbb{U}}$.
\end{lemma}

Since $(\sqrt{\varepsilon}W+\int_0^{\cdot}\psi_\varepsilon(s)ds,\varepsilon N^{\varepsilon^{-1}\varphi_{\varepsilon}})$ under $\mathbb{Q}_T^\varepsilon$ has the same law as that of $(\sqrt{\varepsilon}W,\varepsilon N^{\varepsilon^{-1}})$ under $\bar{\mathbb{P}}^{\bar{\mathbb{U}}}$, there exists a unique solution $\tilde{u}^\varepsilon$ to the following controlled stochastic fractional differential equation,
\begin{eqnarray}\label{eq6-57}
&&\tu(t)=u_0-\int_0^t(-\Delta)^{\gamma}\tu(s)ds-\delta\int_0^t\tu(s)ds-\int_0^tF(\tu(s))ds\notag\\
&&~~+\int_0^tg(s,\tu(s))\psi_\varepsilon(s)ds+\sqrt{\varepsilon}\int_0^tg(s,\tu(s))dW(s)\notag\\
&&~~+\varepsilon\int_0^t\int_Eh(s,\tu(s-),\xi)\left(N^{\varepsilon^{-1}\varphi_{\varepsilon}}(dsd\xi)-\varepsilon^{-1}\lambda(d\xi)ds\right)\notag\\
&&\qquad=u_0-\int_0^t(-\Delta)^{\gamma}\tu(s)ds-\delta\int_0^t\tu(s)ds-\int_0^tF(\tu(s))ds\\
&&~~+\int_0^tg(s,\tu(s))\psi_\varepsilon(s)ds+\sqrt{\varepsilon}\int_0^tg(s,\tu(s))dW(s)+\int_0^t\int_Eh(s,\tu(s-),\xi)(\varphi_\varepsilon(s,\xi)-1)\lambda(d\xi)ds\notag\\
&&~~+\varepsilon\int_0^t\int_Eh(s,\tu(s-),\xi)\left(N^{\varepsilon^{-1}\varphi_{\varepsilon}}(dsd\xi)-\varepsilon^{-1}\varphi_\varepsilon(s,\xi)\lambda(d\xi)ds\right).\notag
\end{eqnarray}
Moreover, we have 
\begin{equation}\label{eq6-58}
\mathcal{G}^\varepsilon\left(\sqrt{\varepsilon}W+\int_0^\cdot\psi_\varepsilon(s)ds,\varepsilon N^{\varepsilon^{-1}\varphi_{\varepsilon}}\right)=\tu.
\end{equation}

The following estimates will be used later.
\begin{lemma}\label{lem6-17}
Assume $(F_1)$-$(F_3)$, $(g_1)$-$(g_3)$ and Condition \ref{con6-4} hold. 
Let {$u_0\in\hh$}. Then
there exists $0<\varepsilon_0\leq \frac{1}{(1+16C_b^2C_g T+32C_b^2C_{0,2}^\Upsilon)^2}$, such that
\begin{equation}\label{eq6-59}
\sup_{0<\varepsilon<\varepsilon_0}\left[\be\sup_{t\in[0,T]}|\tu(t)|^2+\be\int_0^T\|\tu(t)\|^2dt+\be\int_0^T\|\tu(t)\|_{p+1}^{p+1}dt\right]<\infty,
\end{equation}
where $C_b$ is the constant obtained from the Burkholder-Davis-Gundy inequality and $\be$ is the expectation operator corresponding to {$\bar{\mathbb{P}}:=\bar{\mathbb{P}}^{\bar{\mathbb{U}}}$}.
Moreover, for $\alpha\in(0,1/2)$, there exists $C_\alpha>0$, such that
\begin{equation}\label{eq6-60}
\sup_{0<\varepsilon<\varepsilon_0}\be\|\tu\|^2_{W^{\alpha,q}(0,T;V^*)}\leq C_\alpha.
\end{equation}
Thus, the family $\{\tu,0<\varepsilon<\varepsilon_0\}$ is tight in {$L^q(0,T;\hh)$}.
\end{lemma}
\begin{proof}  The details to prove \eqref{eq6-59} is shown in Appendix 1.
Notice that, \eqref{eq6-57} is equivalent to 
\begin{equation*}
\begin{split}
\tu(t)&=u_0-\left(\int_0^t(-\Delta)^{\gamma}\tu(s)ds+\delta\int_0^t\tu(s)ds\right)-\int_0^tF(\tu(s))ds\\
&~~+\int_0^tg(s,\tu(s))\psi_\varepsilon(s)ds+\sqrt{\varepsilon}\int_0^tg(s,\tu(s))dW(s)\\
&~~+\int_0^t\int_Eh(s,\tu(s-),\xi)(\varphi_\varepsilon(s,\xi)-1)\lambda(d\xi)ds+\varepsilon\int_0^t\int_Eh(s,\tu(s-),\xi)\tilde{N}^{\varepsilon^{-1}\varphi_{\varepsilon}}(dsd\xi)\\
&:=J^1+J_\varepsilon^2(t)+J_\varepsilon^3(t)+J_\varepsilon^4(t)+J_\varepsilon^5(t)+J_\varepsilon^6(t)+J_\varepsilon^7(t).
\end{split}
\end{equation*}
By the same arguments as in the proof of \cite[Theorem 3.1]{F1}, we know there exists $\mathrm{C}^1>0$, such that
\begin{equation}\label{eq6-61}
\sup_{0<\varepsilon<\varepsilon_0}\be|J^1|^2\leq \mathrm{C}^1.
\end{equation}
For $J_\varepsilon^2$, using the same method as in \cite[Theorem 2.3]{W2} and the H\"older inequality,  we infer there exists a constant $\mathrm{C}^2>0$ such that
\begin{equation}\label{eq6-62}
\begin{split}
\sup_{0<\varepsilon<\varepsilon_0}\be\|J^2_\varepsilon\|^2_{W^{1,2}(0,T;\vv^*)}&=\sup_{0<\varepsilon<\varepsilon_0}\left(\be\int_0^T\left\|\int_0^tA\tu(s)ds\right\|_*^2dt+\be\int_0^T\|A\tu(t)\|^2_*dt\right)\\
&\leq \sup_{0<\varepsilon<\varepsilon_0}(T^2+1)\int_0^T\be\|A\tu(t)\|_*^2dt\leq\mathrm{C}^2.
\end{split}
\end{equation} 
For $J_\varepsilon^3$, similar to \eqref{eq6-24}, by condition $(F_1)$ and the H\"older inequality, we know there exists a constant $\mathrm{C}^3>0$, such that
\begin{eqnarray}\label{eq6-63}
&&\sup_{0<\varepsilon<\varepsilon_0}\be\|J_\varepsilon^3\|^q_{W^{1,q}(0,T;L^q(\mathcal{O}))}=\sup_{0<\varepsilon<\varepsilon_0}\left(\be\int_0^T\left\|\int_0^tF(\tu(s))ds\right\|^q_{L^q(\mathcal{O})}dt+\be\int_0^T\|F(\tu(t))\|^2_{L^q(\mathcal{O})}dt\right)\notag\\
&&\quad\hskip4cm\leq \sup_{0<\varepsilon<\varepsilon_0}\left(T^{\frac{1}{p}+1}+1\right)\be\int_0^T\|F(\tu(t))\|^q_qdt\leq \mathrm{C}^3.
\end{eqnarray}
To estimate $J_\varepsilon^4$, we  apply condition $(g_2)$ and the H\"older inequality for $0\leq s<t\leq T$, 
\begin{equation*}
\begin{split}
\be|J_\varepsilon^4(t)-J_\varepsilon^4(s)|^2&=\be\left|\int_s^tg(r,\tu(r))\psi_\varepsilon(r)dr\right|^2\leq \be\left(\int_s^t|g(r,\tu(r))\psi_\varepsilon(r)|dr\right)^2\\
&\leq \be\left(\int_s^t\|g(r,\tu(r))\|_{\mathcal{L}_2(U;\hh)}\|\psi_\varepsilon(r)\|_Udr\right)^2\leq C_g\be\left(\int_s^t\sqrt{1+|\tu(r)|^2}\|\psi_\varepsilon(r)\|_Udr\right)^2\\
&\leq 2C_g\Upsilon(t-s)\be\left(1+\sup_{t\in[0,T]}|\tu(t)|^2\right),
\end{split}
\end{equation*}
the last inequality holds since $(\psi_{\varepsilon},\varphi_{\varepsilon})\in\tilde{\mathcal{U}}^{\Upsilon}$.
Consequently, by the above estimate and the H\"older inequality, for $\alpha\in(0,1/2)$, we have
\begin{equation*}
\begin{split}
&\be\|J_\varepsilon^4\|^2_{W^{\alpha,2}(0,T;\hh)}=\be\int_0^T|J^4_\varepsilon(t)|^2dt+\be\int_0^T\int_0^T\frac{|J_\varepsilon^4(t)-J_\varepsilon^4(s)|^2}{|t-s|^{1+2\alpha}}dsdt\\
&\leq T^2\be\int_0^T\|g(r,\tu(r))\|^2_{\mathcal{L}_2(U;\hh)}\|\psi_\varepsilon(r)\|_U^2dr\\
&~~+2C_g\Upsilon \be\left(1+\sup_{t\in[0,T]}|\tu(t)|^2\right)\int_0^T\int_0^T\frac{1}{|t-s|^{2\alpha}}dsdt\\
&\leq 2 C_g\Upsilon \be\left(1+\sup_{t\in[0,T]}|\tu(t)|^2\right)\left(T^2+\frac{1}{(1-2\alpha)(2-2\alpha)}T^{2-2\alpha}\right).
\end{split}
\end{equation*}
Therefore, there exists a constant $\mathrm{C}^4>0$ such that
\begin{equation}\label{eq6-64}
\sup_{0<\varepsilon<\varepsilon_0}\be\|J^4_\varepsilon\|^2_{W^{\alpha,2}(0,T;\hh)}\leq \mathrm{C}^4.
\end{equation}
For $J_\varepsilon^5$, similar to $J_\varepsilon^4$,  by It\^o's isometry and  condition $(g_2)$, for $0\leq s<t\leq T$, we find
\begin{equation*}
\begin{split}
\be |J_\varepsilon^5(t)-J_\varepsilon^5(s)|&=
\be\left|\sqrt{\varepsilon}\int_s^tg(r,\tu(r))dW(r)\right|^2\leq \varepsilon\be\int_s^t\|g(r,\tu(r))\|^2_{\mathcal{L}_2(U;\hh)}dr\\
&\leq \varepsilon C_g(t-s)\be\left(1+\sup_{t\in[0,T]}|\tu(t)|^2\right).
\end{split}
\end{equation*}
Thus,  for $\alpha\in(0,1/2)$, there exists a constant $\mathrm{C}^5>0$ such that
\begin{equation}\label{eq6-65}
\begin{split}
&\sup_{0<\varepsilon<\varepsilon_0}\be\|J^5_\varepsilon\|^2_{W^{\alpha,2}(0,T;\hh)}\leq\sup_{0<\varepsilon<\varepsilon_0}\left(\be\int_0^T|J_\varepsilon^5(t)|^2dt+\be\int_0^T\int_0^T\frac{|J_\varepsilon^5(t)-J^5_\varepsilon(s)|^2}{|t-s|^{1+2\alpha}}dsdt\right)\\
&\leq \sup_{0<\varepsilon<\varepsilon_0}\left[\varepsilon C_g\be\left(1+\sup_{t\in[0,T]}|\tu(t)|^2\right)\left(T^2+\frac{T^{2-2\alpha}}{(1-2\alpha)(2-2\alpha)}\right)\right]\leq \mathrm{C}^5.
\end{split}
\end{equation}
For $J_\varepsilon^6$ and $0\leq s<t\leq T$, we have
\begin{equation*}
\begin{split}
\be|J_\varepsilon^6(t)-J_\varepsilon^6(s)|^2&=\be\left|\int_s^t\int_Eh(r,\tu(r-),\xi)(\varphi_\varepsilon(r,\xi)-1)\lambda(d\xi)dr\right|^2\\
&\leq \be\left(\int_s^t\int_E\|h(r,\xi)\|_{0,\hh}(1+|\tu(r)|)|\varphi_\varepsilon(r,\xi)-1|\lambda(d\xi)dr\right)^2\\
&\leq \be\left[\left(1+\sup_{t\in[0,T]}|\tu(t)|\right)^2\left(\int_s^t\int_E\|h(r,\xi)\|_{0,\hh}|\varphi_\varepsilon(r,\xi)-1|\lambda(d\xi)dr\right)^2\right]\\
&\leq C_{0,1}^\Upsilon\be\left[\left(1+\sup_{t\in[0,T]}|\tu(t)|\right)^2\int_s^t\int_E\|h(r,\xi)\|_{0,\hh}|\varphi_\varepsilon(r,\xi)-1|\lambda(d\xi)dr\right],
\end{split}
\end{equation*}
where $C_{0,1}^\Upsilon$ appears in Lemma \ref{lem6-7} 
(see \eqref{eq6-8}). Using the above estimate to $J_\varepsilon^6$, we obtain
\begin{equation*}
\begin{split}
&\be\|J_\varepsilon^6\|^2_{W^{\alpha,2}(0,T;\hh)}
=\be\int_0^T|J_\varepsilon^6(t)|^2dt+\e\int_0^T\int_0^T\frac{|J_\varepsilon^6(t)-J_\varepsilon^6(s)|^2}{|t-s|^{1+2\alpha}}dsdt\\
&\leq\be\int_0^T\left|\int_0^t\int_Eh(s,\tu(s-),\xi)(\varphi_\varepsilon(s,\xi)-1)\lambda(d\xi)ds\right|^2dt\\
&~~+C_{0,1}^\Upsilon\be\left[\left(1+\sup_{t\in[0,T]}|\tu(t)|\right)^2\int_0^T\int_0^T\int_s^t\int_E\frac{\|h(r,\xi)\|_{0,\hh}|\varphi_\varepsilon(r,\xi)-1|}{|t-s|^{2\alpha+1}}\lambda(d\xi)drdtds\right]\\
&\leq (C_{0,1}^\Upsilon)^2\be\left(1+\sup_{t\in[0,T]}|\tu(t)|\right)^2\left( T+\frac{T^{1-2\alpha}}{2\alpha(1-2\alpha)}\right).
\end{split}
\end{equation*}
Therefore, there exists a constant $\mathrm{C}^6>0$ such that
\begin{equation}\label{eq6-66}
\sup_{0<\varepsilon<\varepsilon_0}\be\|J_\varepsilon^6\|^2_{W^{\alpha,2}(0,T;\hh)}\leq\mathrm{C}^6.
\end{equation}
For the last term $J_\varepsilon^7$ and $0\leq s<t\leq T$, by Lemma \ref{lem6-7}, we derive
\begin{equation*}
\begin{split}
\be|J_\varepsilon^7(t)-J_\varepsilon^7(s)|^2&=\varepsilon^2\be\left|\int_s^t\int_Eh(r,\tu(r-),\xi)\tilde{N}^{\varepsilon^{-1}\varphi_\varepsilon}(drd\xi)\right|^2\leq \varepsilon\be\int_s^t\int_E|h(r,\tu(r-),\xi)|^2\varphi_\varepsilon(r,\xi)\lambda(d\xi)dr\\
&\leq \varepsilon\be\int_s^t\int_E\|h(r,\xi)\|^2_{0,\hh}|\varphi_\varepsilon(r,\xi)|(1+|\tu(r)|)^2\lambda(d\xi)dr\\
&\leq\varepsilon\be\left[\left(1+\sup_{t\in[0,T]}|\tu(t)|\right)^2\int_s^t\int_E\|h(r,\xi)\|^2_{0,\hh}|\varphi_\varepsilon(r,\xi)|\lambda(d\xi)dr\right].
\end{split}
\end{equation*}
Using similar arguments as for the bound of $J_\varepsilon^6$ and Lemma \ref{lem6-7}, we infer there exists $\mathrm{C}^7>0$ such that,
\begin{equation*}
\sup_{0<\varepsilon<\varepsilon_0}\be\|J_\varepsilon^7\|^2_{W^{\alpha,2}(0,T;\hh)}\leq\mathrm{C}^7,
\end{equation*}
which, combining with \eqref{eq6-61}-\eqref{eq6-66}, proves \eqref{eq6-60}. 
\end{proof}

To obtain the main results, we need to prove that $\{\tu,0<\varepsilon<\varepsilon_0\}$ is tight in the vector valued Skorokhod space $\mathbf{D}([0,T];D(A^{-r}))$ for some  $r>\frac{d}{4\gamma}(1-\frac{2}{p+1})$, such that $ D(A^r)\subset\vv\subset L^{p+1}(\mathcal{O})\subset \hh:=\hh\subset L^q(\mathcal{O})\subset\vv^*\subset D(A^{-r})$ (see \cite[Lemma 2.1]{W2} for more details). To that end, we  first recall the following two lemmas (see \cite{A1, J2} and the references therein).
\begin{lemma}\label{lem6-18}
Let $\verb"H"$ be a separable Hilbert space with  inner product $(\cdot,\cdot)$. For an orthonormal basis $\{\chi_k\}_{k\in\mathbb{N}}$ in $\verb"H"$,
define the function $r_L^2:\verb"H"\rightarrow\mathbb{R}^+$ by
$ r_L^2(x)=\sum_{k\geq L+1}(x,\chi_k)^2,\ L\in\mathbb{N}.$
Let $\textsf{B}$ be a total and closed under addition subset of $\verb"H"$. Then, a sequence $\{u_\varepsilon\}_{\varepsilon\in(0,1)}$ of stochastic processes with trajectories in $\mathbf{D}([0,T];\verb"H")$ is tight if and only if the following two conditions hold:
\begin{enumerate}
\item [$(i)$] $\{u_{\varepsilon}\}_{\varepsilon\in[0,1]}$ is $\textsf{B}$-weakly tight, that is, for every $l\in\textsf{B}$, $\{(u_\varepsilon,l)\}_{\varepsilon\in(0,1)}$ is tight in $\mathbf{D}([0,T];\mathbb{R})$;
\item [$(ii)$]  For every $\upsilon>0$,
\begin{equation}\label{eq6-69}
\lim_{L\rightarrow\infty}{\limsup_{\varepsilon\rightarrow 0}}\mathbb{P}\left(r_L^2(u_\varepsilon(s))>\upsilon~\mbox{for some}~s\in[0,T]\right)=0.
\end{equation}
\end{enumerate}
\end{lemma}

Let $\tu$ be defined as in \eqref{eq6-57}, then we have the following result.
\begin{lemma}\label{lem6-19}
The set $\{\tu,0<\varepsilon<\varepsilon_0\}$ is tight in $\mathbf{D}([0,T];D(A^{-r}))$.
\end{lemma}
\begin{proof}  Notice that, $\{\lambda_j^re_j\}_{j\in\mathbb{N}}$ is a complete orthonormal system of $D(A^{-r})$ (see, for example, \cite{W2}). Since
\begin{equation*}
\begin{split}
&\lim_{L\rightarrow\infty}\limsup_{\varepsilon\rightarrow 0}\be\sup_{t\in[0,T]}r_L^2(\tu(t))=\lim_{L\rightarrow\infty}\limsup_{\varepsilon\rightarrow 0}\be\sup_{t\in[0,T]}\sum_{j=L+1}^{\infty}\langle\tu(t),\lambda_j^re_j\rangle_{D(A^{-r})}^2\\
&= \lim_{L\rightarrow\infty}\limsup_{\varepsilon\rightarrow 0} \be\sup_{t\in[0,T]}\sum_{j=L+1}^{\infty}(A^{-r}\tu(t),e_j)_{\hh}^2=\lim_{L\rightarrow\infty}\limsup_{\varepsilon\rightarrow 0}\be\sup_{t\in[0,T]}\sum_{j=L+1}^{\infty}\frac{(\tu(t),e_j)_{\hh}^2}{\lambda_j^{2r}}\\
&\leq \lim_{L\rightarrow\infty}\frac{{\limsup_{\varepsilon\rightarrow 0}} \be[\sup_{t\in[0,T]}|\tu(t)|^2]}{\lambda_{L+1}^{2r}}=0.
\end{split}
\end{equation*}
Therefore, \eqref{eq6-69} holds with $\verb"H"=D(A^{-r})$ by using the Markov inequality.  

Choose $\textsf{B}=D(A^r)$. We claim that $\{\tu\}_{\varepsilon\in[0,1]}$ is $D(A^r)$-weakly tight by using the same method as \cite[Lemma 4.4]{Z1}. That is, for every $l\in D(A^r)$, let $(\beta_\varepsilon,d_\varepsilon)$ be a stopping time with respect to the natural $\bar{\sigma}$-field taking finitely many values, and an interval on $[0,T]$, respectively, satisfying  $d_\varepsilon\rightarrow 0$ as $\varepsilon\rightarrow 0$. By Lemma \ref{lem6-17}, it is easy to check 
$\{\langle \tu,l\rangle_{D(A^r)},0<\varepsilon<\varepsilon_0\}$ is tight
on the real line for all $t\in[0,T]$. Hence, we end up this proof by showing  $\langle\tu(\beta_\varepsilon+d_\varepsilon)-\tu(\beta_\varepsilon),l\rangle_{D(A^r)}\rightarrow 0$ in probability as $\varepsilon \rightarrow 0$ ( see Appendix 2 for the details). 
\end{proof}

We proceed likewise as in \cite[Proposition 3.1]{R2}, there exists a unique solution $\ty(t)$ ($t\geq0$) to the following equation with initial value $0$,
\begin{equation}\label{eq6-70}
\begin{split}
d\ty(t)&=-\left((-\Delta)^{\gamma}\ty(t)+\delta\ty(t)\right)dt+\sqrt{\varepsilon}g(t,\tu(t))dW(t)+\varepsilon\int_Eh(t,\tu(t-),\xi)\tilde{N}^{\varepsilon^{-1}\varphi_{\varepsilon}}(dtd\xi),
\end{split}
\end{equation}
and $\ty\in\mathbf{D}([0,T];\hh)\cap L^2(0,T;\vv)$, $\bar{\mathbb{P}}$-a.s.

\begin{lemma}\label{lem6-20}
There exist  some constants $\tilde{C}>0$ and $\tilde{\varepsilon}_0=\frac{1}{4(C_gT+2C_{0,2}^\Upsilon+4C_b^2C_gT+8C_b^2C_{0,2}^{\Upsilon})}$, such that for any $0<\varepsilon\leq\tilde{\varepsilon}_0$, the solution of equation \eqref{eq6-70} with initial value $\ty(0)=0$ satisfies
$$\be\left[\sup_{t\in[0,T]}|\ty(t)|^2\right]+\be\int_0^T\|\ty(t)\|^2dt\leq\tilde{C}\varepsilon.$$
\end{lemma}
\begin{proof}  By It\^o's formula, similar to \eqref{eq3-5}, we derive
\begin{equation}\label{eq6-71}
\begin{split}
|\ty(t)|^2+2\eta\int_0^t\|\ty(s)\|^2ds&\leq 2\int_0^t(\ty(s),\sqrt{\varepsilon}g(s,\tu(s))dW(s))+\varepsilon\int_0^t\|g(s,\tu(s))\|^2_{\mathcal{L}_2(U;\hh)}ds\\
&~~+2\varepsilon\int_0^t\int_E(\ty(s-), h(s,\tu(s-),\xi))\tilde{N}^{\varepsilon^{-1}\varphi_{\varepsilon}}(dsd\xi)\\
&~~+\varepsilon^2\int_0^t\int_E|h(s,\tu(s-),\xi)|^2N^{\varepsilon^{-1}\varphi_{\varepsilon}}(dsd\xi)\\
&:=I_y^1+I_y^2+I_y^3+I_y^4.
\end{split}
\end{equation}
By assumption $(g_2)$ and Lemma \ref{lem6-7}, we obtain
\begin{equation}\label{eq6-72}
\be\left[\sup_{t\in[0,T]}I_y^2\right]\leq \varepsilon\be \int_0^T\|g(t,\tu(t))\|^2_{\mathcal{L}_2(U;\hh)}dt\leq\varepsilon C_gT\left[1+\be\left(\sup_{t\in[0,T]}|\tu(t)|^2\right)\right],
\end{equation}
and 
\begin{equation}\label{eq6-73}
\begin{split}
\be\left[\sup_{t\in[0,T]}I_y^4\right]&\leq\varepsilon\be\int_0^T\int_E\|h(t,\xi)\|_{0,\hh}^2(1+\tu(t))^2\varphi_\varepsilon(t,\xi)\lambda(d\xi)dt\\
&\leq2\varepsilon C_{0,2}^\Upsilon\left[1+\be\left(\sup_{t\in[0,T]}|\tu(t)|^2\right)\right],
\end{split}
\end{equation}
respectively. As for $I_y^1$ and $I_y^3$, by the similar estimates as for \eqref{eq3-8}-\eqref{eq3-9}, we have 
\begin{equation}\label{eq6-74}
\be\left[\sup_{t\in[0,T]}I_y^1\right]\leq \frac{1}{4}\be\left[\sup_{t\in[0,T]}|\ty(t)|^2\right]+4C_b^2C_gT\varepsilon\left[1+\be\left(\sup_{t\in[0,T]}|\tu(t)|^2\right)\right],
\end{equation}
and
\begin{equation}\label{eq6-75}
\begin{split}
\be\left[\sup_{t\in[0,T]}I_y^3\right]&\leq\frac{1}{4}\be\left[\sup_{t\in[0,T]}|\ty(t)|^2\right]+8C_b^2C_{0,2}^\Upsilon\varepsilon\left[1+\be\left(\sup_{t\in[0,T]}|\tu(t)|^2\right)\right],
\end{split}
\end{equation}
separately. By means of Lemma \ref{lem6-17}, taking supremum with respect to $t$ and expectation on both sides of \eqref{eq6-71}, 
collecting  \eqref{eq6-72}-\eqref{eq6-75} and picking up $\tilde{\varepsilon}_0=\frac{1}{4(C_gT+2C_{0,2}^\Upsilon+4C_b^2C_gT+8C_b^2C_{0,2}^{\Upsilon})}$ such that, for all
$0<\varepsilon\leq\tilde{\varepsilon}_0$, one has
$ \varepsilon \left(C_gT+2C_{0,2}^\Upsilon+4C_b^2C_gT+8C_b^2C_{0,2}^{\Upsilon}\right)\leq 1/4.$
Therefore, there exists a  constant $\tilde{C}_1$, such that for any $0<\varepsilon\leq\tilde{\varepsilon}_0$,
$\be\left[\sup_{t\in[0,T]}|\ty(t)|^2\right]+2\eta\be\int_0^T\|\ty(t)\|^2dt\leq\tilde{C}_1\varepsilon.$
We finish the proof of this lemma. 
\end{proof} 

\begin{theorem}\label{thm6-21}({\bf Verifying Condition \ref{con6-1}$(b)$}) Fix $\Upsilon\in\mathbb{N}$, let $\phi_\varepsilon=(\psi_\varepsilon,\varphi_\varepsilon)$, $\phi=(\psi,\varphi)\in\tilde{\mathcal{U}}^\Upsilon$ be such that $\phi_\varepsilon$ converges in distribution to $\phi$ as $\varepsilon\rightarrow 0$. Then
$$\mathcal{G}^\varepsilon\left(\sqrt{\varepsilon}W+\int_0^{\cdot}\psi_\varepsilon(s)ds,\varepsilon N^{\varepsilon^{-1}\varphi_{\varepsilon}}\right)\Rightarrow \mathcal{G}^0\left(\int_0^{\cdot}\psi(s)ds,\lambda^\varphi_T\right).$$
\end{theorem}
\begin{proof}  Note that $\tu=\mathcal{G}^\varepsilon(\sqrt{\varepsilon}W+\int_0^{\cdot}\psi_\varepsilon(s)ds,\varepsilon N^{\varepsilon^{-1}\varphi_{\varepsilon}})$, $\varepsilon\in(0,\varepsilon_0)$. Lemmas \ref{lem6-17}, \ref{lem6-19} and \ref{lem6-20} imply
\begin{enumerate}
\item [$(i)$] $\{\tu,\varepsilon\in(0,\varepsilon_0)\}$ is tight in $L^q(0,T;\hh)\cap \mathbf{D}([0,T];D(A^{-r}))$;
\item [$(ii)$] $\lim_{\varepsilon\rightarrow 0}\be[\sup_{t\in[0,T]}|\ty(t)|^2]+\be\int_0^T\|\ty(t)\|^2dt=0$,
\end{enumerate}
where $\ty$ is the solution of  \eqref{eq6-70}. 

Set
$
\Sigma=\left( L^q(0,T;\hh)\cap\mathbf{D}([0,T];D(A^{-r}));\tilde{\mathcal{U}}^\Upsilon; L^2(0,T;\vv)\cap\mathbf{D}([0,T];\hh)\right).
$
Let $(\tilde{u},(\psi,\varphi),0)$ be any limit point of the tight family $\{(\tu,(\psi_\varepsilon,\varphi_\varepsilon),\ty),\varepsilon\in(0,\varepsilon_0)\}$. We must show that $\tilde{u}$ has the same law as $\mathcal{G}^0(\int_0^{\cdot}\psi(s)ds,\lambda^{\varphi}_T)$, and actually $\tu\Rightarrow \tilde{u}$ in the smaller space $\mathbf{D}([0,T];\hh)$.

By the Skorokhod representation theorem, there exists a probability space $(\tilde{\Omega},\tilde{\mathcal{F}},\{\tilde{\mathcal{F}}_t\}_{t\geq0},\tilde{\mathbb{P}})$ with expectation $\tilde{\e}$, $\Sigma$-valued random variables $(\tilde{u}_1,(\psi_1,\varphi_1),0)$ and $(\tu_1,(\psi_\varepsilon^1,\varphi_\varepsilon^1),\ty_1)$, $\varepsilon\in (0,\varepsilon_0)$, such that on this basis,  $(\tu_1, (\psi_\varepsilon^1,\varphi_\varepsilon^1),\ty_1)$ (respectively $(\tilde{u}_1,(\psi^1,\varphi^1),0)$) has the same law as $(\tu,(\psi_\varepsilon,\varphi_\varepsilon),\ty)$ (respectively, $(\tilde{u},(\psi,\varphi),0)$). Moreover, $(\tu_1,(\psi_\varepsilon^1,\varphi_\varepsilon^1),\ty_1)\rightarrow (\tilde{u}_1, (\psi^1,\varphi^1),0)$ in $\Sigma$, $\tilde{\mathbb{P}}$-a.s.

From the equation satisfied by $(\tu, (\psi_\varepsilon,\varphi_\varepsilon),\ty)$, we see that $(\tu_1,(\psi^1_\varepsilon,\varphi^1_\varepsilon),\ty_1)$ satisfies the following integral equation,
\begin{equation*}
\begin{split}
\tu_1(t)-\ty_1(t)&=u_0-\int_0^t(-\Delta)^\gamma(\tu_1(s)-\ty_1(s))ds-\delta\int_0^t(\tu_1(s)-\ty_1(s))ds-\int_0^tF(\tu_1(s))ds\\
&~~+\int_0^tg(s,\tu(s))\psi_\varepsilon^1(s)ds+\int_0^t\int_Eh(s,\tu(s-),\xi)(\varphi^1_\varepsilon(s,\xi)-1)\lambda(d\xi)ds,
\end{split}
\end{equation*}
and
\begin{equation*}
\begin{split}
&~\quad\tilde{\mathbb{P}}\left(\tu_1-\ty_1\in C([0,T];\hh)\cap L^2(0,T;\vv)\cap L^{p+1}(0,T;L^{p+1}(\mathbb{R}^d))\right)\\
&=\bar{\mathbb{P}}\left(\tu-\ty\in C([0,T];\hh)\cap L^2(0,T;\vv)\cap L^{p+1}(0,T;L^{p+1}(\mathbb{R}^d))\right)=1.
\end{split}
\end{equation*}
Let $\tilde{\tilde{\Omega}}$ be the subset of $\tilde{\Omega}$ such that $(\tu_1,(\psi^1_\varepsilon,\varphi^1_\varepsilon),\ty_1)\rightarrow (\tilde{u}_1,(\psi^1,\varphi^1),0)$ in $\Sigma$, then $\tilde{\mathbb{P}}(\tilde{\tilde{\Omega}})=1$. Now, we will prove that, for any fixed $\tilde{\omega}\in\tilde{\tilde{\Omega}}$,
\begin{equation}\label{eq6-76}
\sup_{t\in[0,T]}\left|\tu_1(\tilde{\omega},t)-\tilde{u}_1(\tilde{\omega},t)\right|^2\rightarrow 0,\qquad \mbox{as}~\varepsilon\rightarrow 0.
\end{equation}
Let $\tr=\tu_1-\ty_1$, then $\tr(\tilde{\omega})\in C([0,T];\hh)\cap L^2(0,T;\vv)\cap L^{p+1}(0,T;L^{p+1}(\mathbb{R}^d))$, and  satisfies 
\begin{equation*}
\begin{split}
\tr(t)&=u_0-\int_0^t(-\Delta)^\gamma\tr(s)ds-\delta\int_0^t\tr(s)ds-\int_0^tF(\tr(s)+\ty_1(s))ds\\
&~~+\int_0^tg(s,\tr(s)+\ty_1(s))\psi_\varepsilon^1(s)ds+\int_0^t\int_Eh(s,\tr(s-)+\ty_1(s-),\xi)(\varphi^1_\varepsilon(s,\xi)-1)\lambda(d\xi)ds.
\end{split}
\end{equation*}
Since
$\lim_{\varepsilon\rightarrow 0}\left[\sup_{t\in[0,T]}|\ty_1(\tilde{\omega},t)|^2+\int_0^T\|\ty_1(\tilde{\omega},t)\|^2dt\right]=0,
$
by similar arguments as in the proof of Proposition \ref{pro6-14}, we infer that
\begin{equation}\label{eq6-77}
\lim_{\varepsilon\rightarrow 0}\left[\sup_{t\in[0,T]}|\tu_1(\tilde{\omega},t)-\hat{u}(\tilde{\omega},t)|^2\right]=0,
\end{equation}
where 
\begin{equation*}
\begin{split}
\hat{u}(t)&=u_0-\int_0^t(-\Delta)^\gamma \hat{u}(s)ds-\delta\int_0^t\hat{u}(s)ds-\int_0^tF(\hat{u}(s))ds\\
&~~+\int_0^tg(s,\hat{u}(s))\psi^1(s)ds+\int_0^t\int_Eh(s,\hat{u}(s-),\xi)(\varphi^1(s,\xi)-1)\lambda(d\xi)ds.
\end{split}
\end{equation*}
Hence, $\tilde{u}_1=\hat{u}=\mathcal{G}^0(\int_0^{\cdot}\psi^1(s)ds,\lambda_T^{\varphi^1})$, and $\tilde{u}$ has the same law as $\mathcal{G}^0(\int_0^{\cdot}\psi(s)ds,\lambda_T^{\varphi})$. Since $\tu=\tu_1$ in law, \eqref{eq6-77} further implies that
$ \tu\Rightarrow \mathcal{G}^0\left(\int_0^{\cdot}\psi(s)ds,\lambda_T^{\varphi}\right).$
Thus, we complete the proof of this theorem. 
\end{proof} 

{\bf Continuation of the proof of Theorem \ref{thm6-13}}: We need to check Condition \ref{con6-1} is fulfilled. The verification of Condition \ref{con6-1}$(a)$ is given by Proposition \ref{pro6-14}, and the verification   of Condition \ref{con6-1}$(b)$ is proved in Theorem \ref{thm6-21}. The proof of the main theorem of this section is finished. $\Box$

\section{Example: Fractional stochastic Chafee-Infante equations}\label{s7}

Consider the following fractional stochastic Chafee-Infante equations driven by L\'evy noise and Brownian motion,
\begin{alignat*}{3}
\begin{cases}
\displaystyle du(t)+(-\Delta)^{\gamma}u(t)dt+\nu(u^3(t)-u(t))dt=g(u(t))dW(t)+\int_{E}h(u(t-),\xi)\tilde{N}(dt,d\xi),\\
u(t,x)=0,\\
u(0,x)=u_0(x),\\
\end{cases}
\begin{aligned}
&\mbox{in}~~\mathcal{O}\times (0,\infty),\\
&\mbox{on}~\partial\mathcal{O}\times(0,\infty),\\
&\mbox{in}~~\mathcal{O},\\
\end{aligned}
\end{alignat*}
where $\mathcal{O}$ is an open, bounded subset of $\mathbb{R}^d$
($d\leq 3$) with smooth boundary $\partial\mathcal{O}$. To put the above equation in the form of the abstract way, we only take the nonlinear term 
$ F(u)=\nu u^3-u,$
where $p+1=3$ and the conjugate number is $q=\frac{3}{2}$. In order to use the result $D(A^r)$ is continuously embedded into $L^3(\mathcal{O})$ \cite[Lemma 2.1]{W2}, we need to take $r>\frac{d}{12\gamma}$, where $\gamma$ is the index of fractional Laplacian operator. Assuming the same assumptions as the previous sections for the stochastic terms, we can straightforwardly apply our theory to this interesting example.

\section{Appendix}\label{s8}
{\it A.1. Proof of \eqref{eq6-59}}\\
\begin{proof}  Applying It\^o's formula to $|\tu|^2$, by \eqref{eq6-57}, we obtain
\begin{equation}\label{eq7-1}
\begin{split}
|\tu(t)|^2&=|u_0|^2-2\int_0^t<(-\Delta)^\gamma\tu(s)+\delta\tu(s)+F(\tu(s)),\tu(s)>ds\\
&~~+2\int_0^t(g(s,\tu(s))\psi_\varepsilon(s),\tu(s))ds+2\int_0^t\int_E(h(s,\tu(s-),\xi),\tu(s))(\varphi_\varepsilon(s,\xi)-1)\lambda(d\xi)ds\\
&~~+2\sqrt{\varepsilon}\int_0^t(\tu(s),g(s,\tu(s))dW(s))+\varepsilon\int_0^t\|g(s,\tu(s))\|^2_{\mathcal{L}_2(U;\hh)}ds\\
&~~+2\varepsilon\int_0^t\int_E(h(s,\tu(s-),\xi),\tu(s))\left(N^{\varepsilon^{-1}\varphi_{\varepsilon}}(dsd\xi)-\varepsilon^{-1}\varphi_\varepsilon(s,\xi)\lambda(d\xi)ds\right)\\
&~~+\varepsilon\int_0^t\int_E| h(s,\tu(s-),\xi)|^2\varphi_\varepsilon(s,\xi)\lambda(d\xi)ds\\[0.8ex]
&:=|u_0|^2+J_1+J_2+J_3+J_4+J_5+J_6+J_7.
\end{split}
\end{equation}
Similar to estimates \eqref{eq3-5}-\eqref{eq3-6}, we have 
\begin{equation}\label{eq7-2}
\begin{split}
J_1&:=-2\int_0^t<(-\Delta)^\gamma\tu(s)+\delta\tu(s)+F(\tu(s)),\tu(s)>ds\\
&\leq -2\eta\int_0^t\|\tu(s)\|^2ds+2k_2|\mathcal{O}|t-2k_3\int_0^t\|\tu(s)\|_{p+1}^{p+1}ds.
\end{split}
\end{equation}
By assumption $(g_2)$ and the Young inequality, we derive
\begin{equation}\label{eq7-3}
\begin{split}
J_2+J_5&:=2\int_0^t(g(s,\tu(s))\psi_\varepsilon(s),\tu(s))ds+\varepsilon\int_0^t\|g(s,\tu(s))\|^2_{\mathcal{L}_2(U;\hh)}ds\\
&\leq2\int_0^t\|g(s,\tu(s))\|_{\mathcal{L}_2(U;\hh)}\|\psi_\varepsilon(s)\|_U|\tu(s)|ds+\varepsilon\int_0^t\|g(s,\tu(s))\|^2_{\mathcal{L}_2(U;\hh)}ds\\
&\leq\frac{1}{2}\int_0^t\|g(s,\tu(s))\|^2_{\mathcal{L}_2(U;\hh)}\|\psi_\varepsilon(s)\|^2_Uds+2\int_0^t|\tu(s)|^2ds+\varepsilon\int_0^t\|g(s,\tu(s))\|^2_{\mathcal{L}_2(U;\hh)}ds\\
&\leq \frac{1}{2}C_g\int_0^t(1+|\tu(s)|^2)\|\psi_\varepsilon(s)\|_U^2ds+2\int_0^t|\tu(s)|ds+\varepsilon C_g\int_0^t(1+|\tu(s)|^2)ds\\
&\leq \frac{1}{2}C_g\int_0^t\|\psi_\varepsilon(s)\|_U^2ds+\varepsilon C_gt+\int_0^t\left(\frac{1}{2}C_g\|\psi_\varepsilon(s)\|_U^2+2+\varepsilon C_g\right)|\tu(s)|^2ds.
\end{split}
\end{equation}
By  Lemma \ref{lem6-7}, we infer that
\begin{equation}\label{eq7-4}
\begin{split}
J_3&:=2\int_0^t\int_E(h(s,\tu(s-),\xi),\tu(s))(\varphi_\varepsilon(s,\xi)-1)\lambda(d\xi)ds\\
&\leq 2\int_0^t\int_E|h(s,\tu(s-),\xi)||\tu(s)||\varphi_\varepsilon(s,\xi)-1|\lambda(d\xi)ds\\
&\leq 2\int_0^t\int_E\|h(s,\xi)\|_{0,\hh}(1+|\tu(s)|)|\tu(s)||\varphi_\varepsilon(s,\xi)-1|\lambda(d\xi)ds\\
&\leq 2\int_0^t(1+2|\tu(s)|^2)\left(\int_E\|h(s,\xi)\|_{0,\hh}|\varphi_\varepsilon(s,\xi)-1|\lambda(d\xi)\right)ds\\
&\leq 2C_{0,1}^{\Upsilon}+4\int_0^t|\tu(s)|^2\left(\int_E\|h(s,\xi)\|_{0,\hh}|\varphi_\varepsilon(s,\xi)-1|\lambda(d\xi)\right)ds.
\end{split}
\end{equation}
Again making use of Lemma \ref{lem6-7}, $J_7$ can be bounded by
\begin{equation}\label{eq85}
\begin{split}
J_7&=\varepsilon\int_0^t\int_E|h(s,\tu(s-),\xi)|^2\varphi_\varepsilon(s,\xi)\lambda(d\xi)ds\\
&\leq \varepsilon\int_0^t\int_E\|h(s,\xi)\|^2_{0,\hh}\left(1+|\tu(s)|\right)^2\varphi_\varepsilon(s,\xi)\lambda(d\xi)ds\\
&\leq 2\varepsilon C_{0,2}^\Upsilon+2\varepsilon \int_0^t\int_E \|h(s,\xi)\|^2_{0,\hh}
|\tu(s)|^2\varphi_\varepsilon(s,\xi)\lambda(d\xi)ds.
\end{split}
\end{equation}
For $J_4$, it follows from assumption $(g_2)$, the Burkholder-Davis-Gundy and Young inequalities that
\begin{equation}\label{eq7-9}
\begin{split}
\be\left(\sup_{t\in[0,T]}|J_4(t)|\right)&\leq
2C_b\sqrt{\varepsilon}\be\left(\int_0^T\|g(t,\tu(t))\|^2_{\mathcal{L}_2(U;\hh)}|\tu(t)|^2dt\right)^{\frac{1}{2}}\\
&\leq \frac{1}{4}\sqrt{\varepsilon}\be\left[\sup_{t\in[0,T]}|\tu(t)|^2\right]+4C_b^2\sqrt{\varepsilon}\be\int_0^T\|g(t,\tu(t))\|^2_{\mathcal{L}_2(U;\hh)}dt\\
&:=\left(\frac{1}{4}\sqrt{\varepsilon}+4C_b^2C_g\sqrt{\varepsilon} T\right)\be\left[\sup_{t\in[0,T]}|\tu(t)|^2\right]+4C_b^2C_g\sqrt{\varepsilon} T.
\end{split}
\end{equation}
For $J_6$, by Lemma \ref{lem6-7}, the Burkholder-Davis-Gundy and Young inequalities, we also obtain
\begin{equation}\label{eq7-10}
\begin{split}
\be\left(\sup_{t\in[0,T]}|J_6(t)|\right)&\leq 2C_b\be\left(\int_0^T\int_E\varepsilon^2|h(t,\tu(t-),\xi)|^2|\tu(t)|^2\varepsilon^{-1}\varphi_\varepsilon(t,\xi)\lambda(d\xi)dt\right)^{\frac{1}{2}}\\
&\leq \frac{1}{4}\be\left[\sup_{t\in[0,T]}|\tu(t)|^2\right]+4C_b^2\be\int_0^T\int_E\varepsilon|h(t,\tu(t-),\xi)|^2|\varphi_\varepsilon(t,\xi)|\lambda(d\xi)dt\\
&\leq \frac{1}{4}\be\left[\sup_{t\in[0,T]}|\tu(t)|^2\right]+4C_b^2\varepsilon\be\int_0^T\int_E\|h(t,\xi)\|^2_{0,\hh}(1+|\tu(t)|)^2|\varphi_\varepsilon(t,\xi)|\lambda(d\xi)dt\\
&\leq\frac{1}{4}\be\left[\sup_{t\in[0,T]}|\tu(t)|^2\right] +8\varepsilon C_b^2C_{0,2}^\Upsilon\be\left[1+\sup_{t\in[0,T]}|\tu(t)|^2\right].
\end{split}
\end{equation}
Combining \eqref{eq7-1}-\eqref{eq85}, for all $t\in[0,T]$, we arrive at
\begin{equation*}
\begin{split}
|\tu(t)|^2&\leq |u_0|^2+2\kappa_2|\mathcal{O}|t+C_g\Upsilon+\varepsilon C_gt+\int_0^t\left(\frac{1}{2}C_g\|\psi_\varepsilon(s)\|_U^2+2+\varepsilon C_g\right)|\tu(s)|^2ds\\
&~~+2C_{0,1}^{\Upsilon}+4\int_0^t|\tu(s)|^2\left(\int_E\|h(s,\xi)\|_{0,\hh}|\varphi_\varepsilon(s,\xi)-1|\lambda(d\xi)\right)ds+\sup_{t\in[0,T]}|J_4(t)|\\
&~~+2\varepsilon C_{0,2}^\Upsilon+2\varepsilon \int_0^t\int_E \|h(s,\xi)\|^2_{0,\hh}
|\tu(s)|^2\varphi_\varepsilon(s,\xi)\lambda(d\xi)ds+\sup_{t\in[0,T]}|J_6(t)|\\
&\leq  \left(|u_0|^2+2\kappa_2|\mathcal{O}|T+C_g\Upsilon+\varepsilon C_gT+2C_{0,1}^{\Upsilon}+2\varepsilon C_{0,2}^\Upsilon+\sup_{t\in[0,T]}|J_4(t)|+\sup_{t\in[0,T]}|J_6(t)|\right)\\
&~~+\int_0^t\bigg[\frac{1}{2}C_g\|\psi_\varepsilon(s)\|_U^2+2+\varepsilon C_g+\bigg(
\int_E\left(4\|h(s,\xi)\|_{0,\hh}|\varphi_\varepsilon(s,\xi)-1| \right.\\[0.8ex]
&~~\left.+2\varepsilon\|h(s,\xi)\|^2_{0,\hh}\varphi_\varepsilon(s,\xi)\right)\lambda(d\xi)\bigg)\bigg]|\tu(s)|^2ds,\qquad \bar{\mathbb{P}}\mbox{-a.s.}
\end{split}
\end{equation*}
Let $M_1=2\kappa_2|\mathcal{O}|T+C_g\Upsilon+\varepsilon C_gT+2C_{0,1}^{\Upsilon}+2\varepsilon C_{0,2}^\Upsilon$. Using the Gronwall Lemma,  we have
\begin{equation*}
\begin{split}
|\tu(t)|^2&\leq \left(|u_0|^2+M_1+\sup_{t\in[0,T]}|J_4(t)|+\sup_{t\in[0,T]}|J_6(t)|\right)
 \exp\left(C_g\Upsilon+2T+\varepsilon C_gT+4C_{0,1}^{\Upsilon}+2\varepsilon C_{0,2}^\Upsilon\right).
\end{split}
\end{equation*}
Denote by $M_2=\exp\left(C_g\Upsilon+2T+\varepsilon C_gT+4C_{0,1}^{\Upsilon}+2\varepsilon C_{0,2}^\Upsilon\right)$ which does not depend on $\omega$. It follows from \eqref{eq7-9}-\eqref{eq7-10} that
\begin{equation*}
\begin{split}
\be\left[\sup_{t\in[0,T]}|\tu(t)|^2\right]&\leq M_2\left(\be|u_0|^2+M_1+4C_b^2C_g\sqrt{\varepsilon} T+8\varepsilon C_b^2C_{0,2}^\Upsilon\right)\\
&~~+M_2\left(\frac{1}{4}\sqrt{\varepsilon}+4C_b^2C_g\sqrt{\varepsilon} T+\frac{1}{4}+
8\varepsilon C_b^2C_{0,2}^\Upsilon\right)\be\left[\sup_{t\in[0,T]}|\tu(t)|^2\right].
\end{split}
\end{equation*}
We choose $\varepsilon_0\leq \frac{1}{(1+16C_b^2C_g T+32C_b^2C_{0,2}^\Upsilon)^2}$  such that
$\frac{1}{4}\sqrt{\varepsilon}+4C_b^2C_g\sqrt{\varepsilon} T+\frac{1}{4}+
8\varepsilon C_b^2C_{0,2}^\Upsilon\leq \frac{1}{2}.$
Therefore, 
\begin{equation*}
\be\left[\sup_{t\in[0,T]}|\tu(t)|^2\right]\leq 2M_2\left(\be|u_0|^2+M_1+4C_b^2C_g\sqrt{\varepsilon} T+8\varepsilon C_b^2C_{0,2}^\Upsilon\right).
\end{equation*}
The proof is complete. 
\end{proof} 

\bigskip

{\it A.2.  Proof of Lemma \ref{lem6-19}: $\langle\tu(\beta_\varepsilon+d_\varepsilon)-\tu(\beta_\varepsilon),l\rangle_{D(A^r)}\rightarrow 0$ in probability as $\varepsilon \rightarrow 0$ for every $l\in D(A^r)$, where $(\beta_\varepsilon,d_\varepsilon)$ are a stopping time with respect to the natural $\bar{\sigma}$-field taking only finitely many values, and an interval on $[0,T]$, respectively, satisfying  $d_\varepsilon\rightarrow 0$ as $\varepsilon\rightarrow 0$. }

\begin{proof}  With a slight abuse of notation, we will use the inner product $(\cdot,\cdot)$ instead of $\langle\cdot,\cdot\rangle_{D(A^r)}$. For simplicity, denoted by $\bar{d}:=d_\varepsilon$ and $\beta:=\beta_\varepsilon$. By \eqref{eq6-57}, we have 
\begin{equation*}
\begin{split}
\tu(\beta+\bar{d})-\tu(\beta)&=-\left(\int_{\beta}^{\bar{d}+\beta}(-\Delta)^{\gamma}\tu(s)ds+\delta\int_{\beta}^{\bar{d}+\beta}\tu(s)ds\right)-\int_{\beta}^{\bar{d}+\beta}F(\tu(s))ds\\
&~~+\int_{\beta}^{\bar{d}+\beta}g(s,\tu(s))\psi_\varepsilon(s)ds+\sqrt{\varepsilon}\int_{\beta}^{\bar{d}+\beta}g(s,\tu(s))dW(s)\\
&~~+\int_\beta^{\bar{d}+\beta}\int_Eh(s,\tu(s-),\xi)(\varphi_\varepsilon(s,\xi)-1)\lambda(d\xi)ds\\
&~~+\varepsilon\int_\beta^{\bar{d}+\beta}\int_Eh(s,\tu(s-),\xi)\left(N^{\varepsilon^{-1}\varphi_{\varepsilon}}(dsd\xi)-\varepsilon^{-1}\varphi_\varepsilon(s,\xi)\lambda(d\xi)ds\right)\\
&:=I_1^\varepsilon+I_2^\varepsilon+I_3^\varepsilon+I_4^\varepsilon+I_5^\varepsilon+I_6^\varepsilon.
\end{split}
\end{equation*}
For $I_1^\varepsilon$, since $A\tu\in L^2(0,T;\vv^*)$ and $l\in D(A^r)$, by the H\"older inequality and \eqref{eq6-59}, we have
 \begin{equation*}
 \begin{split}
 \lim_{\varepsilon\rightarrow 0}\be\left|\int_\beta^{\bar{d}+\beta}(A\tu(s),l)ds\right|&\leq \lim_{\varepsilon\rightarrow 0}\be\int_\beta^{\bar{d}+\beta}\|A\tu(s)\|_*\|l\|ds
 \leq  \lim_{\varepsilon\rightarrow 0}C\|l\| \sqrt{\bar{d}}=0.
 \end{split}
 \end{equation*}
For $I_2^\varepsilon$, since $F\in L^q(0,T;L^q(\mathcal{O}))$, combining with \eqref{eq6-59} and the H\"older inequality, we obtain
\begin{equation*}
\begin{split}
 \lim_{\varepsilon\rightarrow 0}\be \left|\int_\beta^{\bar{d}+\beta} F(\tu(s)),l)ds\right| &\leq  \lim_{\varepsilon\rightarrow 0}\be \int_\beta^{\bar{d}+\beta}\int_\mathcal{O}|F(\tu(s))||l|dxds\\
 & \leq  \lim_{\varepsilon\rightarrow 0}\be
|l|\bar{d}^{\frac{1}{p+1}}\left(\int_\beta^{\bar{d}+\beta}\|F(\tu(s))\|_qds\right)^{\frac{1}{q}}=0.
\end{split}
\end{equation*}
For $I_3^\varepsilon$, by condition $(g_2)$, the H\"older inequality and \eqref{eq6-59}, we infer
\begin{equation*}
\begin{split}
&\quad~ \lim_{\varepsilon\rightarrow 0}\be\left|\int_\beta^{\bar{d}+\beta}(g(s,\tu(s))\psi_\varepsilon(s),l)ds\right|\leq \lim_{\varepsilon\rightarrow 0}\be |l|
\int_\beta^{\bar{d}+\beta}\|g(s,\tu(s))\|_{\mathcal{L}_2(U;\hh)}\|
\psi_\varepsilon(s)\|_Uds\\
&\leq \lim_{\varepsilon\rightarrow 0}\be |l|\sqrt{C_g}\left(\int_\beta^{\bar{d}+\beta}\left(1+|\tu(s)|^2\right) ds\right)^{\frac{1}{2}}\left(\int_\beta^{\bar{d}+\beta}||\psi_\varepsilon(s)\|_U^2ds\right)^{\frac{1}{2}}\\[1.0ex]
&\leq  \lim_{\varepsilon\rightarrow 0} |l|\sqrt{2C_g}\sqrt{\bar{d}}\sqrt{\Upsilon}\be\sqrt{1+\sup_{t\in[0,T]}|\tu(t)|^2}=0.
\end{split}
\end{equation*}
For $I_5^\varepsilon$,  by Lemma \ref{lem6-7}$(ii)$, together with \eqref{eq6-59}, we derive 
\begin{equation*}
\begin{split}
&~\quad\lim_{\varepsilon\rightarrow 0}\be\left|\int_\beta^{\bar{d}+\beta}\int_E(h(s,\tu(s-),\xi),l)(\varphi_\varepsilon(s,\xi)-1)\lambda(d\xi)ds\right|\\
&\leq \lim_{\varepsilon\rightarrow 0}\be|l|\int_\beta^{\bar{d}+\beta}\int_E\|h(s,\xi)\|_{0,\hh}(1+|\tu(s)|)
|\varphi_\varepsilon(s,\xi)-1|\lambda(d\xi)ds\\
&\leq \lim_{\varepsilon\rightarrow 0}\be |l|\left[\left(1+\sup_{t\in[0,T]}|\tu(t)|\right)\left(\int_\beta^{\bar{d}+\beta}\int_E\|h(s,\xi)\|_{0,\hh}
|\varphi_\varepsilon(s,\xi)-1|\lambda(d\xi)ds\right)\right]=0.
\end{split}
\end{equation*}
Moreover, for $I_4^\varepsilon$ and $I_6^\varepsilon$, by the Burkholder-Davis-Gundy and Young inequalities,  Lemma \ref{lem6-7}, condition $(g_2)$ and \eqref{eq6-59}, we arrive at
$
\lim_{\varepsilon\rightarrow 0}\be |(I_4^\varepsilon, l)|=0
\ \mbox{and}\ \lim_{\varepsilon\rightarrow 0}\be |(I_6^\varepsilon, l)|=0$, respectively.
Therefore, collecting all the estimates above, we conclude the proof. 
\end{proof}

\section*{Acknowledgement}
The research  has been partially supported by the Spanish Ministerio de Ciencia e Innovaci\'on, Agencia Estatal de Investigaci\'on (AEI) and Fondo Europeo de Desarrollo Regional (FEDER) under the projects PGC2018-096540-I00 and PID2021-122991NB-C21, by Junta de Andaluc\'{\i}a and FEDER under the project P18-FR-4509, by the
Generalitat Valenciana, project PROMETEO/2021/063, and the Nature Science
Foundation of Jiangsu Province (Grant No. BK20220233).




\end{document}